\documentclass[10pt,letterpaper]{amsart}
\usepackage{charter}
\usepackage[OT2,OT1]{fontenc}
\usepackage[latin9]{luainputenc}
\usepackage{color}
\usepackage{float}
\usepackage{textcomp}
\usepackage{bm}
\usepackage{amstext}
\usepackage{amsthm}
\usepackage{amssymb}
\usepackage{cancel}
\usepackage{stmaryrd}
\usepackage[all]{xy}

\makeatletter

\providecommand{\tabularnewline}{\\}

\numberwithin{equation}{section}
\numberwithin{figure}{section}
\theoremstyle{plain}
\newtheorem{thm}{\protect\theoremname}[section]
\theoremstyle{plain}
\newtheorem{lem}[thm]{\protect\lemmaname}
\theoremstyle{definition}
\newtheorem{defn}[thm]{\protect\definitionname}
\theoremstyle{plain}
\newtheorem{prop}[thm]{\protect\propositionname}
\theoremstyle{plain}
\newtheorem{cor}[thm]{\protect\corollaryname}
\theoremstyle{remark}
\newtheorem{rem}[thm]{\protect\remarkname}

\usepackage{amsfonts}
\usepackage{mathrsfs}
\usepackage{graphicx}
\usepackage{amscd}

\newcommand{\cyr}{
\renewcommand\rmdefault{wncyr} \renewcommand\sfdefault{wncyss} \renewcommand\encodingdefault{OT2} \normalfont
\selectfont
}
\DeclareTextFontCommand{\textcyr}{\cyr}    

   \def\@settitle
{
\begin{center} \baselineskip14\p@\relax \LARGE\@title \end{center}
}

\usepackage[normalem]{ulem}
\usepackage[mathscr]{eucal}
\numberwithin{equation}{section}

\usepackage{verbatim}
\usepackage{amscd}
\usepackage{bm}\usepackage[mathscr]{eucal}
\usepackage[all]{xy}
\usepackage{enumerate}
\usepackage{color}
\usepackage{colortbl}

\input xyimport.tex

\font\eu=eusm10 at 10pt

\def\eF{\text{\eu F}}

\title{On Classification of $\mQ$-Fano 3-folds \\ of
Gorenstein index 2. III}

\author{Hiromichi Takagi}

\address{Department of Mathematics, Gakushuin University, 
Mejiro, Toshima-ku, Tokyo 171-8588, Japan}
\email{hiromici@math.gakushuin.ac.jp}

\newcommand{\sA}{\mathcal{A}}
\newcommand{\sB}{\mathcal{B}}

\newcommand{\sD}{\mathcal{D}}
\newcommand{\sE}{\mathcal{E}}
\newcommand{\sF}{\eF}

\newcommand{\sL}{\mathcal{L}}
\newcommand{\sN}{\mathcal{N}}

\newcommand{\sO}{\mathcal{O}}

\newcommand{\sR}{\mathcal{R}}
\newcommand{\sQ}{\mathcal{Q}}

\newcommand{\sU}{\mathcal{U}}
\newcommand{\sV}{\mathcal{V}}

\newcommand{\mC}{\mathbb{C}}

\newcommand{\rG}{\mathrm{G}}

\newcommand{\mP}{\mathbb{P}}
\newcommand{\mQ}{\mathbb{Q}}

\newcommand{\Bs}{\mathrm{Bs}\,}

\newcommand{\Pic}{\mathrm{Pic}\,}

\newcommand{\Sing}{\mathrm{Sing}\,}

\newcommand{\SL}{\mathrm{SL}\,}

\newcommand{\rank}{\mathrm{rank}\,}

\numberwithin{equation}{section}

 \newcounter{myparagraph}[subsection]

\makeatother

\providecommand{\corollaryname}{Corollary}
\providecommand{\definitionname}{Definition}
\providecommand{\lemmaname}{Lemma}
\providecommand{\propositionname}{Proposition}
\providecommand{\remarkname}{Remark}
\providecommand{\theoremname}{Theorem}

\begin{document}
\maketitle 

\begin{center}

\textit{\small{}Dedicated to Professor Yujiro Kawamata on the occasion
of his 70th birthday}{\small\par}

\end{center}
\begin{abstract}
We classified prime $\mQ$-Fano $3$-folds $X$ with only $1/2(1,1,1)$-singularities
and with $h^{0}(-K_{X})\geq4$ a long time ago. The classification
was undertaken by blowing up each $X$ at one $1/2(1,1,1)$-singularity
and constructing a Sarkisov link. The purpose of this paper is to
reveal the geometries behind the Sarkisov links for $X$ in 5 classes.
The main result asserts that any $X$ in the 5 classes can be embedded
as linear sections into bigger dimensional $\mQ$-Fano varieties called
key varieties, where the key varieties are constructed by extending
partially the Sarkisov link in higher dimensions. 
\end{abstract}

\maketitle
\markboth{$\mQ$-Fano $3$-folds III}{Hiromichi Takagi} {\small{}{\tableofcontents{}}}{\small\par}

2020\textit{ Mathematics subject classification}: 14J45, 14E05 .


\section{\textbf{Introduction}}

\subsection{Background}

In this paper, we work over $\mC$, the complex number field.

This is a continuation of the papers \cite{Tak1} after a long time.
A projective variety $X$ is called a \textit{$\mQ$-Fano variety}
if $X$ has only terminal singularities and $-K_{X}$ is ample. A
$\mQ$-Fano variety $X$ is called \textit{prime} if $-K_{X}$ generates
the group of numerical equivalence classes of $\mQ$-Cartier divisors
on $X$. In \cite{Tak1}, we classified prime $\mQ$-Fano 3-folds
$X$ with only $1/2(1,1,1)$-singularities and with $h^{0}(-K_{X})\geq4$.
In this paper, we further study $X$ in the 5 classes No.1.1,~1.4,~1.9,~1.10,
and 1.13 among \cite[Table 1]{Tak1}.

\subsection{Prime $\mQ$-Fano 3-fold and Sarkisov link}

In this subsection, we explain our method of the classification of
prime $\mQ$-Fano 3-folds in \cite{Tak1} only in the five classes.
The result is presented in the following table:

\begin{table}[H]
\begin{tabular}{|c|c|c|c|c|c|c|}
\hline 
No. & $g(X)$ & $N$ & $e$ & $\deg C$ & $g(C)$ & $X'$\tabularnewline
\hline 
\hline 
1.1 & $4$ & $2$ & $7$ & $7$ & $8$ & $\mP(1^{3},2)$\tabularnewline
\hline 
\hline 
1.4 & $5$ & $1$ & $6$ & $9$ & $9$ & $\mP^{3}$\tabularnewline
\hline 
1.9 & $6$ & $1$ & $6$ & $3$ & $0$ & $B_{3}$\tabularnewline
\hline 
1.10 & $6$ & $1$ & $5$ & $9$ & $6$ & $Q^{3}$\tabularnewline
\hline 
1.13 & $8$ & $1$ & $4$ & $7$ & $2$ & $B_{5}$\tabularnewline
\hline 
\end{tabular}\label{tab:Table 1}
\end{table}
\noindent The number $g(X)$ in the second column of the table is
the \textit{genus} of $X$ defined to be $h^{0}(-K_{X})-2.$ The number
$N$ in the third column is the number of $1/2(1,1,1)$-singularities
of $X$. We explain the data in 4th--7th column below. For $X$'s
in the 5 classes, we classify them by constructing the following Sarkisov
links:

\begin{equation}\label{eq:Sarkisov} \xymatrix{& Y\ar@{-->}[r]\ar[dl]_f & Y'\ar[dr]^{f'}\\
X & & & X',}
\end{equation}where $f\colon Y\to X$ is the blow-up of $X$ at a $1/2(1,1,1)$-singularity,
$Y\dashrightarrow Y'$ is a flop, and $f'$ is the blow-up of a $\mQ$-Fano
$3$-fold $X'$ along a smooth curve $C$ in $X'\setminus\Sing X'$
with the genus $g(C)$ and the degree $\deg C$ as in the 5th and
6th column of the table, where the degree of $C$ is measured by the
primitive \textsl{Cartier} divisor on $X'$. We denote by $E$ the
$f$-exceptional divisor and by $\widetilde{E}$ the strict transform
on $Y'$ of $E$. The flop $Y\dashrightarrow Y'$ is the $E$-flop
in the sense of \cite{Ko}. The number $e$ in the 4th column is defined
to be $E^{3}-\widetilde{E}^{3}$, which roughly measures how many
flopping curves the flop $Y\dashrightarrow Y'$ has (see \cite{Tak1}
for more details).\textit{ }In the 7th column, $B_{3}$ is a smooth
cubic 3-fold in $\mP^{4}$, $B_{5}$ is a codimension $3$ smooth
linear section of ${\rm G}(2,5)$, and $Q^{3}$ is a smooth quadric
3-fold. Let $h\colon Y\to W$ be the $E$-flopping contraction. We
call $W$ the \textit{mid point} (for No. 1.1, however, see Caution
in the subsection \ref{subsec:Genus-4mid}). 

For simplicity, we call a $\mQ$-Fano 3-fold $X$ as in the table
as follows:\vspace{3pt}

\noindent\textbf{ Names of $\mQ$-Fano 3-folds:}\textit{ a prime
$\mQ$-Fano 3-fold of genus 4, of genus 5, of genus 6 and }$\mathtt{C}$-type,\textit{
of genus 6 and }$\mathtt{Q}$-type,\textit{ of genus 8, }for $X$
of No.1.1,~1.4,~1.9,~1.10, and 1.13 respectively, where the word
$\mathtt{C}$-type (resp. $\mathtt{Q}$-type) comes from the fact
that $X'$ is isomorphic to the cubic 3-fold for No.1.9 (resp. the
quadric 3-fold for No.1.10).

\subsection{Main result}

The main result of this paper is a classification of $\mQ$-Fano 3-folds
in the 5 classes in different nature to that in \cite{Tak1}. The
prototype of this result is the following: in \cite{Gu}, Gushel\!\'\empty~
shows that any smooth prime Fano 3-fold of genus 8 is a linear section
of $\rG(2,6)$. In \cite{Mu1,Mu2}, Mukai shows that any prime smooth
prime Fano 3-fold of genus 7, 9, 10 is a linear section of the orthogonal
Grassmanian ${\rm OG}(5,10)$, the symplectic Grassmannian ${\rm Sp}(3,6)$,
and the adjoint homogeneous variety of type $G_{2}$, respectively.
Here we say that a projective variety $X$ is a \textit{linear section}
of a projective variety $\Sigma$ with respect to a linear system
$|M|$ if it holds that $X=\Sigma\cap D_{1}\cap\dots\cap D_{k}$ for
$k=\dim\Sigma-\dim X$ and some $D_{1},\dots,D_{k}\in|M|$. In the
case of Gushel\!\'\empty~ and Mukai, the linear system is the one
generated by the primitive very ample divisor of a homogeneous space.
We usually do not mention the linear system $|M|$ if $M$ generates
the group of the numerical equivalence classes of $\mQ$-Cartier divisors
on $\Sigma$.
\begin{thm}[Embedding theorem]
 \label{thm:main1} For each one of the $5$ classes, there is a
unique rational $\mQ$-Fano variety $\Sigma$ of Picard number $1$
such that any prime $\mQ$-Fano $3$-fold $X$ in the class is a linear
section of $\Sigma$. The $\mQ$-Fano varieties $\Sigma$ are of $11$-,
$12$-, $9$-, $8$-, and $5$-dimensional for $X$ of genus $4$,
$5$, of genus $6$ and $\mathtt{Q}$-type, of genus $6$ and $\text{\ensuremath{\mathtt{C}}}$-type,
and of genus $8$, respectively. 
\end{thm}

For a prime $\mQ$-Fano 3-fold $X$ in each of the 5 classes, we will
call the variety $\Sigma$ \textit{the key variety} for $X$. Theorem
\ref{thm:main1} is proved separately in each case; Theorem \ref{thm:emb8}
(genus 8), Theorem \ref{thm:embg4g6} (genus 4, 6), and Theorem \ref{thm:embg5}
(genus 5). We refer for more detailed descriptions of $\Sigma$ (constructions,
birational geometries, singularities, etc) to the section \ref{sec:Genus8}
(genus 8), the sections \ref{sec:Extending Mid G46} and \ref{sec:Common-prescription-for}
(genus 4, 6), and the section \ref{sec:Embedding-theoremG5} (genus
5). 

The equations of the key varieties are also available; see \cite[Ex.6.8]{R}
in the genus 4 case, and \cite{Ha} in the genus 5 case. In the genus
6 and 8 cases, we will publish them in separated papers (cf. \cite{Tak4}). 

\subsection{Structure of the paper}

\noindent The section \ref{sec:Preliminaries}: After showing some
miscellaneous results in the subsection \ref{subsec:Miscellaneous-results},
we investigate in the subsections \ref{subsec:IndecMid}--\ref{subsec:Mid-point-in G8}
the mid point of the Sarkisov link (\ref{eq:Sarkisov}) in details
in each of 5 cases. 

\noindent The section \ref{sec:Genus8}: In this section we concentrate
in studying the genus 8 case. In the subsection \ref{subsec:Extending-the-mid G8},
we extend the mid point based upon the result in the subsection \ref{subsec:Mid-point-in G8}.
In the subsection \ref{subsec:Construction-of-the key G8}, we construct
the key variety modifying birationally the extension of the mid point
(Theorem \ref{thm:key8}). In the subsection \ref{subsec:Embedding-theorem G8},
we show Theorem \ref{thm:main1} in the genus 8 case (Theorem \ref{thm:emb8}).

In the other cases, constructions of the key varieties and proofs
of Theorem \ref{thm:main1} are similar to those in the genus 8 case
but are more involved. 

The genus 4 and 6 cases are treated in a unified way in part in the
sections \ref{sec:Extending Mid G46} and \ref{sec:Common-prescription-for};
The mid points are extended in Propositions \ref{prop:ExtMidG4},
\ref{prop:ExtMidG6Q} and \ref{prop:ExtMidG6C}. Together with compensations
in the subsection \ref{subsec:Compensation}, the key varieties are
constructed in Theorem \ref{thm:common}, and Theorem \ref{thm:main1}
is proved in Theorem \ref{thm:embg4g6}. 

The genus 5 case is treated in the section \ref{sec:Embedding-theoremG5};
The mid point is extended in the subsection \ref{subsec:Extending-the-midG5},
the key variety is constructed in the subsection \ref{subsec:ConstructionKey G5},
and Theorem \ref{thm:main1} is proved in the subsection \ref{subsec:Embedding-theorem G5}.

\subsection{Flow of the construction of the key variety \label{subsec:Flow}}

Roughly speaking, the key variety in any case is constructed in the
following manner: First we extend the mid point of the Sarkisov link
(\ref{eq:Sarkisov}) to an appropriate variety $\overline{\Sigma}$
(in the genus 4 case, the mid point is replaced by another 3-fold
in the subsection \ref{subsec:Genus-4mid}). The extension $\overline{\Sigma}$
is found more or less naturally from the equation of the mid point.
Second we construct a good resolution of $\overline{\Sigma}$. Except
in the case of genus 6 and $\mathtt{C}$-type, we can construct a
small crepant resolution $\Sigma'\to\overline{\Sigma}$ with a projective
bundle structure over certain Fano manifold. Except in the case of
genus 4, the Fano manifold is $X'$ as in (\ref{eq:Sarkisov}). The
small resolution is actually an extension of $Y'$ in (\ref{eq:Sarkisov}).
Ideally, as the third step, we would construct a small birational
map $\Sigma'\dashrightarrow\widetilde{\Sigma}$ such that $\widetilde{\Sigma}$
is an extension of $Y$ and find a contraction $\widetilde{\Sigma}\to\Sigma$
such that $\Sigma$ is the desired key variety. This strategy works
in the genus 8 case. Even in the other cases except the case of genus
6 and $\mathtt{C}$-type, this works but the construction of $\Sigma'\dashrightarrow\widetilde{\Sigma}$
is slightly involved (we refer a more detailed explanation of this
to Remark \ref{rem:flopdiff}). Hence we choose another resolution
$\widehat{\Sigma}\to\overline{\Sigma}$ except in the genus 8 case.
The variety $\widehat{\Sigma}$ has a structure of a projective bundle
over certain Fanifold. Then we perform a small birational map $\Sigma'\dashrightarrow\widetilde{\Sigma}$
which is a composite of a flop and a flip such that $\widetilde{\Sigma}$
is an extension of $Y$ and find a contraction $\widetilde{\Sigma}\to\Sigma$
such that $\Sigma$ is the desired key variety. the advantage of this
construction is that the flop and the flip can be described easily.
Moreover, this works also in the genus 6 and $\mathtt{C}$-type. Actually,
the cases except the genus 5 case can be treated in a unified way
in part. For unified treatement, we take a bit roundabout way. We
refer for this to the sections \ref{sec:Extending Mid G46} and \ref{sec:Common-prescription-for}.
The case of genus 5 can be treated more or less in a straightforward
way. We refer for this to the section \ref{sec:Embedding-theoremG5}.

\subsection{Future plan}

In this subsection, we use the notation as in the subsection \ref{subsec:Flow}.
The construction of the key variety in each case is slightly involved
but we have a significant application; in the forthcoming paper \cite{Tak3},
we construct a projective bundle which, in a certain sense, is dual
to $\widehat{\Sigma}$ in the case of genus 6 and $\mathtt{C}$-type,
or is dual to $\Sigma'$ in the other cases. Using these dual varities,
we can describe the cubic 3-fold $X'$ in the case of genus 6 and
$\mathtt{C}$-type, or the curve $C$ as in the Sarkisov link (\ref{eq:Sarkisov})
in the other cases.

\vspace{5pt}

\noindent\textbf{ Notation and Conventions} 
\begin{itemize}
\item \textit{Conventions on projective bundle}: Let $\sE$ be a vector
bundle on a variety $X$, or a vector space. The notation $\mP(\sE)$
is just the projectivization of $\sE$ (\textit{We don't use the Grothendieck
notation}). Setting $\Sigma=\mP(\sE)$, we often denote by $\sO_{\Sigma}(1)$,
or $H_{\Sigma}$ \textit{the tautological line bundle associated to
the vector bundle $\sE$} without mentioning $\sE$.
\item \textit{Point of a projective space}: Let $V$ be a vector space.
For a nonzero vector $\bm{x}\in V$ and a 1-dimensional subspace $V^{1}\subset V$,
we denote by $[\bm{x}]$ and $[V^{1}]$ the point of $\mP(V)$ corresponding
to $\bm{x}$ and $V^{1}$ respectively.
\item \textit{Cartier divisor and invertible sheaf}: We sometimes abuse
notation of a Cartier divisor and an invertible sheaf. For example,
we sometimes use the expression like $D=f^{*}\sO_{X}(1)$.
\item \textit{Join $X_{a}*X_{b}$} : The projective variety in $\mP(V_{a}\oplus V_{b}$)
which is the union of all the lines joining two projective varieties
$X_{a}\subset\mP(V_{a}\oplus0)$ and $X_{b}\subset\mP(0\oplus V_{b})$.
If $X_{b}$ is a projective space, $X_{a}*X_{b}$ is just the cone
over $X_{a}$ with the vertex $X_{b}$.
\item \textit{Flopping contraction of Atiyah type}: The small contraction
$f\colon X\to Y$ such that, for a sufficiently small analytic neighborhood
$U$ of any point $y$ of $Y$ in the image of the $f$-exceptional
locus, $f|_{f^{-1}(U)}\colon f^{-1}(U)\to U$ is isomorphic to the
product of some open subset $U_{0}\subset\mC^{n-3}$ with a small
resolution $X_{0}\to\left\{ xy+zw=0\right\} \subset\mC^{4}$. 
\item \textit{$1/2\,(1^{n})$-singularity}: The singularity of an $n$-dimensional
variety analytically isomorphic to that of the origin of the quotient
of $\mC^{n}$ by the involution defined by $\mC^{n}\ni\bm{x}\mapsto-\bm{x}\in\mC^{n}$.
We often call this singularity a \textit{1/2-singularity} for simplicity.
\item $Q^{n}$: The \textit{smooth} quadric $n$-fold.
\end{itemize}
\textbf{\noindent Acknowledgment}: I am grateful to Professor Shinobu
Hosono for his encouragement while writing this paper. From personal
conversations with Professor Mukai, I learned a lot of things around
his articles \cite{Mu1,Mu2}, and I have been strongly motivated to
get some results on $\mQ$-Fano 3-fold similar to his results. I appreciate
his generosity that gave me a lot of ideas. I got a sprout of the
research in this paper while I was staying at the Max-Planck-Institut
f\"ur Mathematik in 2007--2008. I appreciate after a long time the
institute providing a nice research environment with a free atmosphere.
Finally, I sincerely thank Professor Yujiro Kawamata, my thesis advisor,
for his appropriate guidance, encouragement and patience in the doctoral
course, which has been supporting my life as a mathematician. This
work is supported in part by Grant-in Aid for Scientific Research
(C) 16K05090.

\section{\textbf{Preliminaries \label{sec:Preliminaries}}}

\subsection{Miscellaneous results\label{subsec:Miscellaneous-results} }

The results in this subsection are frequently used in the sequel.
Proofs for them are omitted since they are elementary.
\begin{lem}
\label{lem:LinAlg}Let $V$ be a vector space, and $V=V^{1}\oplus V'$
be a direct sum decomposition with a $1$-dimensional subspace $V^{1}$
and a complementary subspace $V'$. The vector space $\wedge^{2}V$
has the following direct sum decomposition: 
\[
\wedge^{2}V=V'\oplus\wedge^{2}V',
\]
where we identify the subspace $V^{1}\wedge V'$ with $V'$. Let $U\subset\wedge^{2}V'$
be a subspace. For $\bm{x}\in V'$ and $\bm{y}\in U$, the following
are equivalent:

\begin{enumerate}[(1)]

\item $\bm{x}+\bm{y}\in\rG(2,V)\cap\mP(V'\oplus U).$

\item There exists a $2$-dimensional subspace $V^{2}\subset V'$
such that $\wedge^{2}V^{2}\subset U$, $\bm{x}\in V^{2}$ and $\bm{y}\in\wedge^{2}V^{2}$.

\end{enumerate}
\end{lem}

\begin{lem}
\label{lem:ABFib} Let $S$ be a projective manifold and $\sA,\sB$
vector bundles on $S$ whose dual bundles are globally generated.
Let $U_{\sA}:=H^{0}(S,\sA^{*})^{*}$and $U_{\sB}:=H^{0}(S,\sB^{*})^{*}$.
Let $p\colon\mP_{S}(\sA\oplus\sB)\to S$ be the natural morphism and
$\mu\colon\mP_{S}(\sA\oplus\sB)\to\mP(U_{\sA}\oplus U_{\sB})$ the
morphism defined by the tautological linear system $|H_{\mP(\sA\oplus\sB)}|$.
The following assertions hold:

\vspace{3pt}

\noindent $(1)$ The projective bundle $\mP_{S}(\sA\oplus\sB)$ is
contained in $\mP(U_{\sA}\oplus U_{\sB})\times S$ as a subbundle,
and the morphism $\mu$ is nothing but the composite $\mP_{S}(\sA\oplus\sB)\hookrightarrow\mP(U_{\sA}\oplus U_{\sB})\times S\to\mP(U_{\sA}\oplus U_{\sB}).$
The pull-back of $\sO_{\mP(U_{\sA}\oplus U_{\sB})}(1)$ by this morphism
is the tautological line bundle of $\mP_{S}(\sA\oplus\sB)$.

\vspace{3pt}

\noindent $(2)$ For a point $s\in S$, let $\sA_{s}$ and $\sB_{s}$
the fibers of $\sA$ and $\sB$ at $s$ respectively, which are subspaces
of $U_{\sA}$ and $U_{\sB}$ respectively. The $\mu$-image coincides
the locus
\[
\left\{ [\bm{x}+\bm{y}]\in\mP(U_{\sA}\oplus U_{\sB})\mid\exists_{s\in S},\bm{x}\in\sA_{s},\bm{y}\in\sB_{s}\right\} 
\]
 and the $\mu$-fiber over a point $[\bm{x}+\bm{y}]$ coincides with
the locus $\left\{ s\in S\mid\bm{x}\in\sA_{s},\bm{y}\in\sB_{s}\right\} .$
\end{lem}

Lemma \ref{lem:ABFib} also holds for a direct sum of three or more
vector bundles.

\subsection{Indecomposability of the mid point \label{subsec:IndecMid}}

In this subsection, we quickly review the classification of the mid
point $W$ of the Sarkisov link (\ref{eq:Sarkisov}) with a few compensation.
An important concept for the classification is indecomposability of
an effective divisor due to Mukai.
\begin{defn}
Let $X$ be a normal projective variety and $D$ a Weil divisor on
$X$. We say that $D$ is \textit{indecomposable} if there exists
no Weil divisors $A$ and $B$ such that $D\sim A+B$ and $h^{0}(A)\geq2$
and $h^{0}(B)\geq2$. If $-K_{X}$ is indecomposable, then we say
$X$ is indecomposable.
\end{defn}

An indecomposable $\mQ$-Fano variety generalizes a prime $\mQ$-Fano
variety for possibly non $\mQ$-factorial $\mQ$-Fano variety. In
our context, we have the following:
\begin{prop}
\label{prop: indec} The mid point $W$ is indecomposable. In the
genus 5,6, or 8 case, the anti-canonical divisor $-K_{W}$ is very
ample. The image $\Pi$ of $E$ on $W$ is a plane. 
\end{prop}

\begin{proof}
Assume by contradiction that there exist Weil divisors $A$ and $B$
such that $-K_{W}\sim A+B$ and $h^{0}(A)\geq2$ and $h^{0}(B)\geq2$.
Then we have $-K_{Z}\sim h_{*}^{-1}(A)+h_{*}^{-1}(B)+\Delta$ with
an effective $h$-exceptional divisor $\Delta$ (possibly equal to
$0)$ since $g$ is crepant. Further we have $-K_{X}\sim g_{*}h_{*}^{-1}(A)+g_{*}h_{*}^{-1}(B)+g_{*}(\Delta)$
with $h^{0}(g_{*}h_{*}^{-1}(A))\geq2$ and $h^{0}(g_{*}h_{*}^{-1}(B))\geq2$.
This is a contradiction since $X$ is prime. 

In the genus 5,6, or 8 case, $W$ is Gorenstein. Hence, by \cite[Thm.6.5 (2) and Prop.7.8]{Mu2},
$-K_{W}$ is very ample. Since $(-K_{W})^{2}E=1$, we see that $\Pi$
is a plane.
\end{proof}
\begin{cor}
\label{Cor:projPi} In the genus 5, 6, or 8 case, the rational map
$W\dashrightarrow X'$ in the Sarkisov link (\ref{eq:Sarkisov}) is
the projection of $W$ from the plane $\Pi$.
\end{cor}

\begin{proof}
In the genus 5, 6, or 8 case, we have $E'\sim z(-K_{Y'})-(z+1)\widetilde{E}$,
where $z+1$ is the Fano index of $X'$ by \cite[Part I, Table 1]{Tak1}.
Since $f'$ is the blow-up along the curve $C,$ we have $-K_{Y'}=f'^{*}(-K_{X'})-E'$.
Therefore we obtain $-K_{Y'}-\widetilde{E}=f'^{*}H_{X'}$ where $H_{X'}$
is the primitive Cartier divisor. This implies the assertion since
$-K_{Y'}$ is the pull-back of $\sO_{W}(1)$. 
\end{proof}

\subsection{Mid point in the genus 4 case \label{subsec:Genus-4mid}}

In the genus 4 case, we slightly modify the Sarkisov link (\ref{eq:Sarkisov})
partially. 

Let $g\colon Z\to X$ be the blow-up at the two $1/2$-singularities,
and $E_{1}$ and $E_{2}$ be the $g$-exceptional divisors. By \cite[Part II, Thm.1.0]{Tak1},
$-K_{Z}$ is nef and big. Let $h\colon Z\to\overline{Z}$ be the anti-canonical
model and $\Pi_{1}$ and $\Pi_{2}$ the images on $\overline{Z}$
of $E_{1}$ and $E_{2}$ respectively. In the same way as in the proof
of Proposition \ref{prop: indec}, we see that $\overline{Z}$ is
indecomposable, $-K_{\overline{Z}}$ is very ample and $\Pi_{1}$
and $\Pi_{2}$ are planes on $\overline{Z}$. By \cite[Thm.6.5 (2) and Prop.7.8]{Mu2},
$\overline{Z}$ is the intersection of a quadric and a cubic. 
\begin{prop}
\label{Prop:disj} The following assertions hold:

\vspace{3pt}

\noindent $(1)$ It holds that $\Pi_{1}\cap\Pi_{2}=\emptyset$. 

\vspace{3pt}

\noindent $(2)$ $h\colon Z\to\overline{Z}$ is a crepant small contraction,
and hence $\overline{Z}$ has only terminal singularities.
\end{prop}

\begin{proof}
First of all, we show that $h$ has no exceptional divisor $G$ such
that $h(G)$ is a point. Assume by contraction that there is such
a divisor $G$. Since $G$ is a crepant divisor while $-K_{X}$ is
ample, we see that $G$ intersects $E_{1}$ or $E_{2}.$ We may assume
that $G\cap E_{1}\not=\emptyset$. Then any irreducible component
of $G\cap E_{1}$ is a curve which is at the same time numerically
trivial and negative for $K_{Y}$, a contradiction. 

Assume by contradiction that $\Pi_{1}\cap\Pi_{2}\not=\emptyset$.
Let $\mathsf{t}$ be a point of $\Pi_{1}\cap\Pi_{2}$. Let $\gamma_{1}\cup\dots\cup\gamma_{k}$
be the irreducible decomposition of the fiber of $Z\to W$ over $\mathsf{t}$,
where any $\gamma_{l}$ is a curve by the first paragraph. For any
$l$, it holds that $E_{1}\cdot\gamma_{l}>0$ or $E_{2}\cdot\gamma_{l}>0$,
and there exist $i$ and $j$ such that $E_{1}\cdot\gamma_{i}>0$
and $E_{2}\cdot\gamma_{j}>0$. We fix a curve $\gamma_{j}$ such that
$E_{2}\cdot\gamma_{j}>0$. Since $E_{2}$ is mapped isomorphically
onto the plane $\Pi_{2}$, we have $E_{2}\cdot\gamma_{j}=1$. Now
we consider $Y$ in the Sarkisov link (\ref{eq:Sarkisov}) as the
target of the contraction of $E_{2}$ from $Z$ and $E$ as the image
of $E_{1}$. Let $\gamma_{j}'$ be the image of $\gamma_{j}$ on $Y$.
We have $-K_{Y}\cdot\gamma_{j}'=\frac{1}{2}$ since $E_{2}\cdot\gamma_{j}=1$.
Let $\gamma_{j}''$ be the strict transform of $\gamma_{j}'$ on $Y'$.
By a property of flop (cf.\cite{Ko}), we have $-K_{Y'}\cdot\gamma_{j}''=-K_{Y}\cdot\gamma_{j}'=\frac{1}{2}$.
By \cite[Part I, Table 1]{Tak1}, we have $E'\sim4(-K_{Y'})-5\widetilde{E}$.
This implies that $E'\cdot\gamma_{j}''=2-5\widetilde{E}\cdot\gamma_{j}''$.
If $\widetilde{E}\cdot\gamma_{j}''>0$, then $E'\cdot\gamma_{j}''\leq-3$
and hence $\gamma_{j}''\subset E'$. This is impossible since $E'$
does not contain the $\frac{1}{2}$-singularity of $Y'$ while $\gamma_{j}''$
contains it. Thus $\widetilde{E}\cdot\gamma_{j}''=0$, and hence $E_{1}\cdot\gamma_{j}=0$
and $\gamma_{j}$ does not intersect any $\gamma_{l}$ such that $E_{1}\cdot\gamma_{l}>0$.
This implies that the fiber of $Z\to W$ over $\mathsf{t}$ is disconnected,
a contradiction. 

If $h$ is a crepant divisorial contraction, then, by the first paragraph,
the $h$-exceptional locus contains a prime divisor, say $G$, such
that $h(G)$ is a curve, and $G$ intersects $E_{1}$ or $E_{2}$.
Assume that $E_{1}\cap G\not=\emptyset$ and $E_{2}\cap G\not=\emptyset$.
By the argument of the first paragraph, $E_{i}\cap G$ ($i=1,2$)
cannot contain an $h$-exceptional curve, and hence $E_{i}\cap G$
dominates $h(G)$. This implies that $\Pi_{1}\cap\Pi_{2}\not=\emptyset$,
a contradiction to (1). Therefore, we may assume that $E_{1}\cap G\not=\emptyset$
and $E_{2}\cap G=\emptyset$. Again, we consider $Y$ in the Sarkisov
link (\ref{eq:Sarkisov}) as the target of the contraction of $E_{2}$
from $Z$ and $E$ as the image of $E_{1}$. Then the image of $G$
on $Y$ is a crepant divisor. This is impossible since $Y$ has no
crepant divisorial contraction by \cite[Part I, Table 1]{Tak1}.
\end{proof}
\noindent \textbf{Caution (change of notation):} In the genus 4 case,
henceforth we set 
\[
W:=\overline{Z}
\]
 for notational convenience. We also call this $W$ the \textit{mid
point} in the genus 4 case.

\subsection{Mid point in the genus 6 case \label{subsec:g6setup}}

By \cite[Thm.6.5 (2) and Prop.7.8]{Mu2}, the mid point $W$ is a
quadric section of a del Pezzo $4$-fold $W_{0}$ with only canonical
singularities. The indecomposability simplifies the situation as follows:
\begin{prop}
\label{prop:sm} The del Pezzo $4$-fold $W_{0}$ is smooth.
\end{prop}

\begin{proof}
Assume that $W_{0}$ is singular and is not a cone over the smooth
quintic del Pezzo $3$-fold $B_{5}$ (we do not exclude the possibility
that $W_{0}$ is a cone over a singular quintic del Pezzo $3$-fold).
Then, by \cite[p.160, (6)]{Fuj3}, $W_{0}$ contains a double point,
say, $\mathsf{t}$. By projecting $W_{0}$ from $\mathsf{t}$, $W_{0}$
is mapped onto a non-degenerate cubic $4$-fold $W_{0}'$ in $\mP^{6}$.
By the classification of $W_{0}'$, we see that $\sO_{W_{0}'}(1)$
is decomposable. This implies that $\sO_{W_{0}}(1)$ is, and hence
$-K_{W}$ is decomposable, a contradiction. 

Assume that $W_{0}$ is the cone over \textbf{$B_{5}$}. If $W$ do
not contain the vertex of $W_{0}$, then, by projecting $W_{0}$ from
the vertex, $B_{5}$ contains the plane which is the image of $\Pi$.
This is absurd since $B_{5}$ does not contain a plane. Therefore
$W$contains the vertex of $W_{0}$, and hence $W$ has a non-hypersurface
singularity at the vertex. This is again absurd since $W$ has only
Gorenstein terminal singularities. 
\end{proof}
By \cite{Fuj3}, we can write $W_{0}=\rG(2,V)\cap\mP(U^{8})$ with
$V\simeq\mC^{5}$ and $U^{8}\simeq\mC^{8}\subset\wedge^{2}V$. We
write $W=W_{0}\cap Q$, where $Q$ is a quadric 6-fold in $\mP(U^{8})$.
It is well-known that the 2-plane $\Pi$ has one of the following
description as a subvariety of $\rG(2,V)$:
\begin{enumerate}
\item 
\[
\Pi=\{[\mC^{2}]\mid V^{1}\subset\mC^{2}\subset V^{4}\}\subset\rG(2,V)
\]
with some fixed vector subspaces $V^{1}\simeq\mC$ and $V^{4}\simeq\mC^{4}$
of $V$.
\item 
\[
\Pi=\{[\mC^{2}]\mid\mC^{2}\subset V^{3}\}\subset\rG(2,V)
\]
with a fixed vector subspace $V^{3}\simeq\mC^{3}$ of $V$.
\end{enumerate}
\begin{prop}
\label{Prop:QC} If $X$ is of $\mathtt{Q}$-type, then $\Pi$ satisfies
$(1)$. If $X$ is of $\mathtt{C}$-type, then $\Pi$ satisfies $(2)$.
\end{prop}

\begin{proof}
Since $W\dashrightarrow X'$ is the projection from $\Pi$ by Corollary
\ref{Cor:projPi}, $X'$ is contained in the image $W_{0}'$ of the
projection of $W_{0}$ from $\Pi$. By \cite{Fuj3} (see also \cite[Lem.3.4.4]{IsP}),
$W_{0}'$ is a smooth quadric 3-fold $Q^{3}$ if $\Pi$ satisfies
(1), or is $\mP^{4}$ if $\Pi$ satisfies (2). If $X$ is of $\mathtt{C}$-type,
then $X'$ is a smooth cubic 3-fold, hence $\Pi$ must satisfy (2).
Assume by contradiction that $X$ is of $\mathtt{Q}$-type and $\Pi$
satisfies (2). Then, since $X'$ is a smooth quadric 3-fold $Q^{3}$,
we may choose the quadric $6$-fold $Q$ as the cone over $Q^{3}$
with the vertex $\Pi$ . This implies that $W$ is singular along
$\Pi$, a contradiction.
\end{proof}
\noindent \textbf{$\mathtt{Q}$-type}: 

\vspace{3pt}

\noindent \textbf{Caution} \textbf{(change of notation)}: Hereafter,
for notational convenience, we denote by $\Pi_{0}$ the $2$-plane
$\Pi$ only in this case. 

\vspace{3pt}

Note that $\Pi_{0}$ satisfies (1). The notation $\Pi$ will denote
the unique $3$-plane 

\[
{\Pi}:=\{[\mC^{2}]\mid V^{1}\subset\mC^{2}\}\subset\rG(2,V)
\]
containing $\Pi_{0}$. It holds that $\Pi_{0}={\Pi}\cap\mP(U^{8})$.
By a simple dimension count of linear subspaces, the linear hull $\mP(U^{9})$
of $\Pi\cup\mP(U^{8})$ is a hyperplane of $\mP(\wedge^{2}V).$ We
set 
\[
A_{\mathtt{Q}}:=\rG(2,V)\cap\mP(U^{9}).
\]
Since $W_{0}$ is a linear section of $\rG(2,V)$, so is $A_{\mathtt{Q}}$.
Thus $A_{\mathtt{Q}}$ is not a cone over $W_{0}$ since otherwise
$A_{\mathtt{Q}}$ cannot be contained in ${\rm G}(2,5)$. We can show
that $A_{\mathtt{Q}}$ is actually smooth in the same way as the proof
of Proposition \ref{prop:sm}. Now we produce the situation as in
Lemma \ref{lem:LinAlg}. We choose a direct sum decomposition $V=V^{1}\oplus V'$
with a complementary subspace $V'$ to $V^{1}$. Note that $\Pi=\mP(V^{1}\wedge V'$),
which we identify with $\mP(V')$. There exists a $5$-dimensional
subspace $U^{5}\subset\wedge^{2}V'$ such that $U^{9}=V'\oplus U^{5}$. 

The projection of ${\rm G}(2,V)$ from the $3$-plane $\Pi$ induces
the natural rational map $A_{\mathtt{Q}}\dashrightarrow{\rm G(2,V')}\cap\mP(U^{5})$
and the target ${\rm G(2,V')}\cap\mP(U^{5})$ is nothing but the smooth
quadric 3-fold $Q^{3}$ as in the proof of Proposition \ref{Prop:QC}.

By \cite{Fuj3}, the pair $(A_{\mathtt{Q}},{\Pi})$ is unique up to
projective equivalence. We may take the following coordinates: Let
$\bm{e}_{i}\,(1\leq i\leq5)$ be a basis of $V$. We set $V^{1}:=\mC\bm{e}_{1}$,
and $V':=$ the subspace of $V$ generated by $\bm{e}_{2},\dots,\bm{e}_{5}$.
Let $z_{i}\,(1\leq i\leq5)$ be the coordinate for $\bm{e}_{i}$ and
$x_{ij}$ the Pl\"ucker coordinate for $\bm{e}_{i}\wedge\bm{e}_{j}\,(1\leq i<j\leq5)$.
We set
\begin{align*}
\Pi:= & \left\{ x_{ij}=0\,(2\leq i\leq j\leq5)\right\} \subset\mP(\wedge^{2}V),\\
U^{9}:= & \left\{ x_{24}-x_{35}=0\right\} \subset\wedge^{2}V,\\
U^{5}:= & \left\{ x_{24}-x_{35}=0\right\} \subset\wedge^{2}V'.
\end{align*}
Moreover, dropping the coordinate $x_{24}$ by the equality $x_{24}=x_{35}$,
we consider $A_{\mathtt{Q}}$ as a subvariety of $\mP^{8}$ with coordinates
$z_{2},z_{3},z_{4},z_{5}$ and ${{\empty}^{t}\!\bm{x}}:=\left(\begin{array}{ccccc}
x_{23} & x_{25} & x_{34} & x_{35} & x_{45}\end{array}\right)$ defined by the following equations:

\vspace{3pt}

\noindent \textbf{Equation of $A_{\mathtt{Q}}$}

\begin{equation}
N_{\mathtt{Q}}\bm{x}=\bm{o},\,x_{23}x_{45}-x_{35}^{2}+x_{25}x_{34}=0,\label{eq:EqAQ}
\end{equation}

where we set
\[
N_{\mathtt{Q}}:=\left(\begin{array}{ccccc}
z_{4} & 0 & z_{2} & -z_{3} & 0\\
z_{5} & -z_{3} & 0 & z_{2} & 0\\
0 & -z_{4} & 0 & z_{5} & z_{2}\\
0 & 0 & z_{5} & -z_{4} & z_{3}
\end{array}\right).
\]

In this situation, 
\begin{align*}
\Pi= & \left\{ \bm{x}=\bm{o}\right\} ,\\
Q^{3}= & \left\{ x_{23}x_{45}-x_{35}^{2}+x_{25}x_{34}=0\right\} \subset\mP^{4}.
\end{align*}

\vspace{5pt}

\noindent \textbf{$\mathtt{C}$-type}:

\vspace{5pt}

In this case, we set
\[
A_{\mathtt{C}}:=W_{0}.
\]
Since $\Pi=\mP(\wedge^{2}V^{3})$ in this case, we may write $U^{8}=\wedge^{2}V^{3}\oplus U^{5}$
with some $U^{5}\simeq\mC^{5}$. The projection of ${\rm G}(2,V)$
from the $2$-plane $\Pi$ induces the natural rational map $A_{\mathtt{c}}\dashrightarrow\mP(U^{5})$
and the target $\mP(U^{5})$ is nothing but $\mP^{4}$ as in the proof
of Proposition \ref{Prop:QC}. Let $a\colon\widehat{A}_{\mathtt{C}}\to A_{\mathtt{C}}$
is the blow-up of $A_{\mathtt{C}}$ along $\Pi$. By a general property
of linear projection, a morphism $b\colon\widehat{A}_{\mathtt{C}}\to\mP(U^{5})$
is induced. By \cite[Sect.10]{Fuj3}, $b$ is the blow-up of $\mP(U^{5})$
along a twisted cubic $\gamma_{\mathtt{C}}$.

By \cite[Sect.10]{Fuj3} again, the pair $(A_{\mathtt{C}},{\Pi})$
is unique up to projective equivalence. For choices of coordinates
$x_{1},x_{2},x_{3}$ of $\wedge^{2}V^{3}$ and $y_{1},\dots,y_{5}$
of $U^{5}$, we may write the equation of $A_{\mathtt{C}}$ as follows:

\vspace{3pt}

\noindent \textbf{Equation of $A_{\mathtt{C}}$}

\begin{equation}
\left(\begin{array}{ccc}
y_{4} & y_{3} & y_{2}\\
y_{3} & y_{2} & y_{1}
\end{array}\right)\left(\begin{array}{c}
x_{1}\\
x_{2}\\
x_{3}
\end{array}\right)=\left(\begin{array}{c}
0\\
0
\end{array}\right),\,y_{5}\left(\begin{array}{c}
x_{1}\\
x_{2}\\
x_{3}
\end{array}\right)=\left(\begin{array}{c}
y_{2}^{2}-y_{1}y_{3}\\
y_{1}y_{4}-y_{2}y_{3}\\
y_{3}^{2}-y_{2}y_{4}
\end{array}\right),\label{eq:EqAc}
\end{equation}
 where the twisted cubic $\gamma_{\mathtt{C}}$ is equal to 
\[
\left\{ y_{2}^{2}-y_{1}y_{3}=y_{1}y_{4}-y_{2}y_{3}=y_{3}^{2}-y_{2}y_{4}=y_{5}=0\right\} .
\]

\subsection{Mid point in the genus 8 case \label{subsec:Mid-point-in G8}}

Let $V$ be a $6$-dimensional vector space. By \cite[Thm.6.5 (2) and Prop.7.8]{Mu2},
$W$ is a codimension 5 linear section of $\rG(2,V)$. We write $W=\rG(2,V)\cap\mP(U^{10})$
with a 10-dimensional subspace $U^{10}\subset\wedge^{2}V$. Note that
the image of a fiber of $E'\to C$ on $W$ in the Sarkisov link (\ref{eq:Sarkisov})
is a line intersecting $\Pi$ since the equality $E'\sim-K_{Y'}-2\widetilde{E}$
holds by \cite[Part I, Table 1]{Tak1}. Therefore, by \cite[Lem.5.3]{Tak2},
the 2-plane $\Pi$ has the following description as a subvariety of
$\rG(2,V)$:

\[
\Pi=\{[\mC^{2}]\mid V^{1}\subset\mC^{2}\subset V^{4}\}\subset\rG(2,V)
\]
with some fixed vector subspaces $V^{1}\simeq\mC$ and $V^{4}\simeq\mC^{4}$
of $V$. 

\section{\textbf{Embedding theorem in the genus 8 case }\label{sec:Genus8}}

\subsection{Extending the mid point \label{subsec:Extending-the-mid G8}}

The 2-plane $\Pi$ is contained in the following 4-plane 

\[
\overline{\Pi}:=\{[\mC^{2}]\mid V^{1}\subset\mC^{2}\}\subset\rG(2,V),
\]
and it holds that $\Pi=\overline{\Pi}\cap\mP(U^{10})$. By a simple
dimension count of linear subspaces, the linear hull $\mP(U^{12})$
of $\overline{\Pi}\cup\mP(U^{10})$ is a codimension 3 linear subspace
of $\mP(\wedge^{2}V).$ We set 
\[
\overline{\Sigma}:=\rG(2,V)\cap\mP(U^{12}).
\]
Since $W$ is a linear section of $\rG(2,V)$, so is $\overline{\Sigma}$.
Now we produce the situation as in Lemma \ref{lem:LinAlg}. We choose
a direct sum decomposition $V=V^{1}\oplus V'$ with a complementary
subspace $V'$ to $V^{1}$. Note that $\Pi=\mP(V^{1}\wedge V'$),
which we identify with $\mP(V')$. There exists a $5$-dimensional
subspace $U^{7}\subset\wedge^{2}V'$ such that $U^{12}=V'\oplus U^{7}$. 

The projection of ${\rm G}(2,V)$ from the 4-plane $\overline{\Pi}$
induces the natural rational map $\overline{\Sigma}\dashrightarrow{\rm G(2,V')}\cap\mP(U^{7})$.
This also induces the rational map $W\dashrightarrow{\rm G(2,V')}\cap\mP(U^{7})$,
which is the projection from $\Pi$. By Corollary \ref{Cor:projPi},
we have ${\rm G(2,V')}\cap\mP(U^{7})=X'\simeq B_{5}$. 

\subsection{Construction of the key variety \label{subsec:Construction-of-the key G8}}
\begin{defn}
\label{Def:bundleG8} Let $\sU$ be the rank two universal subbundle
on $\rG(2,V')\simeq\rG(2,5)$. Set 
\[
\Sigma':=\mP_{B_{5}}(\sU|_{B_{5}}\oplus\sO_{B_{5}}(-1)).
\]
 Note that, by a standard computation, it follows that 
\begin{equation}
-K_{\Sigma'}=3H_{\Sigma'}.\label{eq:-3H}
\end{equation}
\end{defn}

To investigate the birational geometry of $\Sigma'$, we need the
following beautiful classical result: 
\begin{lem}
\label{lem:projVeronese} The natural morphism $\mathbb{P}_{B_{5}}(\sU|_{B_{5}})\to\mP(V')$
to $\mP(V')\simeq\mP^{4}$ from the total space of lines in $\mP^{4}$
parameterized by $B_{5}\subset\rG(2,V')$ is the blow-up along the
projected Veronese surface $\sV$.
\end{lem}

\begin{proof}
By \cite[Prop.2.4 (b)]{Il}, $B_{5}$ parameterizes tri-secant lines
of the projected second Veronese surface $\sV\subset\mP^{4}$. Let
$\mathsf{s}$ be a point of $\mP^{4}$. The fiber of $\mathbb{P}_{B_{5}}(\sU|_{B_{5}})\to\mP(V')$
parameterizes the tri-secant lines of $\sV$ through $\mathsf{s}.$
If $\mathsf{s}\not\in\sV$, then there is a unique tri-secant line
of $\sV$ through $\mathsf{s}$ (this is classically known and follows
by \cite[Lem.8.1]{L} for example), and if $\mathsf{s}\in\sV$, then
tri-secant lines of $\sV$ through $\mathsf{s}$ are parameterized
by $\mP^{1}$ (\cite[Prop.2.4 (a)]{Il}). Note that $-K_{\mathbb{P}_{B_{5}}(\sU|_{B_{5}})}=2H_{\mathbb{P}_{B_{5}}(\sU|_{B_{5}})}+L$,
where $L$ is the pull-back of $\sO_{B_{5}}(1)$. Since a fiber of
$\mathbb{P}_{B_{5}}(\sU|_{B_{5}})\to\mP(V')$ is not contained in
a fiber of $\mathbb{P}_{B_{5}}(\sU|_{B_{5}})\to B_{5}$, it is positive
for $-K_{\mathbb{P}_{B_{5}}(\sU|_{B_{5}})}$. Therefore, by \cite[Thm.2.3]{An},
$\mathbb{P}_{B_{5}}(\sU|_{B_{5}})\to\mP^{4}$ is the blow-up along
$\sV$. 
\end{proof}
\begin{prop}
\label{Prop:g8Sigma'} The tautological linear system $|H_{\Sigma'}|$
defines a surjective birational morphism $\Sigma'\to\overline{\Sigma}$,
which we will denote by $\varphi_{|H_{\Sigma'}|}$. It is a flopping
contraction of Atiyah type. The image of the flopping locus on $\overline{\Sigma}$
is a projected second Veronese surface $\mathsf{\mathcal{V}}$ in
$\overline{\Pi}.$ 
\end{prop}

\begin{proof}
Take a point $\mathsf{p}:=[\wedge^{2}V^{2}]\in B_{5}=\rG(2,V')\cap\mP(U^{7})$,
where $V^{2}\subset V'$ is a $2$-dimensional subspace such that
$\wedge^{2}V^{2}\subset U^{7}$. The fiber of the projection $\Sigma'\to B_{5}$
over $\mathsf{p}$ is $\mP(V^{2}\oplus\wedge^{2}V^{2})$, which is
a linear subspace of $\mP(V'\oplus U^{7})$. 

For a $2$-dimensional subspace $V^{2}\subset V'$ such that $\wedge^{2}V^{2}\subset U^{7}$,
we take a point $[\bm{x}+\bm{y}]\in\mP(V^{2}\oplus\wedge^{2}V^{2})$
with $\bm{x}\in V^{2}$ and $\bm{y}\in\wedge^{2}V^{2}$. By Lemma
\ref{lem:LinAlg}, it holds that $[\bm{x}+\bm{y}]\in G(2,V)\cap\mP(V'\oplus U^{7})=\overline{\Sigma}.$
Therefore the image of $\Sigma'\to\mP(V'\oplus U^{7})$ is contained
in $\overline{\Sigma}$, and hence the desired morphism $\Sigma'\to\overline{\Sigma}$
is induced. By Lemma \ref{lem:ABFib} (1), this is defined by the
tautological linear system $|H_{\Sigma'}|$. 

Let $\mathsf{t}:=[\bm{x}+\bm{y}]$ be a point of $\overline{\Sigma}$
with $\bm{x}\in V'$ and $\bm{y}\in U^{7}$. Then the fiber of $\Sigma'\to\overline{\Sigma}$
over t is $\left\{ \mathsf{t}\right\} \times\left\{ [\wedge^{2}V^{2}]\mid\bm{x}\in V^{2},\bm{y}\in\wedge^{2}V^{2}\subset U^{7}\right\} $
by Lemma \ref{lem:ABFib} (2), which is nonempty by Lemma \ref{lem:LinAlg}.
Therefore the morphism $\Sigma'\to\overline{\Sigma}$ is surjective.\textbf{
}If $\bm{y}\not=\bm{o},$ then $V^{2}$ is uniquely determined by
$\wedge^{2}V^{2}=\mC\bm{y}.$ Therefore the morphism $\Sigma'\to\overline{\Sigma}$
is birational. 

If $\bm{y}=\bm{o},$ then $\mathsf{t}$ is a point of $\overline{\Pi}$.
Note that the restriction of the morphism $\Sigma'\to\overline{\Sigma}$
over $\overline{\Pi}$ is $\mathbb{P}_{B_{5}}(\sU|_{B_{5}}\oplus0)\to\mP(V'\oplus0)=\overline{\Pi}\simeq\mP^{4},$
which can be identified with the natural morphism $\mathbb{P}_{B_{5}}(\sU|_{B_{5}})\to\mP(V')$
to $\mP(V')\simeq\mP^{4}$ from the total space of lines in $\mP^{4}$
parameterized by $B_{5}\subset\rG(2,V')$. Let $l$ be the fiber of
$\mathbb{P}_{B_{5}}(\sU|_{B_{5}}\oplus0)\to\overline{\Pi}$ over a
point of the projected Veronese surface $\sV\subset\mP(V')$. We compute
the normal bundle $\sN_{l/\Sigma'}$. Since $\mathbb{P}_{B_{5}}(\sU|_{B_{5}}\oplus0)\to\overline{\Pi}$
is the blow-up of $\overline{\Pi}$ along $\sV$ by Lemma \ref{lem:projVeronese},
we see that $\sN_{l/\mathbb{P}_{B_{5}}(\sU|_{B_{5}}\oplus0)}=\sO_{l}^{\oplus2}\oplus\sO_{l}(-1)$.
Let $L_{\Sigma'}$ be the pull-back of $\sO_{B_{5}}(1)$ on $\Sigma'$.
Since $\mathbb{P}_{B_{5}}(\sU|_{B_{5}}\oplus0)$ is linearly equivalent
to $H_{\Sigma'}-L_{\Sigma'}$ in $\Sigma'$, we see that $\sN_{\mathbb{P}_{B_{5}}(\sU|_{B_{5}}\oplus0)/\Sigma'}|_{l}=\sO_{l}(-1)$.
Therefore, by the normal bundle sequence $0\to\sN_{l/\mathbb{P}_{B_{5}}(\sU|_{B_{5}}\oplus0)}\to\sN_{l/\Sigma'}\to\sN_{\mathbb{P}_{B_{5}}(\sU|_{B_{5}}\oplus0)/\Sigma'}|_{l}\to0$,
we see that $\sN_{l/\Sigma'}=\sO_{l}^{\oplus2}\oplus\sO_{l}(-1)^{\oplus2}$,
and $\Sigma'\to\overline{\Sigma}$ is a flopping contraction of Atiyah
type as desired. 
\end{proof}
Let $\Sigma'\dashrightarrow\widetilde{\Sigma}$ be the flop for this
flopping contraction. Let $\Pi'$ and $\widetilde{\Pi}$ be the strict
transforms of $\overline{\Pi}$ on $\Sigma'$ and $\widetilde{\Sigma}$
respectively. It is well-known that the flop can be constructed by
the blow-up along the flopping locus and the blow-down of the exceptional
divisor along the other direction. From this, we see that the restriction
$\Pi'\dashrightarrow\widetilde{\Pi}$ of the flop is the blow-up $\Pi'\to\overline{\Pi}$
of $\overline{\Pi}$ along $\sV$. Therefore $\widetilde{\Pi}$ is
isomorphic to $\mP^{4}$. Let $H_{\widetilde{\Sigma}}$ be the strict
transform on $\widetilde{\Sigma}$ of $H_{\Sigma'}$.
\begin{lem}
\label{lem:P/PG8} The normal bundle $\sN_{\widetilde{\Pi}/\widetilde{\Sigma}}$
is $\sO_{\mP^{4}}(-2)$. 
\end{lem}

\begin{proof}
We have $H_{\widetilde{\Sigma}}|_{\widetilde{\Pi}}=\sO_{\mP^{4}}(1)$
since $\sO_{\overline{\Sigma}}(1)|_{\overline{\Pi}}=\sO_{\mP^{4}}(1)$.
Therefore we have $-K_{\widetilde{\Sigma}}|_{\widetilde{\Pi}}=\sO_{\mP^{4}}(3)$.
Since $-K_{\widetilde{\Pi}}=\sO_{\mP^{4}}(5)$, we have $\sN_{\widetilde{\Pi}/\widetilde{\Sigma}}\simeq\sO_{\mP^{4}}(-2)$
as desired. 
\end{proof}
\begin{lem}
\label{lem:2H+P G8} $2H_{\widetilde{\Sigma}}+\widetilde{\Pi}$ is
semiample. 
\end{lem}

\begin{proof}
We show that $2H_{\widetilde{\Sigma}}+\widetilde{\Pi}$ is nef. Assume
that $(2H_{\widetilde{\Sigma}}+\widetilde{\Pi})\cdot\gamma<0$ for
an irreducible curve $\gamma$. Then $\widetilde{\Pi}\cdot\gamma<0$
since $H_{\widetilde{\Sigma}}$ is nef, and hence $\gamma\subset\widetilde{\Pi}$.
Since $\widetilde{\Pi}|_{\widetilde{\Pi}}=\sO_{\mP^{4}}(-2)$ and
$H_{\widetilde{\Sigma}}|_{\widetilde{\Pi}}=\sO_{\mP^{4}}(1)$, we
have $(2H_{\widetilde{\Sigma}}+\widetilde{\Pi})\cdot\gamma=0$, a
contradiction. Therefore $2H_{\widetilde{\Sigma}}+\widetilde{\Pi}$
is nef.

Since $-K_{\widetilde{\Sigma}}$ is nef and big, $2H_{\widetilde{\Sigma}}+\widetilde{\Pi}$
is semiample by the Kawamata-Shokurov base point free theorem (\cite{KMM}). 
\end{proof}
\begin{thm}
\label{thm:key8}Let $\mu\colon\widetilde{\Sigma}\to\Sigma$ be the
contraction defined by a sufficient multiple of $2H_{\widetilde{\Sigma}}+\widetilde{\Pi}$.
The exceptional locus of this contraction is $\widetilde{\Pi}$. The
image of $\widetilde{\Pi}$ on $\Sigma$ is a $1/2$-singularity.
$\Sigma$ is a $5$-dimensional rational $\mQ$-Fano variety with
only one $1/2$-singularity and with $\rho(\Sigma)=1$. The image
$M_{\Sigma}$ of $H_{\widetilde{\Sigma}}$ is a primitive Weil divisor
and it holds that $-K_{\Sigma}=3M_{\Sigma}$.
\end{thm}

\begin{proof}
As we have checked in the proof of Lemma \ref{lem:2H+P G8}, $2H_{\widetilde{\Sigma}}+\widetilde{\Pi}$
is numerical trivial for any curve in $\widetilde{\Pi}$. Thus the
image of $\widetilde{\Pi}$ by $\widetilde{\Sigma}\to\Sigma$ is a
$1/2$-singularity by Lemma \ref{lem:P/PG8}. Assume by contradiction
that $(2H_{\widetilde{\Sigma}}+\widetilde{\Pi})\cdot\gamma=0$ for
an irreducible curve $\gamma\not\subset$$\widetilde{\Pi}$. Since
$H_{\widetilde{\Sigma}}$ is nef and $\gamma\not\subset\widetilde{\Pi}$,
we have $H_{\widetilde{\Sigma}}\cdot\gamma=\widetilde{\Pi}\cdot\gamma=0$.
By the condition that $H_{\widetilde{\Sigma}}\cdot\gamma=0$, $\gamma$
is a flopped curve. This is absurd since $\widetilde{\Pi}$ is positive
for a flopped curve. 

We show that $\rho(\Sigma$)=1. Since $\Sigma'\to B_{5}$ is a projective
bundle, we see that $\rho(\Sigma')=\rho(B_{5})+1=2$. Since $\Sigma'\dashrightarrow\widetilde{\Sigma}$
is a flop, we have $\rho(\widetilde{\Sigma})=\rho(\Sigma')=2$. Finally,
since $\widetilde{\Sigma}\to\Sigma$ contracts a divisor, we have
$\rho(\Sigma)\leq\rho(\widetilde{\Sigma})-1=1$. Hence we have $\rho(\Sigma)=1$.

The equality $-K_{\Sigma}=3M_{\Sigma}$ follows from (\ref{eq:-3H}).
We show that $M_{\Sigma}$ is primitive. If $M_{\Sigma}$ were not
primitive, then $M_{\Sigma}$ would be written as $M_{\Sigma}=\alpha M'_{\Sigma}$
with a primitive Weil divisor $M'_{\Sigma}$ and positive integer
$\alpha\geq2.$ Since $\Sigma$ has only a $1/2$-singularity, $2M'_{\Sigma}$
is a Cartier divisor. Therefore we have $2H_{\widetilde{\Sigma}}+F_{\widetilde{\Sigma}}=\alpha\mu^{*}(2M'_{\Sigma})$
and hence there is a Cartier divisor $D$ on $\Sigma'$ such that
$2H_{\Sigma'}+\Pi'=\alpha D$. Since $\Pi'\cdot l=-1$ for a flopping
curve $l$ by the proof of Proposition \ref{Prop:g8Sigma'}, this
implies that $\alpha D\cdot l=-1$, which is impossible if $\alpha\geq2$.
Therefore $M_{\Sigma}$ is primitive. 

The rationality of $\Sigma$ follows since $\Sigma$ is birational
to the projective bundle $\Sigma'$ over the rational Fano 3-fold
$B_{5}$.
\end{proof}

\subsection{Embedding theorem \label{subsec:Embedding-theorem G8}}

Now we show Theorem \ref{thm:main1} for a prime $\mQ$-Fano $3$-fold
$X$ of genus 8.
\begin{thm}
\label{thm:emb8} A $\mQ$-Fano $3$-fold $X$ of genus $8$ is a
linear section of $\Sigma$. 
\end{thm}

\begin{proof}
Note that $W\cap\Sing\overline{\Sigma}$ is $0$-dimensional since
$W$ has only terminal singularities and $W$ is a linear section
of $\overline{\Sigma}$ with respect to $|\sO_{\overline{\Sigma}}(1)|$.
Therefore, since $\widetilde{\Sigma}\to\overline{\Sigma}$ is crepant
and small and nontrivial fibers are 1-dimensional, the strict transform
$W_{\widetilde{\Sigma}}$ of $W$ in $\widetilde{\Sigma}$ is a linear
section of $\widetilde{\Sigma}$ with respect to $|H_{\widetilde{\Sigma}}|$
and hence the restriction $W_{\widetilde{\Sigma}}\to W$ of $\widetilde{\Sigma}\to\overline{\Sigma}$
over $W$ is also crepant and small. Since $W$ has only terminal
singularities and $W_{\widetilde{\Sigma}}\to W$ is crepant, we see
that $W_{\widetilde{\Sigma}}$ is normal and has only terminal singularities
by \cite[the proof of Prop.16.4]{CKM}. Note that $\widetilde{\Pi}$
is relatively ample for $W_{\widetilde{\Sigma}}\to W$. Since $Y\to W$
is the unique small extraction such that the strict transform of $\Pi$
is relatively ample, we see that $Y=W_{\widetilde{\Sigma}}$. Since
we may write $Y=W_{\widetilde{\Sigma}}=\widetilde{H}_{1}\cap\widetilde{H}_{2}$
with $\widetilde{H}_{i}\in|H_{\widetilde{\Sigma}}|$ ($i=1,2$), we
see that $X=M_{1}\cap M_{2}$ with the images $M_{i}\in|M_{\Sigma}|$
of $\widetilde{H}_{i}$ as desired.
\end{proof}

\subsection{Extension of the Sarkisov link}

By the proof of Theorem \ref{thm:emb8}, we have the following:
\begin{cor}
\label{cor:ExtSarkisovG8} The following diagram is an extension of
the Sarkisov link (\ref{eq:Sarkisov}):\begin{equation}\label{eq:SarkisovG8} \xymatrix{& \widetilde{\Sigma}\ar@{-->}[rr]^{\rm{flop}}\ar[dl]\ar[dr] & & \Sigma'\ar[dr]\ar[dl]\\
\Sigma & &\overline{\Sigma} & & B_5,}
\end{equation}where $\Sigma$, $\widetilde{\Sigma},$ $\overline{\Sigma}$ and $\Sigma'$
are extensions of $X$, $Y$, $W$ and $Y'$ respectively.
\end{cor}

\section{\textbf{Extending the mid point in the genus 4 or 6 case \label{sec:Extending Mid G46}}}

In this section, we extend the mid point $W$ to a certain variety
$\overline{\Sigma}$ in the genus 4 or 6 case. We also construct a
crepant small resolution of $\overline{\Sigma}$ in case of genus
4 or genus 6 and $\mathtt{Q}$-type, which will be the main ingredient
for Theorem \ref{thm:main2} in each of these cases. 

\subsection{Genus 4 \label{subsec:ExtGenus-4}}

\subsubsection{\textbf{Extending the mid point}}

As we have seen in the subsection \ref{subsec:Genus-4mid}, the mid
point $W$ is a complete intersection of a quadric and a cubic in
$\mP^{5}$ containing two disjoint planes $\Pi_{1}$ and $\Pi_{2}$.

To extend $W$, we start from mutually disjoint two planes $\Pi_{1}$
and $\Pi_{2}$ in $\mP^{5}$. By a coordinate change, we may assume
that $\Pi_{1}=\{x_{1}=x_{2}=x_{3}=0\}$ and $\Pi_{2}=\{y_{1}=y_{2}=y_{3}=0\}$
in $\mP^{5}$ with coordinates $x_{1},x_{2},x_{3},y_{1},y_{2},y_{3}$.
Set ${{\empty}^{t}\!\bm{x}}=\begin{pmatrix}x_{1} & x_{2} & x_{3}\end{pmatrix}$
and ${\empty}^{t}\!\bm{y}=\begin{pmatrix}y_{1} & y_{2} & y_{3}\end{pmatrix}$.
Take a quadric $\mathsf{\mathsf{Q_{0}}}$ and a cubic $\mathtt{\mathsf{C_{0}}}$
both of which contain $\Pi_{1}$ and $\Pi_{2}$. Then we may write
\[
\mathsf{Q_{0}}=\{{\empty}^{t}\!\bm{y}M_{0}\bm{x}=0\},\,\mathsf{C_{0}}=\{{\empty}^{t}\!\bm{y}M_{1}\bm{x}=0\},
\]
where $M_{0}$ (resp.~$M_{1}$) is a $3\times3$ matrix with constant
entries (resp.~linear entries).

Remarkably, the indecomposability simplifies the situation as follows:
\begin{lem}
\label{lem:rank3} It holds that $\rank M_{0}=3$. We may assume that
$M_{0}$ is the identity matrix by a coordinate change, namely, 
\[
\mathsf{Q_{0}}=\{{\empty}^{t}\!\bm{y}\bm{x}=0\}.
\]
 
\end{lem}

\begin{proof}
If $\rank M_{0}=1$, then $\mathsf{Q_{0}}\cap\mathsf{C_{0}}$ is reducible,
a contradiction. Assume that $\rank M_{0}=2$. Then $\mathsf{Q}_{0}$
is the cone over $\mP^{1}\times\mP^{1}$. Therefore $\mathsf{Q_{0}}$
contains two families of $3$-planes $\{\mathsf{P}_{a}\}_{a\in\mP^{1}}$
and $\{\mathsf{P}_{b}\}_{b\in\mP^{1}}$ such that the sums $\mathsf{P}_{a}+\mathsf{P}_{b}$
are hyperplane sections of $\mathsf{Q_{0}}$. Hence $(\mathsf{P}_{a}+\mathsf{P}_{b})\cap\mathsf{Q_{0}}\cap\mathsf{C_{0}}$
are anti-canonical divisors, a contradiction to the indecomposability
of $\mathsf{Q_{0}}\cap\mathsf{C_{0}}$. Therefore $\rank M_{0}=3$.
The latter assertion is obvious.
\end{proof}
Moreover, subtracting $(1/3\,{\rm Tr}\,M_{1})$ times the equation
of $\mathsf{Q_{0}}$ from the equation of $\mathsf{C_{0}}$, we may
assume that 
\[
\mathrm{tr}\,M_{1}=0.
\]

With the considerations as above, we will see that the mid point $W=\mathsf{Q_{0}}\cap\mathsf{C_{0}}$
can be extended to the following $11$-dimensional complete intersection
$\overline{\Sigma}$ of a quadric and a cubic:
\begin{defn}[\textbf{Extension of $W$}]
 \label{defSigmabarG4} Fix a $3$-dimensional vector space $U^{3}$.
Let $\mathrm{S}^{1,0,-1}U^{3}$ be the $8$-dimensional irreducible
component of $U^{3}\otimes(U^{3})^{*}$ as representation space of
$\SL\!(U^{3})$. We define $\overline{\Sigma}$ to be the complete
intersection 
\[
\mathsf{Q}\cap\mathsf{C}\subset\mP\left((U^{3})^{*}\oplus U^{3}\oplus\mathrm{S}^{1,0,-1}U^{3}\right)\simeq\mP^{13}
\]
 with 
\[
\text{the quadric \ensuremath{\mathsf{\mathsf{Q}}:=\{(\bm{y},\bm{x})=0\}} and the cubic \ensuremath{\mathsf{\mathsf{C}}:=\{\langle\bm{y},M,\bm{x}\rangle=0\},}}
\]
where 
\begin{itemize}
\item $\bm{x}\in U^{3}$, $\bm{y}\in(U^{3})^{*}$, and $M\in\mathrm{S}^{1,0,-1}U^{3}$, 
\item $(\ \ ,\ \ )$ is the dual pairing between $(U^{3})^{*}$ and $U^{3}$,
and 
\item $\langle\ \ ,\ \ ,\ \ \rangle$ is the natural tri-linear form induced
by the contraction 
\[
(U^{3})^{*}\times(U^{3}\otimes(U^{3})^{*})\times U^{3}\to\mC.
\]
\end{itemize}
\end{defn}

We set
\[
\overline{\Pi}_{1}:=\mP\left((U^{3})^{*}\oplus0\oplus\mathrm{S}^{1,0,-1}U^{3}\right)=\left\{ \bm{x}=\bm{o}\right\} ,\,\overline{\Pi}_{2}:=\mP\left(0\oplus U^{3}\oplus\mathrm{S}^{1,0,-1}U^{3}\right)=\left\{ \bm{y}=\bm{o}\right\} ,
\]
 which are $10$-planes contained in $\overline{\Sigma}.$
\begin{prop}
\label{prop:ExtMidG4} The mid point $W$ is a codimension $8$ linear
section of $\overline{\Sigma}$ such that $\overline{\Pi}_{1}\cap W=\Pi_{1}$
and $\overline{\Pi}_{2}\cap W=\Pi_{2}$. 
\end{prop}

\begin{proof}
By taking a basis of $U^{3}$, we may describe $\overline{\Sigma}$
explicitly as follows:

Let $\bm{e}_{1},\,\bm{e}_{2},\,\bm{e}_{3}$ be a basis of $U^{3}$
and $x_{1},x_{2},x_{3}$ the coordinates of $U^{3}$ associated to
this basis. Let $y_{1},y_{2},y_{3}$ be the coordinates of the dual
space $(U^{3})^{*}$ associated to the dual basis $\bm{e}_{1}^{*},\,\bm{e}_{2}^{*},\,\bm{e}_{3}^{*}$.
Set ${\empty}^{t}\!\bm{x}=\begin{pmatrix}x_{1} & x_{2} & x_{3}\end{pmatrix}$
and ${\empty}^{t}\!\bm{y}=\begin{pmatrix}y_{1} & y_{2} & y_{3}\end{pmatrix}$.
These notation will be compatible with the above. Let $z_{ij}$ ($1\leq i,j\leq3$)
be the coordinates of $U^{3}\otimes(U^{3})^{*}$ associated to the
basis $\bm{e}_{i}\otimes\bm{e}_{j}^{*}$. Then 
\begin{itemize}
\item $\mathrm{S}^{1,0,-1}U^{3}$ is the subspace $\{\sum_{i=1}^{3}z_{ii}=0\}$
of $U^{3}\otimes(U^{3})^{*}$, 
\item $\mathsf{\mathsf{Q}}=\{{\empty}^{t}\!\bm{y}\bm{x}=0\}$, 
\item $\mathsf{C}=\{{\empty}^{t}\!\bm{y}M\bm{x}=0\}$, where $M$ is the
$3\times3$ matrix with the entries $(z_{ij})$, and
\item $\overline{\Pi}_{1}=\left\{ \bm{x}=\bm{o}\right\} ,\,\overline{\Pi}_{2}=\left\{ \bm{y}=\bm{o}\right\} $.
\end{itemize}
Therefore we have the assertion by the discussion above. 
\end{proof}
By elementary calculations, we obtain the singular locus of $\overline{\Sigma}$
as follows:
\begin{prop}
\label{Prop:The-singular-locusG4} The singular locus of $\overline{\Sigma}$
is contained in $\overline{\Pi}_{1}\sqcup\overline{\Pi}_{2}$ and
is equal to 
\[
\left\{ [\bm{y},\bm{o},M]\in\overline{\Pi}_{1}\mid\rank\left(\begin{array}{c}
{\empty}^{t}\!\bm{y}\\
{\empty}^{t}\!\bm{y}M
\end{array}\right)\leq1\right\} \cup\left\{ [\bm{o},\bm{x},M]\in\overline{\Pi}_{2}\mid\rank\left(\begin{array}{c}
{\empty}^{t}\!\bm{x}\\
{\empty}^{t}\!\bm{x\,}{\empty}^{t}\!M
\end{array}\right)\leq1\right\} .
\]
 In particular, $\overline{\Sigma}$ is Gorenstein and normal.
\end{prop}

\subsubsection{\textbf{Crepant small resolution \label{subsec:g4Crep}}}

We take the coordinates of $U^{3}$ and $(U^{3})^{*}$ as in the proof
of Proposition \ref{prop:ExtMidG4}. We consider 
\[
B_{6}=\{{\empty}^{t}\!\bm{y}\bm{x}=0\}\subset\mP((U^{3})^{*})\times\mP(U^{3}).
\]
We also identify $B_{6}$ with its image by the Segre embedding 
\begin{align*}
S\colon\mP((U^{3})^{*})\times\mP(U^{3}) & \hookrightarrow\mP((U^{3})^{*}\otimes U^{3})\\{}
[\bm{y}]\times[\bm{x}] & \mapsto[\bm{y}\otimes\bm{x}].
\end{align*}
Then $B_{6}$ spans $\mP(\mathrm{S}^{-1,0,1}U^{3})$, where $\mathrm{S}^{-1,0,1}U^{3}$
is the $8$-dimensional irreducible component of $(U^{3})^{*}\otimes U^{3}$
as ${\rm SL}(U^{3})$-representation space. We denote by $p_{ij}$
the coordinate of $\mP((U^{3})^{*}\otimes U^{3})$ corresponding to
$\bm{e}_{i}^{*}\otimes\bm{e}_{j}.$ The subspace $\mathrm{S}^{-1,0,1}U^{3}$
is nothing but $\left\{ \sum_{i=1}^{3}p_{ii}=0\right\} .$ We denote
the natural projections by $p_{1}\colon B_{6}\to\mP((U^{3})^{*})$
and $p_{2}\colon B_{6}\to\mP(U^{3})$, and set $\sO_{B_{6}}(1,0):=p_{1}^{*}\sO_{\mP((U^{3})^{*})}(1)$
and $\sO_{B_{6}}(0,1):=p_{2}^{*}\sO_{\mP(U^{3})}(1)$. 
\begin{defn}
\label{Def:bundleG4} We set 
\[
\Sigma':=\mP_{B_{6}}(\sO_{B_{6}}(-1,0)\oplus\sO_{B_{6}}(0,-1)\oplus\Omega_{\mP(\mathrm{S}^{-1,0,1}U^{3})}^{1}(1)|_{B_{6}}).
\]
 Note that, by a standard computation, it follows that $-K_{\Sigma'}=9H_{\Sigma'}$. 
\end{defn}

\begin{prop}
\label{prop:graph} 

The following assertions holds:
\begin{enumerate}
\item The tautological linear system $|H_{\Sigma'}|$ defines a surjective
birational morphism $\Sigma'\to\overline{\Sigma}$, which we will
denote by $\varphi_{|H_{\Sigma'}|}$. 
\item The morphism $\varphi_{|H_{\Sigma'}|}$ is an isomorphism outside
of $\Sing\overline{\Sigma}$ (note that $\overline{\Pi}_{1}\cap\overline{\Pi}_{2}$
is contained in $\Sing\overline{{\Sigma}}$ by Proposition \ref{Prop:The-singular-locusG4}).
Moreover, the $\varphi_{|H_{\Sigma'}|}$-fiber over a point ${\rm {\mathsf{t}}\in\Sing\overline{\Sigma}}$
is
\[
\begin{cases}
\mP^{1}: & {\mathsf{t}}\not\in\Sing\overline{\Sigma}\setminus(\overline{\Pi}_{1}\cap\overline{\Pi}_{2}),\\
\text{a sextic del Pezzo surface:} & {\rm {\mathsf{t}}\in\overline{\Pi}_{1}\cap\overline{\Pi}_{2}.}
\end{cases}
\]
\item The morphism $\varphi_{|H_{\Sigma'}|}$ is a crepant small resolution.
\end{enumerate}
\end{prop}

\begin{proof}
Take a point $\mathsf{p}:=[W^{1}\otimes U^{1}]\in B_{6}$, where $W^{1}\subset(U^{3})^{*}$
and $U^{1}\subset U^{3}$ are $1$-dimensional subspaces such that
$W^{1}\subset(U^{1})^{\perp}$ with respect to the dual pairing. We
set $(W^{1}\otimes U^{1})^{\perp}=(\mathrm{S}^{-1,0,1}U^{3}/(W^{1}\otimes U^{1}))^{*}$.
The fiber of the projective bundle $\Sigma'\to B_{6}$ over ${\rm p}$
is 
\[
\mP(W^{1}\oplus U^{1}\oplus(W^{1}\otimes U^{1})^{\perp}),
\]
which is a linear subspace of $\mP((U^{3})^{*}\oplus U^{3}\oplus\mathrm{S}^{1,0,-1}U^{3})$.
By Lemma \ref{lem:ABFib} (1), the tautological linear system defines
a map $\Sigma'\to\mP((U^{3})^{*}\oplus U^{3}\oplus\mathrm{S}^{1,0,-1}U^{3})$.
By the descriptions of fibers of $\Sigma'\to B_{6}$ and the definition
of $\overline{\Sigma}$, we see that the image of this map is contained
in $\overline{\Sigma}$. 

Let $\mathsf{t}=[\bm{y}+\bm{x}+M]$ be a point of $\overline{\Sigma}$
with $\bm{y}\in(U^{3})^{*},\bm{x}\in U^{3},M\in\mathrm{S}^{1,0,-1}U^{3}$
. By Lemma \ref{lem:ABFib} (2), the fiber of $\Sigma'\to\overline{\Sigma}$
over $\mathsf{t}$ is 
\[
\left\{ \mathsf{t}\right\} \times\left\{ [W^{1}\otimes U^{1}]\mid W^{1}\subset(U^{1})^{\perp},\bm{y}\in W^{1},\bm{x}\in U^{1},M\in(W^{1}\otimes U^{1})^{\perp}\right\} .
\]
 We check the condition for $(\mathsf{t},[W^{1}\otimes U^{1}])$ to
be in the fiber of $\Sigma'\to\overline{\Sigma}$ over $\mathsf{t}$.
If $\mathsf{t}\not\in\overline{\Pi}_{1}\cup\overline{\Pi}_{2}$, then
$W^{1}$ and $U^{1}$ are uniquely determined as $W^{1}=\mC\bm{y}$
and $U^{1}=\mC\bm{x}$. Therefore the morphism $\Sigma'\to\overline{\Sigma}$
is an isomorphism outside of $\overline{\Pi}_{1}\cup\overline{\Pi}_{2}$.
In particular, the morphism $\Sigma'\to\overline{\Sigma}$ is surjective
and birational. Assume that $\mathsf{t}\in\overline{\Pi}_{1}\setminus(\overline{\Pi}_{1}\cap\overline{\Pi}_{2}),$
equivalently, $\bm{x=o}$ and $\bm{y}\not=\bm{o}$. Then $W^{1}$
is uniquely determined as $W^{1}=\mC\bm{y}$. We set $M_{\bm{y}}:=\left(\begin{array}{c}
{\empty}^{t}\!\bm{y}\\
{\empty}^{t}\!\bm{y}M
\end{array}\right).$ Note that, by Proposition \ref{Prop:The-singular-locusG4}, $\mathsf{t}\in(\Sing\overline{\Sigma})\cap\overline{{\Pi}}_{1}$
if and only if $\rank M_{\bm{y}}=1$. The condition for $U^{1}$ is
that $U^{1}\subset\left\{ \bm{z}\in U^{3}\mid M_{\bm{y}}\bm{z}=\bm{0}\right\} $.
Therefore, if $\rank M_{\bm{y}}=2$, then $U^{3}$ is uniquely determined,
and if $\rank M_{\bm{y}}=1$, $U^{1}$'s are parameterized by $\mP((W^{1})^{\perp})\simeq\mP^{1}$.
From this, the description of the fiber $\Sigma'\to\overline{\Sigma}$
over $\mathsf{t}$ follows. We can also describe the fiber over $\mathsf{t}\in\overline{\Pi}_{2}\setminus(\overline{\Pi}_{1}\cap\overline{\Pi}_{2})$
in the same way. Finally, assume that $\mathsf{t}\in\overline{\Pi}_{1}\cap\overline{\Pi}_{2}$.
Since $\overline{\Pi}_{1}\cap\overline{\Pi}_{2}=\mP(\mathrm{S}^{1,0,-1}U^{3})$,
its inverse image in $\Sigma'$ is the projective subbundle 
\[
S_{\Sigma'}:=\mP_{B_{6}}(\Omega_{\mP(\mathrm{S}^{-1,0,1}U^{3})}^{1}(1)|_{B_{6}}).
\]
Note that $S_{\Sigma'}\subset\mP(\mathrm{S}^{1,0,-1}U^{3})\times B_{6}$,
and $S_{\Sigma'}\to\mP(\mathrm{S}^{1,0,-1}U^{3})$ is the universal
family of hyperplane sections of $B_{6}$, which is a fibration of
sextic del Pezzo surfaces. Hence the description of the fiber of $\Sigma'\to\overline{\Sigma}$
over $\mathsf{t}\in\overline{\Pi}_{1}\cap\overline{\Pi}_{2}$ follows. 

The assertion (3) follows from (2) since it holds that $-K_{\Sigma'}=9H_{\Sigma'}$.
\end{proof}
\begin{rem}
\label{rem:flopdiff} It is possible to construct the flop for $\Sigma'\to\overline{\Sigma}$
but the construction is slightly involved (and produce singularities)
since the flopping contraction $\Sigma'\to\overline{\Sigma}$ has
jumping fibers as in Proposition \ref{prop:graph} (2). We will see
that the construction as in Proposition \ref{Prop:SigmaflopG46} is
very close to and is easier than the construction of the flop for
$\Sigma'\to\overline{\Sigma}$.
\end{rem}

\subsection{Genus 6 }

In this subsection, we use the notation as in the subsection \ref{subsec:g6setup}.

\subsubsection{\textbf{Extending the mid point in the case of $\mathsf{\mathtt{Q}}$-type}
\label{subsec:ExtQuadric-type}}
\begin{defn}[\textbf{Extension of $W$}]
\label{Def:SigmaBarG6Q}  We denote by $Q_{\mathtt{q}}$ the quadric
in the projective space $\mP(V'\oplus U^{5}\oplus(U^{5})^{*})$ defined
by the dual pairing $U^{5}\times(U^{5})^{*}\to\mC$. We set 
\[
\overline{\Sigma}:=\left(A_{\mathtt{Q}}*\mP\left((U^{5})^{*}\right)\right)\cap Q_{\mathtt{q}},
\]
and 
\[
\overline{\Pi}:=\mP(V'\oplus0\oplus(U^{5})^{*})\subset\mP(V'\oplus U^{5}\oplus(U^{5})^{*}).
\]
Note that $\overline{\Pi}\simeq\mP^{8}$ and $\overline{\Pi}\subset\overline{\Sigma}$. 
\end{defn}

\begin{prop}
\label{prop:ExtMidG6Q} The pair $(W,\Pi_{0})$ is projectively equivalent
to the pair of a linear section $W'$ of $\overline{\Sigma}$ and
the $2$-plane $\overline{\Pi}\cap W'$. 
\end{prop}

\begin{proof}
We take coordinates $x_{1},\dots,x_{4}$ of $V'$ and $y_{1},\dots,y_{5}$
of $U^{5}$ respectively. Recall that $A_{\mathtt{Q}}={\rm G}(2,V)\cap\mP(V'\oplus U^{5})$,
and $W$ is a quadric section of $A_{\mathtt{Q}}\cap\mP(U^{8})$ with
$U^{8}\subset V'\oplus U^{5}$. We may assume that $\Pi_{0}={\Pi}\cap\{x_{1}=0\}=\{x_{1}=y_{1}=\cdots=y_{5}=0\}\subset\mP(V'\oplus U^{5})$.
Then we may write $U^{8}=(V'\oplus U^{5})\cap\left\{ l(x,y)=0\right\} $
and $W=A_{\mathtt{Q}}\cap\{l(x,y)=q(x,y)=0\}$, where $l(x,y)=x_{1}+l'(y)$
with a linear form $l'(y)$ and $q(x,y)$ is a quadratic form. Since
$\Pi_{0}\subset W$, we can write $q(x,y)=x_{1}m(x)+q'(x,y)$ with
a linear form $m(x)$ and a quadric form $q'(x,y)\in(y_{1},\dots,y_{5})$.
Replacing $q(x,y)$ with $q(x,y)-l(x,y)m(x)=-l'(y)m(x)+q'(x,y)$,
we may assume that ${\Pi}\subset\{q(x,y)=0\}$. Therefore we may write
\[
q(x,y)=l_{1}y_{1}+\dots+l_{5}y_{5}
\]
with linear forms $l_{1},\dots,l_{5}$.

Now we consider the projective space $\mP(V'\oplus U^{5}\oplus(U^{5})^{*})$
and the quadric $Q_{\mathtt{q}}$ as in Definition \ref{Def:SigmaBarG6Q}.
Explicitly, let $z_{1},\dots,z_{5}$ be the coordinates of $(U^{5})^{*}$
dual to $y_{1},\dots,y_{5}$. Then 
\begin{equation}
Q_{\mathtt{q}}=\{y_{1}z_{1}+\dots+y_{5}z_{5}=0\}\subset\mP(V'\oplus U^{5}\oplus(U^{5})^{*}).\label{eq:Qq}
\end{equation}
Then, by the above construction, we see that the pair $(W,\Pi_{0})$
is projectively equivalent to the pair of $W':=\overline{\Sigma}\cap\{z_{1}-l_{1}=\cdots=z_{5}-l_{5}=0,l(x,y)=0\}$
and $\overline{\Pi}\cap W'=\left\{ x_{1}=y_{1}=\cdots=y_{5}=0,z_{1}-l_{1}=\cdots=z_{5}-l_{5}=0\right\} $. 
\end{proof}
Now we use the coordinates and the equation for $A_{\mathtt{Q}}$
as in (\ref{eq:EqAQ}). We denote by $y_{23},y_{25},y_{34},y_{35},y_{45}$
be the coordinates of $(U^{5})^{*}$ dual to that of $U^{5}$. Then
$\overline{\Sigma}$ is defined by the equations of $A_{\mathtt{Q}}$
and 
\[
Q_{\mathtt{q}}=\left\{ x_{23}y_{23}+x_{25}y_{25}+x_{34}y_{34}+x_{35}y_{35}+x_{45}y_{45}=0\right\} .
\]
Using these, we obtain the following by an explicit calculation:
\begin{prop}
\label{Prop:SIngg6Q} We set 
\begin{align*}
 & {{\empty}^{t}\!\bm{x}}:=\left(\begin{array}{ccccc}
x_{23} & x_{25} & x_{34} & x_{35} & x_{45}\end{array}\right),{{\empty}^{t}\!\bm{y}}:=\left(\begin{array}{ccccc}
y_{23} & y_{25} & y_{34} & y_{35} & y_{45}\end{array}\right),\\
 & M_{\mathtt{Q}}:=\begin{pmatrix}0 & z_{2}^{2} & z_{3}^{2} & z_{2}z_{3} & z_{2}z_{4}-z_{3}z_{5}\\
-z_{2}^{2} & 0 & z_{3}z_{5}+z_{2}z_{4} & z_{2}z_{5} & -z_{5}^{2}\\
-z_{3}^{2} & -(z_{3}z_{5}+z_{2}z_{4}) & 0 & -z_{3}z_{4} & -z_{4}^{2}\\
-z_{2}z_{3} & -z_{2}z_{5} & z_{3}z_{4} & 0 & -z_{4}z_{5}\\
-(z_{2}z_{4}-z_{3}z_{5}) & z_{5}^{2} & z_{4}^{2} & z_{4}z_{5} & 0
\end{pmatrix}.
\end{align*}
 The singular locus of $\overline{\Sigma}$ is contained in $\overline{\Pi}$
and is equal to $\left\{ \bm{x}=\bm{o},M_{\mathtt{Q}}\bm{y}=\bm{o}\right\} $.
In particular, $\overline{\Sigma}$ is Gorenstein and normal.
\end{prop}

\subsubsection{\textbf{Crepant small resolution of the mid point in the case of
$\mathtt{Q}$-type \label{subsec:g6Crep}}}
\begin{defn}
\label{Def:BundleG6Q} We set 
\[
\Sigma':=\mP_{Q^{3}}(\sU|_{Q^{3}}\oplus\sO_{Q^{3}}(-1)\oplus\Omega_{\mP(U^{5})}(1)|_{Q^{3}})\to Q^{3},
\]
where $\sU$ is the rank 2 universal subbundle on ${\rm G(2,V')}$.
By a standard computation, it follows that $-K_{\Sigma'}=7H_{\Sigma'}$. 
\end{defn}

\begin{prop}
\label{prop:graphg6} 

The following assertions hold:
\end{prop}

\begin{enumerate}
\item The tautological linear system $|H_{\Sigma'}|$ defines a surjective
birational morphism $\Sigma'\to\overline{\Sigma}$, which we will
denote by $\varphi_{|H_{\Sigma'}|}$. 
\item The morphism $\varphi_{|H_{\Sigma'}|}$ is an isomorphism outside
of $\Sing\overline{\Sigma}$ (note that $\mP\left((U^{5})^{*}\right)=\mP\left(0\oplus0\oplus(U^{5})^{*}\right)$
is contained in $\Sing\overline{{\Sigma}}$ by Proposition \ref{Prop:SIngg6Q}).
Moreover, the $\varphi_{|H_{\Sigma'}|}$-fiber over a point ${\rm {\mathsf{t}}\in\Sing\overline{\Sigma}}$
is
\[
\begin{cases}
\mP^{1}: & {\mathsf{t}}\not\in\Sing\overline{\Sigma}\setminus\mP\left((U^{5})^{*}\right),\\
\text{a quadric surface:} & {\mathsf{t}}\in\mP\left((U^{5})^{*}\right).
\end{cases}
\]
\item The morphism $\varphi_{|H_{\Sigma'}|}$ is a crepant small resolution.
\end{enumerate}
\begin{proof}
Take a point $\mathsf{p}:=[\wedge^{2}W^{2}]\in Q^{3}={\rm {G}(2,V')\cap\mP}(U^{5})$,
where $W^{2}\subset V'$ is a $2$-dimensional subspace such that
$\wedge^{2}W^{2}\subset U^{5}$. We set $(\wedge^{2}W^{2})^{\perp}=(U^{5}/\wedge^{2}W^{2})^{*}$.
The fiber of $\Sigma'\to Q^{3}$ over $\mathsf{p}$ is $\mP(W^{2}\oplus\wedge^{2}W^{2}\oplus(\wedge^{2}W^{2})^{\perp})$,
which is a linear subspace of $\mP(V'\oplus U^{5}\oplus(U^{5})^{*})$.
By Lemma \ref{lem:ABFib} (1), the tautological linear system $|H_{\Sigma'}|$
defines a morphism $\Sigma'\to\mP(V'\oplus U^{5}\oplus(U^{5})^{*}).$
By the descriptions of fibers of $\Sigma'\to Q^{3}$, Lemma \ref{lem:LinAlg}
and Definition \ref{Def:SigmaBarG6Q}, we see that the image of this
map is contained in $\overline{\Sigma}$. 

Let $\mathsf{t}=[\bm{x}+\bm{y}+\bm{z}]\in\overline{\Sigma}\subset\mP(V'\oplus U^{5}\oplus(U^{5})^{*})$
be a point with $\bm{x}\in V'$, $\bm{y}\in U^{5}$ and $\bm{z}\in(U^{5})^{*}$.
By Lemma \ref{lem:ABFib} (2), the $\varphi_{|H_{\Sigma'}|}$-fiber
over $\mathsf{t}$ is 
\[
\left\{ \mathsf{t}\right\} \times\left\{ [\wedge^{2}W^{2}]\in\rG(2,V')\mid\wedge^{2}W^{2}\subset U^{5},\bm{x}\in W^{2},\bm{y}\in\wedge^{2}W^{2},\bm{z}\in(\wedge^{2}W^{2})^{\perp}\right\} .
\]
If $\bm{y}\not=\bm{o}$, namely, $\mathsf{t}\in\overline{\Sigma}\setminus\overline{\Pi}$,
then $\bm{y}$ uniquely determines the $2$-dimensional subspace $W^{2}\subset V'$
by $\mC\bm{y}=\wedge^{2}W^{2}$. Then the fiber is nonempty by Lemma
\ref{lem:LinAlg} and (\ref{eq:Qq}), and consists of one point. Therefore
the morphism $\Sigma'\to\overline{\Sigma}$ is an isomorphism outside
of $\overline{\Pi}$. In particular, the morphism $\Sigma'\to\overline{\Sigma}$
is surjective and birational. We assume that $\mathsf{t}\in\overline{\Pi}\setminus\mP\left((U^{5})^{*}\right)$,
namely, $\bm{x}\not=\bm{o}$ and $\bm{y=\bm{o}}$. The condition of
$W^{2}$ so that $\left\{ \mathsf{t}\right\} \times\left\{ [\wedge^{2}W^{2}]\right\} $
is contained in the $\varphi_{|H_{\Sigma'}|}$-fiber over $\mathsf{t}$
is that $\bm{x}\in W^{2}$, $\wedge^{2}W^{2}\subset U^{5}$ and $\bm{z}\in(\wedge^{2}W^{2})^{\perp}$.
Using this with the coordinates and the equation (\ref{eq:EqAQ})
for $A_{\mathtt{Q}}$ and the description of $\Sing\overline{\Sigma}$
as in Proposition \ref{Prop:SIngg6Q}, we see that the $\varphi_{|H_{\Sigma'}|}$-fiber
over $\mathsf{t}$ consists of one point if $\mathsf{t}\in\overline{\Pi}\setminus\Sing\overline{\Sigma}$,
or is isomorphic to $\mP^{1}$ if $\mathsf{t}\in\Sing\overline{\Sigma}\setminus\mP\left((U^{5})^{*}\right).$
Finally, assume that $\mathsf{t}\in\mP\left((U^{5})^{*}\right)$.
The inverse image in $\Sigma'$ of $\mP\left((U^{5})^{*}\right)$
is the projective subbundle 
\[
S_{\Sigma'}:=\mP_{Q^{3}}(\Omega_{\mP((U^{5})^{*})}^{1}(1)|_{Q^{3}}).
\]
Note that $S_{\Sigma'}\subset\mP\left((U^{5})^{*}\right)\times Q^{3}$,
and $S_{\Sigma'}\to\mP\left((U^{5})^{*}\right)$ is the universal
family of hyperplane sections of $Q^{3}$, which is a fibration of
quadric surfaces. Hence the description of the fiber of $\Sigma'\to\overline{\Sigma}$
over $\mathsf{t}\in\mP\left((U^{5})^{*}\right)$ follows.

The assertion (3) follows from (2) since it holds that $-K_{\Sigma'}=7H_{\Sigma'}$.
\end{proof}

\subsubsection{\textbf{Extending the mid point in the case of $\mathsf{\mathtt{C}}$-type}
\label{subsec:ExtCubic-type}}
\begin{defn}[\textbf{Extension of $W$}]
 We denote by $Q_{\mathtt{c}}$ the quadric in the projective space
$\mP(\wedge^{2}V^{3}\oplus U^{5}\oplus(U^{5})^{*})$ defined by the
dual pairing $U^{5}\times(U^{5})^{*}\to\mC$. We set 
\[
\overline{\Sigma}:=\left(A_{\mathtt{C}}*\mP\left((U^{5})^{*}\right)\right)\cap Q_{\mathtt{c}},
\]
and 
\[
\overline{\Pi}:=\mP\left(\wedge^{2}V^{3}\oplus0\oplus(U^{5})^{*}\right)\subset\mP\left(\wedge^{2}V^{3}\oplus U^{5}\oplus(U^{5})^{*}\right).
\]
Note that $\overline{\Pi}\simeq\mP^{7}$ and $\overline{\Pi}\subset\overline{\Sigma}$. 
\end{defn}

\begin{prop}
\label{prop:ExtMidG6C} The pair $(W,\Pi)$ is projectively equivalent
to the pair of a linear section $W'$ of $\overline{\Sigma}$ and
the $2$-plane $\overline{\Pi}\cap W'$. 
\end{prop}

\begin{proof}
We can show this in a similar (and simpler) way to Proposition \ref{prop:ExtMidG6Q},
hence we omit a proof. 
\end{proof}

\section{\textbf{Embedding theorem} \textbf{in the genus 4 and 6 cases \label{sec:Common-prescription-for}}}

In this section, we show Theorem \ref{thm:main1} in the genus 4 and
6 cases (Theorem \ref{thm:embg4g6}). To show the theorem in a unified
way, we proceed in the the following two subsections \ref{subsec:Basic-set-up}
and \ref{subsec:Prescription} under a more general setting. 

\subsection{Basic set-up\label{subsec:Basic-set-up}}

Let $A$ be a Fano manifold. We denote by $f_{A}$ the Fano index
of $A$, and by $L_{A}$ the ample divisor such that $-K_{A}=f_{A}L_{A}$. 

\vspace{3pt}

\noindent \textbf{Assumption 1.} We assume that $L_{A}$ is very
ample, $\dim A\geq4$, and
\begin{equation}
d:=f_{A}-(\dim A-2)>0.\label{eq:d>0}
\end{equation}

\vspace{3pt}

We embed $A$ by $|L_{A}|$ the projective space denoted by $\mP(U_{A})$.
Sometimes we also denote $L_{A}$ by $\sO_{A}(1)$. 

\vspace{3pt}

\noindent \textbf{Assumption 2.} We assume moreover that $A$ contains
mutually disjoint codimension two linear spaces $\Pi_{1}=\mP(U_{(1)}),\dots,\Pi_{l}=\mP(U_{(l)})$,
where $U_{(1)},\dots,U_{(l)}$ are linear subspaces of $U_{A}$ of
dimension $\dim A-1$. 

\vspace{3pt}

We set 
\[
\Pi:=\Pi_{1}\sqcup\cdots\sqcup\Pi_{l}.
\]

Let $a\colon\widehat{A}\to A$ be the blow-up along $\Pi$ and $F_{a}$
the $a$-exceptional divisor, which consists of $l$ connected components. 

\vspace{3pt}

\noindent \textbf{Assumption 3.} We further assume that there exists
a morphism $b'\colon\widehat{A}\to\mP(U_{B})$ to a projective space
$\mP(U_{B})\simeq\mP^{N}$ such that 
\[
{b'}^{*}\sO_{\mP(U_{B})}(1)=a^{*}(dL_{A})-F_{a},
\]
where $d$ is defined as in (\ref{eq:d>0}).

\vspace{3pt}

We denote by $B$ the image of $b'$ and by $b\colon\widehat{A}\to B$
the induced morphism. We also set $L_{B}:=\sO_{\mP(U_{B})}(1)|_{B}$.
Therefore we have 
\begin{equation}
b^{*}L_{B}=a^{*}(dL_{A})-F_{a},\label{eq:LBLA}
\end{equation}

\begin{rem}
We can classify the situation as above but we omit a proof since we
will just apply the construction in the section \ref{sec:Common-prescription-for}
to the situations appearing for prime $\mathbb{Q}$-Fano 3-folds of
genus 4 or 6.
\end{rem}

\subsection{Construction of the key varieties \label{subsec:Prescription} }
\begin{defn}
\label{Def:BundleAhat}We define 
\begin{equation}
\widehat{\Sigma}:=\mP_{\widehat{A}}(a^{*}\sO_{A}(-1)\oplus b^{*}(\Omega_{\mP(U_{B})}(1)|_{B})).\label{eq:SigmaHat}
\end{equation}
We denote by $\pi\colon\widehat{\Sigma}\to\widehat{A}$ the natural
projection. 
\end{defn}

The linear system $|H_{\widehat{\Sigma}}|$ defines a morphism $\varphi_{|H_{\widehat{\Sigma}}|}\colon\widehat{\Sigma}\to\overline{\Sigma}$
since $\Bs|H_{\widehat{\Sigma}}|=\emptyset$. Note that $\overline{\Sigma}\subset\mP(U_{A}\oplus(U_{B})^{*})$.

\vspace{3pt}

Summarizing the above constructions, we obtain the following diagram: 

\begin{equation}
\xymatrix{ & \widehat{\Sigma}\ar[d]^{\pi=\text{proj. bundle}}\ar[r]^{\varphi_{|H_{\widehat{\Sigma}}|}\quad\quad\quad\quad} & \overline{\Sigma}\subset\mP(U_{A}\oplus(U_{B})^{*})\\
 & \widehat{A}\ar[dl]_{a=\text{bl.up along}\,\Pi\,}\ar[dr]^{b}\\
\mP(U_{A})\supset A\quad\quad &  & \quad\quad B\subset\mP(U_{B}).
}
\label{eq:basicg6Q}
\end{equation}

Using $-K_{\widehat{A}}=a^{*}(-K_{A})-F_{a}$ and (\ref{eq:LBLA}),
we have 
\begin{align}
-K_{\widehat{\Sigma}}=(N+1)H_{\widehat{\Sigma}}-\pi^{*}(K_{\widehat{A}}+a^{*}L_{A}+b^{*}L_{B})=\label{align:KSigma}\\
(N+1)H_{\widehat{\Sigma}}+(f_{A}-d-1)\pi^{*}a^{*}L_{A} & =(N+1)H_{\widehat{\Sigma}}+(\dim A-3)\pi^{*}a^{*}L_{A},\nonumber 
\end{align}
where we also use (\ref{eq:d>0}) in the last equality. By this calculation,
we see that $\Bs|-K_{\widehat{\Sigma}}|=\emptyset$. We denote by
$\nu\colon\widehat{\Sigma}\to\widehat{\Sigma}'$ the anti-canonical
model.

The $\mathbb{Q}$-Fano variety \textbf{$\Sigma$ }will be constructed
as in the following several steps summarized in the diagram: 
\begin{equation}
\xymatrix{ & \widehat{\Sigma}\ar[dd]_{b\circ\pi}\ar[ddrr]\ar@{-->}[rr]^{{\tiny\text{Atiyah flop}}} &  & \Sigma^{+}\ar@{-->}[rr]^{\text{standard flip}}\ar[dd] &  & \widetilde{\Sigma}\ar[dd]^{\text{contr. of the str. trans. of the div.\,\ensuremath{E_{\widehat{\Sigma}}} and \ensuremath{F_{\widehat{\Sigma}}}}}\ar[ddll]\\
\\
 & B &  & \overline{\Sigma} &  & \Sigma,
}
\label{eq:KeyConst}
\end{equation}
where the divisors $E_{\widehat{\Sigma}}$ and $F_{\widehat{\Sigma}}$
are defined in the sequel. We will finally achieve the construction
of $\Sigma$ in Theorem \ref{thm:common}.
\begin{rem}
This is not a Sarkisov link since the relative Picard number of $b\circ\pi$
is greater than or equal to $3$ in any case. We will, however, explain
later that this is close to an extension of the Sarkisov link for
a $\mQ$-Fano 3-fold of genus 4 or 6.
\end{rem}

\vspace{5pt}

\noindent \textbf{Two divisors $E_{\widehat{\Sigma}}$ and $F_{\widehat{\Sigma}}$
on $\widehat{\Sigma}$.} 

\vspace{5pt}

\noindent \indent We define the following two divisors $E_{\widehat{\Sigma}}$
and $F_{\widehat{\Sigma}}$ on $\widehat{\Sigma}$: 
\[
E_{\widehat{\Sigma}}:=\mP_{\widehat{A}}(0\oplus b^{*}\Omega_{\mP(U_{B})}(1)|_{B}),\,F_{\widehat{\Sigma}}:=\pi^{*}F_{a}.
\]

We remark that 
\begin{equation}
\text{\ensuremath{E_{\widehat{\Sigma}}\sim H_{\widehat{\Sigma}}-\pi^{*}a^{*}L_{A}},\ \text{and} \ensuremath{F_{\widehat{\Sigma}}\sim\pi^{*}(a^{*}(dL_{A})-b^{*}L_{B})}},\label{eq:EF}
\end{equation}
where the former is a standard equation as for projective subbundle
and the latter follows from (\ref{eq:LBLA}). 

In the step by step construction of the birational map from $\widehat{\Sigma}$
to $\Sigma$ in the sequel, it is useful to describe how these two
divisors $E_{\widehat{\Sigma}}$ and $F_{\widehat{\Sigma}}$ on $\widehat{\Sigma}$
are transformed by birational maps.

\vspace{5pt}

\noindent \textbf{Flop $\boldmath{\widehat{\Sigma}\dashrightarrow\Sigma^{+}}$.}
\begin{prop}
\label{Prop:SigmaflopG46} The following assertions hold:

\vspace{3pt}

\noindent $(1)$ The anti-canonical model $\nu\colon\widehat{\Sigma}\to\widehat{\Sigma}'$
is defined over $\overline{\Sigma}$, and is a flopping contraction
of Atiyah type. The divisor $\pi^{*}b^{*}L_{B}$ is relatively ample
for the flopping contraction.

\vspace{3pt}

\noindent $(2)$ The anti-canonical model $E_{\widehat{\Sigma}}\to E_{\widehat{\Sigma}}'$
is the restriction of $\nu$ and is also a flopping contraction of
Atiyah type. It is defined over $\mP(0\oplus(U_{B})^{*})$. 
\end{prop}

\begin{proof}
$(1)$. Let $l\subset\widehat{\Sigma}$ be an irreducible $\nu$-exceptional
curve (the existence of such an $l$ will be verified below). By (\ref{align:KSigma}),
we have $H_{\widehat{\Sigma}}\cdot l=\pi^{*}a^{*}L_{A}\cdot l=0$,
where we also use the assumption $\dim A\geq4$. By the former condition,
$\nu$ is defined over $\overline{\Sigma}$. Since $\pi\colon\widehat{\Sigma}\to\widehat{A}$
is a projective bundle and $H_{\widehat{\Sigma}}=\sO(1)$ in a fiber,
$l$ is not contracted by $\pi$. By the latter condition, $l$ is
contracted also by $a\circ\pi\colon\widehat{\Sigma}\to A$. Therefore
the image $\gamma$ of $l$ on $\widehat{A}$ is an exceptional curve
of the blow-up $a\colon\widehat{A}\to A$ along $\Pi$. Thus $\gamma\simeq\mP^{1}$
and is mapped isomorphically to a line $\gamma'$ on $B$ by the equation
(\ref{eq:LBLA}). Then, by (\ref{eq:SigmaHat}), the restriction $\widehat{\Sigma}_{\gamma}$
of $\widehat{\Sigma}$ over $\gamma$ is isomorphic to $\mP_{\mP^{1}}(\sO_{\mP^{1}}^{\oplus N}\oplus\sO_{\mP^{1}}(-1))$.
Since $l$ is contained in $\widehat{\Sigma}_{\gamma}$ and the map
defined by $|H_{\widehat{\Sigma}}|_{\widehat{\Sigma}_{\gamma}}|$
is the blow-up of $\mP^{N+1}$ along a $(N-1)$-plane, we see that
$l$ is an exceptional curve of this blow-up and hence $l\simeq\mP^{1}$
(now the existence of the curve $l$ has been verified). 

To show the first assertion, it suffices to show $\sN_{l/\widehat{\Sigma}}\simeq\sO_{\mP^{1}}^{\oplus\dim A+N-3}\oplus\sO_{\mP^{1}}(-1)^{\oplus2}$.
Note that the normal bundle $\sN_{\gamma/\widehat{A}}$ is $\sO_{\mP^{1}}(-1)\oplus\sO_{\mP^{1}}^{\oplus\dim A-2}$
since $a\colon\widehat{A}\to A$ is the blow-up along $\Pi$ and $\gamma$
is one of its fiber. Therefore the restriction to $l$ of $\sN_{\widehat{\Sigma}_{\gamma}/\widehat{\Sigma}}$
is also $\sO_{\mP^{1}}(-1)\oplus\sO_{\mP^{1}}^{\oplus\dim A-2}$.
Moreover, $\sN_{l/\widehat{\Sigma}_{\gamma}}\simeq\sO_{\mP^{1}}(-1)\oplus\sO_{\mP^{1}}^{\oplus N-1}$
since $l$ is a fiber of the blow-up $\widehat{\Sigma}_{\gamma}\to\mP^{N+1}$
along a $(N-1)$-plane. Therefore, by the normal bundle sequence $0\to\sN_{l/\widehat{\Sigma}_{\gamma}}\to\sN_{l/\widehat{\Sigma}}\to\sN_{\widehat{\Sigma}_{\gamma}/\widehat{\Sigma}}|_{l}\to0$,
we see that $\sN_{l/\widehat{\Sigma}}\simeq\sO_{\mP^{1}}^{\oplus\dim A+N-3}\oplus\sO_{\mP^{1}}(-1)^{\oplus2}$
as desired. 

The divisor $L_{B}$ is relatively ample for the flopping contraction
since the image of a flopping curve on $B$ is a line as we have seen
above.

\vspace{3pt}

\noindent $(2).$ Note that, by (\ref{align:KSigma}) and (\ref{eq:EF}),
we have 
\begin{align}
-K_{E_{\widehat{\Sigma}}}=\{(N+1)H_{\widehat{\Sigma}}+(\dim A-3\pi^{*}a^{*}L_{A}\}|_{E_{\widehat{\Sigma}}}-(H_{\widehat{\Sigma}}-\pi^{*}a^{*}L_{A})|_{E_{\widehat{\Sigma}}}=\\
N(H_{E_{\widehat{\Sigma}}})+(\dim A-2)\pi^{*}a^{*}L_{A}.\nonumber 
\end{align}
Thus $\Bs|-K_{E_{\widehat{\Sigma}}}|=\emptyset$. Using these, we
can show (2) in a similar way to (1).
\end{proof}
Let $\widehat{\Sigma}\dashrightarrow\Sigma^{+}$ be the flop for this
flopping contraction $\nu$. It is well-known that the flop can be
constructed by the blow-up along the $\nu$-exceptional locus and
the blow-down of the exceptional divisor along the other direction.
Let $E_{\Sigma^{+}}$ be the strict transform on $\Sigma^{+}$ of
$E_{\widehat{\Sigma}}$. By the construction of the flop, we see that
the restriction $E_{\widehat{\Sigma}}\dashrightarrow E_{\Sigma^{+}}$
of $\widehat{\Sigma}\dashrightarrow\Sigma^{+}$ is also the flop.

\vspace{5pt}

We call a positive dimensional fiber of $\varphi_{|H_{\widehat{\Sigma}}|}\colon\widehat{\Sigma}\to\overline{\Sigma}$
a \textit{$\varphi_{|H_{\widehat{\Sigma}}|}$-exceptional curve,}
and the union of $\varphi_{|H_{\widehat{\Sigma}}|}$-exceptional curves
\textit{the $\varphi_{|H_{\widehat{\Sigma}}|}$-exceptional locus.}
We can identify the $\varphi_{|H_{\widehat{\Sigma}}|}$-exceptional
locus as follows:
\begin{lem}
\noindent \label{lem:excepflopE G46} The $\varphi_{|H_{\widehat{\Sigma}}|}$-exceptional
locus is the union of the flopping locus for $\widehat{\Sigma}\dashrightarrow\Sigma^{+}$
and the divisor $E_{\widehat{\Sigma}}$. The flopping locus is contained
in $F_{\widehat{\Sigma}}$. The $\varphi_{|H_{\widehat{\Sigma}}|}$-image
of $E_{\widehat{\Sigma}}=\mP_{\widehat{A}}(0\oplus b^{*}\Omega_{\mP(U_{B})}(1)|_{B})$
is $\mP(0\oplus(U_{B})^{*})$ and the $\varphi_{|H_{\widehat{\Sigma}}|}$-inverse
image of $\mP(0\oplus(U_{B})^{*})$ coincides with $E_{\widehat{\Sigma}}$.
\end{lem}

\begin{proof}
Let $l$ be a $\varphi_{|H_{\widehat{\Sigma}}|}$-exceptional curve.
Note that $H_{\widehat{\Sigma}}\cdot l=0$. If $\pi^{*}a^{*}L_{A}\cdot l>0$,
then by (\ref{eq:EF}), we have $E_{\widehat{\Sigma}}\cdot l<0$ and
hence $l\subset E_{\widehat{\Sigma}}$. If $\pi^{*}a^{*}L_{A}\cdot l=0$,
then, by the proof of Proposition \ref{Prop:SigmaflopG46}, $l$ is
an exceptional curve of the anti-canonical model $\widehat{\Sigma}\to\widehat{\Sigma}'$,
namely, a flopping curve. Therefore the $\varphi_{|H_{\widehat{\Sigma}}|}$-exceptional
locus is contained in the union of the flopping locus for $\widehat{\Sigma}\dashrightarrow\Sigma^{+}$
and the divisor $E_{\widehat{\Sigma}}$. Since the restriction of
$\varphi_{|H_{\widehat{\Sigma}}|}$ to $E_{\widehat{\Sigma}}$ is
$\mP_{\widehat{A}}(0\oplus b^{*}\Omega_{\mP(U_{B})}(1)|_{B})\to\mP(0\oplus(U_{B})^{*})$,
$E_{\widehat{\Sigma}}$ is contained in $\varphi_{|H_{\widehat{\Sigma}}|}$-exceptional
locus. Thus the first assertion follows. By the second assertion of
Proposition \ref{Prop:SigmaflopG46} (1) and (\ref{eq:EF}), a flopping
curve is negative for $F_{\widehat{\Sigma}}$, hence is contained
in $F_{\widehat{\Sigma}}$. Therefore the second assertion follows.
The final assertion obviously holds.
\end{proof}
We denote by $a'$ and $b'$ the restriction to $F_{a}$ of $a$ and
$b$ respectively. By (\ref{eq:SigmaHat}), we have $F_{\widehat{\Sigma}}=\mP_{F_{a}}({a'}^{*}\sO_{\Pi}(-1)\oplus{b'}^{*}\Omega_{\mP(U_{B})}(1)|_{B})$.
\[
\xymatrix{ & F_{\widehat{\Sigma}}\ar[d]_{\pi_{F}}\\
 & F_{a}\ar[dl]_{a'}\ar[dr]^{b'}\\
\Pi &  & B.
}
\]

\begin{lem}
\noindent \label{Lem:contF} The following assertions hold:

\noindent \vspace{3pt}

\noindent \noindent $(1)$ Let $\nu_{F}$ be the restriction to $F_{\widehat{\Sigma}}$
of the flopping contraction $\nu\,($note that $\nu_{F}$ is the contraction
over $\Pi$ by the proof of Proposition \ref{Prop:SigmaflopG46}).
The $\nu_{F}$-image of $F_{\widehat{\Sigma}}$ is isomorphic to $\mP_{\Pi}(\sO_{\Pi}(-1)\oplus(U_{B})^{*}\otimes\sO_{\Pi})=\sqcup_{i=1}^{l}\mP_{\Pi_{i}}(\sO_{\Pi_{i}}(-1)\oplus(U_{B})^{*}\otimes\sO_{\Pi_{i}})$.
The image of $\mP_{\Pi_{i}}(\sO_{\Pi_{i}}(-1)\oplus(U_{B})^{*}\otimes\sO_{\Pi_{i}})$
on $\overline{\Sigma}$ is $\mP(U_{(i)}\oplus(U_{B})^{*}),$ which
we denote by $\overline{\Pi}_{i}$.

\noindent \vspace{3pt}

\noindent \noindent $(2)$ Let $F_{\Sigma^{+}}$ be the strict transform
on $\Sigma^{+}$ of $F_{\widehat{\Sigma}}$. The restriction $F_{\widehat{\Sigma}}\dashrightarrow F_{\Sigma^{+}}$
of the flop $\widehat{\Sigma}\dashrightarrow\Sigma^{+}$ is identified
with the contraction $F_{\widehat{\Sigma}}\to\mP_{\Pi}(\sO_{\Pi}(-1)\oplus(U_{B})^{*}\otimes\sO_{\Pi})$. 
\end{lem}

\begin{proof}
$(1)$. As we have seen in the proof of Proposition \ref{Prop:SigmaflopG46},
a fiber of $F_{\widehat{\Sigma}}\to\Pi$ ($\widehat{\Sigma}_{\gamma}$
in the proof) is $\mP_{\mP^{1}}(\sO_{\mP^{1}}^{\oplus N}\oplus\sO_{\mP^{1}}(-1))$.
By the proof of Proposition \ref{Prop:SigmaflopG46} again, the restriction
of $\nu_{F}$ to a fiber of $F\to\Pi$ is the blow-up of $\mP^{N}$
along a $(N-1)$-plane. Hence $\nu_{F}(F_{\widehat{\Sigma}})$ is
the $\mP^{N}$-bundle $\mP_{\Pi}(\sE_{F_{\widehat{\Sigma}}})$ with
$\sE_{F_{\widehat{\Sigma}}}:=\{{a'}_{*}{\pi_{F}}_{*}\sO_{F_{\widehat{\Sigma}}}(H_{F_{\widehat{\Sigma}}}))\}^{*}$.
We have 
\[
{a'}_{*}{\pi_{F}}_{*}\sO_{F_{\widehat{\Sigma}}}(H_{F_{\widehat{\Sigma}}})\simeq{a'}_{*}({a'}^{*}\sO_{\Pi}(1)\oplus{b'}^{*}T_{\mP(U_{B})}(-1)|_{B})\simeq\sO_{\Pi}(1)\oplus{a'}_{*}{b'}^{*}T_{\mP(U_{B})}(-1)|_{B}.
\]
To compute ${a'}_{*}{b'}^{*}T_{\mP(U_{B})}(-1)|_{B}$, we consider
the restriction to $B$ of the Euler sequence of $\mP(U_{B})$: 
\[
0\to\sO_{B}(-1)\to U_{B}\otimes\sO_{B}\to T_{\mP(U_{B})}(-1)|_{B}\to0.
\]
Since ${a'}_{*}{b'}^{*}\sO_{B}(-1)=R^{1}{a'}_{*}{b'}^{*}\sO_{B}(-1)=0$
by (\ref{eq:LBLA}), we have ${a'}_{*}{b'}^{*}T_{\mP(U_{B})}(-1)|_{B}\simeq U_{B}\otimes\sO_{\Pi}$.
Therefore we have $\sE_{F_{\widehat{\Sigma}}}\simeq\sO_{\Pi}(-1)\oplus(U_{B})^{*}\otimes\sO_{\Pi}$
as desired. The final assertion obviously holds.

\vspace{3pt}

\noindent $(2)$. The assertion follows from the explicit construction
of the flop of Atiyah type.
\end{proof}
\begin{lem}
\label{Lem:EandFint} 

The following assertions hold:
\end{lem}

\begin{enumerate}
\item $E_{\Sigma^{+}}\cap F_{\Sigma^{+}}=\mP_{\Pi}((U_{B})^{*}\otimes\sO_{\Pi})\simeq\Pi\times\mP((U_{B})^{*})$. 
\item The exceptional locus of ${\Sigma}^{+}\to\overline{\Sigma}$ is the
union of the divisor $E_{\Sigma^{+}}$ and the flopped locus. 
\end{enumerate}
\begin{proof}
(1). By definition, the intersection $E_{\widehat{\Sigma}}\cap F_{\widehat{\Sigma}}$
is equal to $\mP_{F_{a}}({b'}^{*}(\Omega_{\mP(U_{B})}^{1}(1)|_{B}))$.
By the contraction $F_{\widehat{\Sigma}}\to\mP_{\Pi}(\sO_{\Pi}(-1)\oplus(U_{B})^{*}\otimes\sO_{\Pi})\simeq F_{\Sigma^{+}}$,
$E_{\widehat{\Sigma}}\cap F_{\widehat{\Sigma}}$ is mapped onto $\mP_{\Pi}((U_{B})^{*}\otimes\sO_{\Pi})$.
By the explicit construction of the flop $\widehat{\Sigma}\dashrightarrow\Sigma^{+}$,
we see that $E_{\Sigma^{+}}\cap F_{\Sigma^{+}}=\mP_{\Pi}((U_{B})^{*}\otimes\sO_{\Pi})$.

\noindent (2). The assertion follows from Lemma \ref{lem:excepflopE G46}.
\end{proof}
In the following steps, we separate $E_{\Sigma^{+}}$ and $F_{\Sigma^{+}}$
by a flip, and finally contract their strict transforms.

\vspace{5pt}

\noindent\textbf{ Flip $\boldmath{{\Sigma}^{+}\dashrightarrow\widetilde{\Sigma}}$.}

Let $H_{\Sigma^{+}}$ and $L_{A}^{+}$ be the strict transforms on
$\Sigma^{+}$ of $H_{\widehat{\Sigma}}$ and $\pi^{*}a^{*}L_{A}$
respectively.
\begin{prop}
\label{Prop:flipping G46} 
\end{prop}

\begin{enumerate}
\item Let $\Gamma$ be a fiber of $E_{\Sigma^{+}}\cap F_{\Sigma^{+}}\to\mP((U_{B})^{*})$,
which is a copy of $\Pi$ by Lemma \ref{Lem:EandFint} (1). It holds
that
\[
\sN_{\Gamma/\Sigma^{+}}=\sO_{\mP^{\dim A-2}}(-1)^{\oplus2}\oplus\sO_{\mP^{\dim A-2}}^{\oplus N}.
\]
\item There exists a small contraction $\Sigma^{+}\to(\Sigma^{+})'$ contracting
$E_{\Sigma^{+}}\cap F_{\Sigma^{+}}\simeq\Pi\times\mP((U_{B})^{*})$
onto $\mP((U_{B})^{*})$. 
\end{enumerate}
\begin{proof}
(1). We set $G=E_{\Sigma^{+}}\cap F_{\Sigma^{+}}$. To determine the
normal bundle $\sN_{\Gamma/\Sigma^{+}}$, let us consider the normal
bundle sequence $0\to\sN_{\Gamma/F_{\Sigma^{+}}}\to\sN_{\Gamma/\Sigma^{+}}\to{\sN_{F_{\Sigma^{+}}/\Sigma^{+}}}|_{\Gamma}\to0.$
Since $G=\mP_{\Pi}((U_{B})^{*}\otimes\sO_{\Pi})\subset\mP_{\Pi}(\sO_{\Pi}(-1)\oplus(U_{B})^{*}\otimes\sO_{\Pi})=F_{\Sigma^{+}}$,
we have $\sN_{G/F_{\Sigma^{+}}}=(H_{F_{\Sigma^{+}}}-L_{\Pi})|_{G}$,
where $L_{\Pi}$ is the pull-back of $\sO_{\Pi}(1)$. Therefore we
have $\sN_{\Gamma/F_{\Sigma^{+}}}\simeq\sO_{\mP^{\dim A-2}}(-1)\oplus\sO_{\mP^{\dim A-2}}^{\oplus N}$.
By (\ref{align:KSigma}), we have $-K_{\Sigma^{+}}=(N+1)H_{\Sigma^{+}}+(\dim A-3)L_{A}^{+}$.
Since $H_{\Sigma^{+}}|_{\Gamma}=0$ and $L_{A}^{+}|_{\Gamma}=\sO_{\mP^{\dim A-2}}(1)$,
we have $-K_{\Sigma^{+}}|_{\Gamma}=\sO_{\mP^{\dim A-2}}(\dim A-3)$.
Therefore $\deg\sN_{\Gamma/\Sigma^{+}}=\dim A-3-(\dim A-1)=-2$ and
hence by the above normal bundle sequence, $F_{\Sigma^{+}}|_{\Gamma}=\sO_{\mP^{\dim A-2}}(-1)$
and then we have $\sN_{\Gamma/\Sigma^{+}}=\sO_{\mP^{\dim A-2}}(-1)^{\oplus2}\oplus\sO_{\mP^{\dim A-2}}^{\oplus N}$.

\noindent (2). We show that $L_{A}^{+}+F_{\Sigma^{+}}$ is nef over
$\overline{\Sigma}$ and numerically trivial only for fibers of $G\to\mP((U_{B})^{*})$.
Assume that $(L_{A}^{+}+F_{\Sigma^{+}})\cdot\gamma\leq0$ for an exceptional
curve $\gamma$ for $\Sigma^{+}\to\overline{\Sigma}$. It is enough
to show that $\gamma$ is contained in a fiber of $G\to\mP((U_{B})^{*})$,
and $(L_{A}^{+}+F_{\Sigma^{+}})\cdot\gamma=0$. By Lemma \ref{Lem:EandFint}
(2), $\gamma$ is a flopped curve or is contained in $E_{\Sigma^{+}}.$
Assume that $\gamma$ is a flopped curve. Let $\gamma'\subset\widehat{\Sigma}$
be the corresponding flopping curve. Since $\pi^{*}a^{*}L_{A}\cdot\gamma'=0$
on $\widehat{\Sigma}$, we have $L_{A}^{+}\cdot\gamma=0$ on $\Sigma^{+}$.
By the proof of Proposition \ref{Prop:SigmaflopG46}, we have $F_{\widehat{\Sigma}}\cdot\gamma'=-1$,
hence we have $F_{\Sigma^{+}}\cdot\gamma=1$ by a property of the
flop of Atiyah type. Therefore $(L_{A}^{+}+F_{\Sigma^{+}})\cdot\gamma=1>0$,
a contradiction. Thus we have $\gamma\subset E_{\Sigma^{+}}$. If
$F_{\Sigma^{+}}\cdot\gamma<0$, then $\gamma\subset F_{\Sigma^{+}}$,
hence $\gamma\subset F_{\Sigma^{+}}\cap E_{\Sigma^{+}}=\Lambda$.
Since $\gamma$ is exceptional over $\overline{\Sigma}$, $\gamma$
must be contained in a fiber of $G\to\mP((U_{B})^{*})$. To compute
$(L_{A}^{+}+F_{\Sigma^{+}})\cdot\gamma$, we may assume that $\gamma$
is a line. Then we have $L_{A}^{+}\cdot\gamma=1$, and $F_{\Sigma^{+}}\cdot\gamma=-1$
by the proof of (1). Therefore $(L_{A}^{+}+F_{\Sigma^{+}})\cdot\gamma=0$
as desired. Since we are already done if $F_{\Sigma^{+}}\cdot\gamma<0$,
we may assume that $F_{\Sigma^{+}}\cdot\gamma\geq0$ in the sequel.
Then, since $L_{A}^{+}$ is nef, we have $(L_{A}^{+}+F_{\Sigma^{+}})\cdot\gamma\geq0$,
hence $L_{A}^{+}\cdot\gamma=F_{\Sigma^{+}}\cdot\gamma=0$ by the assumption
that $(L_{A}^{+}+F_{\Sigma^{+}})\cdot\gamma\leq0$. By $F_{\Sigma^{+}}\cdot\gamma=0$,
$\gamma$ cannot be a flopped curve. Therefore its strict transform
$\gamma'$ on $\widehat{\Sigma}$ satisfies $\pi^{*}a^{*}L_{A}\cdot\gamma'=0$
and $H_{\widehat{\Sigma}}\cdot\gamma'=0$. However, this implies that
$\gamma'$ is a flopping curve by the proof of Proposition \ref{Prop:SigmaflopG46},
a contradiction. Now we have shown that $L_{A}^{+}+F_{\Sigma^{+}}$
is nef over $\overline{\Sigma}$ and numerically trivial only for
fibers of $G\to\mP((U_{B})^{*})$. 

Note that $-K_{\Sigma^{+}}$ is nef and big since so is $-K_{\widehat{\Sigma}}$
by Proposition \ref{Prop:SigmaflopG46} and $\widehat{\Sigma}\dashrightarrow\Sigma^{+}$
is a flop. Therefore, $L_{A}^{+}+F_{\Sigma^{+}}$ is semiample by
the Kawamata-Shokurov base point free theorem (cf.\cite{KMM}). Thus
the contraction over $\overline{\Sigma}$ defined by $L_{A}^{+}+F_{\Sigma^{+}}$
is the desired one. 
\end{proof}
By Proposition \ref{Prop:flipping G46} (1), the contraction $\Sigma^{+}\to\overline{\Sigma}^{+}$
is of flipping type, and the flip can be constructed by the blow-up
along $G$ and the blow-down of the exceptional divisor along the
other direction (this is a so called a family of standard flips \cite{Ka1}).
Let $\Sigma^{+}\dashrightarrow\widetilde{\Sigma}$ be the flip. By
Proposition \ref{Prop:flipping G46} (1) again, the flipped locus
is a $\mP^{1}$-bundle over $\mP((U_{B})^{*}).$ We denote by $E_{\widetilde{\Sigma}}$,
$F_{\widetilde{\Sigma}}$ and $H_{\widetilde{\Sigma}}$ be the strict
transforms on $\widetilde{\Sigma}$ of $E_{\Sigma^{+}}$, $F_{\Sigma^{+}}$
and $H_{\widehat{\Sigma}}$ respectively. 

\vspace{5pt}

\noindent \textbf{Contracting $E_{\widetilde{\Sigma}}$ and $F_{\widetilde{\Sigma}}$.} 

\vspace{5pt}

\noindent \noindent By the constructions of the flop $\widehat{\Sigma}\dashrightarrow\Sigma^{+}$
and the flip $\Sigma^{+}\dashrightarrow\widetilde{\Sigma}$, and the
description of $E_{\Sigma^{+}}\cap F_{\Sigma^{+}}$ as in Lemma \ref{Lem:EandFint}
(1), we see that $E_{\widetilde{\Sigma}}\cap F_{\widetilde{\Sigma}}=\emptyset$.

By the construction of the flip, we see that the restriction $F_{\Sigma^{+}}\dashrightarrow F_{\widetilde{\Sigma}}$
of the flip $\Sigma^{+}\dashrightarrow\widetilde{\Sigma}$ is the
contraction $F_{\Sigma^{+}}\to\sqcup_{i=1}^{l}\mP(U_{(i)}\oplus U_{B}^{*})$,
where we recall that $U_{(i)}$ are defined in the beginning of the
subsection \ref{subsec:Basic-set-up}. Thus $F_{\widetilde{\Sigma}}$
is the disjoint union of $\widetilde{F}_{i}:=\mP(U_{(i)}\oplus U_{B}^{*})\simeq\mP^{\dim A+N-1}$
$(i=1,\dots,l)$.
\begin{lem}
\label{lem:F|F G46} The normal bundle $\sN_{\widetilde{F}_{i}/\widetilde{\Sigma}}$
is $\sO_{\mP^{\dim A+N-1}}(-2)$ for $i=1,\dots,l$.
\end{lem}

\begin{proof}
Since $E_{\widetilde{\Sigma}}\cap F_{\widetilde{\Sigma}}=\emptyset$,
we see that the restriction to $F_{\widetilde{\Sigma}}$ of the strict
transform on $\widetilde{\Sigma}$ of $\pi^{*}a^{*}L_{A}$ is linearly
equivalent to $H_{\widetilde{\Sigma}}|_{F_{\widetilde{\Sigma}}}$
by (\ref{eq:EF}). Therefore, by (\ref{align:KSigma}), we see that
$-K_{\widetilde{\Sigma}}|_{\widetilde{F}_{i}}=\sO_{\mP^{\dim A+N-1}}(N+\dim A-2)$.
Since $-K_{\widetilde{F}_{i}}=\sO_{\mP^{\dim A+N-1}}(\dim A+N)$,
we have $\sN_{\widetilde{F}_{i}/\widetilde{\Sigma}}=\sO_{\mP^{\dim A+N-1}}(-2)$
as desired. 
\end{proof}
\begin{lem}
\label{Lem:2H+F0} $2H_{\widetilde{\Sigma}}+F_{\widetilde{\Sigma}}$
is semiample. 
\end{lem}

\begin{proof}
We show that $2H_{\widetilde{\Sigma}}+F_{\widetilde{\Sigma}}$ is
nef. Assume that $(2H_{\widetilde{\Sigma}}+F_{\widetilde{\Sigma}})\cdot\gamma<0$
for an irreducible curve $\gamma$. Then $F_{\widetilde{\Sigma}}\cdot\gamma<0$
since $H_{\widetilde{\Sigma}}$ is nef, and hence $\gamma\subset F_{\widetilde{\Sigma}}$.
Since $F_{\widetilde{\Sigma}}|_{\widetilde{F}_{i}}=\sO_{\mP^{\dim A+N-1}}(-2)$
and $H_{\widetilde{\Sigma}}|_{\widetilde{F}_{i}}=\sO_{\mP^{\dim A+N-1}}(1)$,
we have $(2H_{\widetilde{\Sigma}}+F_{\widetilde{\Sigma}})\cdot\gamma=0$,
a contradiction. Therefore $2H_{\widetilde{\Sigma}}+F_{\widetilde{\Sigma}}$
is nef.

To show $2H_{\widetilde{\Sigma}}+F_{\widetilde{\Sigma}}$ is semiample,
we have only to show $m(2H_{\widetilde{\Sigma}}+F_{\widetilde{\Sigma}})-K_{\widetilde{\Sigma}}$
is nef and big for $m\gg0$ by the Kawamata-Shokurov base point free
theorem. Since $-K_{\widehat{\Sigma}}$ is nef and big, and $\widehat{\Sigma}\dashrightarrow\Sigma^{+}$
is a flop, we see that $-K_{\Sigma^{+}}$ is also nef and big. Since
$\Sigma^{+}\dashrightarrow\widetilde{\Sigma}$ is a flip, we see that
$-K_{\widetilde{\Sigma}}$ is big and is negative only for flipped
curves. Let $\gamma$ be a flipped curve. Then $-K_{\widetilde{\Sigma}}\cdot\gamma=-(\dim A+N-2)$
and $F_{\widetilde{\Sigma}}\cdot\gamma=1$ by the construction of
the flip, we have $(m(2H_{\widetilde{\Sigma}}+F_{\widetilde{\Sigma}})-K_{\widetilde{\Sigma}})\cdot\gamma>0$
for $m\gg0$. Therefore $m(2H_{\widetilde{\Sigma}}+F_{\widetilde{\Sigma}})-K_{\widetilde{\Sigma}}$
is nef for $m\gg0$. The bigness is clear since $2H_{\widetilde{\Sigma}}+F_{\widetilde{\Sigma}}$
is nef and $-K_{\widetilde{\Sigma}}$ is big. 
\end{proof}
\begin{thm}
\label{thm:common}Let $\mu\colon\widetilde{\Sigma}\to\Sigma$ be
the contraction defined by a sufficient multiple of $2H_{\widetilde{\Sigma}}+F_{\widetilde{\Sigma}}$.
We recall that $\mP(U_{B})\simeq\mP^{N}$. The following assertion
holds:

\vspace{3pt}

\noindent $(1)$ The $\mu$-exceptional locus is the union of the
two divisors $E_{\widetilde{\Sigma}}$ and $F_{\widetilde{\Sigma}}$. 

\vspace{3pt}

\noindent $(2)$ The $\mu$-image of $F_{\widetilde{\Sigma}}$ consists
of $l$ $1/2$-singularities.

\vspace{3pt}

\noindent $(3)$ The discrepancy of $E_{\widetilde{\Sigma}}$ is
$\dim A-3$ and $\mu(E_{\widetilde{\Sigma}})\simeq\mP((U_{B})^{*})$.
In particular $\Sigma$ has Gorenstein terminal singularities along
$\mP((U_{B})^{*})$.

\vspace{3pt}

\noindent $(4)$ If $\rho(A)=1$, then the $(\dim A+N)$-dimensional
variety $\Sigma$ is a $\mQ$-Fano variety with only terminal singularities
and with $\rho(\Sigma)=1$. 

\vspace{3pt}

\noindent $(5)$ The image $M_{\Sigma}$ of $H_{\Sigma}$ is a primitive
integral ample Weil divisor $M_{\Sigma}$ and it holds that $-K_{\Sigma}=(\dim A+N-2)M_{\Sigma}$. 
\end{thm}

\begin{proof}
As we have checked in the proof of Lemma \ref{Lem:2H+F0}, $2H_{\widetilde{\Sigma}}+F_{\widetilde{\Sigma}}$
is numerical trivial for any curve in $F_{\widetilde{\Sigma}}$. Thus
the image of $F_{\widetilde{\Sigma}}$ by $\widetilde{\Sigma}\to\Sigma$
consists of $l$ $1/2$-singularities by Lemma \ref{lem:F|F G46}.
Since $E_{\widetilde{\Sigma}}\cap F_{\widetilde{\Sigma}}=\emptyset$,
we have $(2H_{\widetilde{\Sigma}}+F_{\widetilde{\Sigma}})|_{E_{\widetilde{\Sigma}}}=2H_{\widetilde{\Sigma}}|_{E_{\widetilde{\Sigma}}}$.
Therefore $E_{\widetilde{\Sigma}}$ is $\mu$-exceptional and $\mu(E_{\widetilde{\Sigma}})$
is isomorphic to $\mP((U_{B})^{*})$.

We show that the $\mu$-exceptional locus is the union of $E_{\widetilde{\Sigma}}$
and $F_{\widetilde{\Sigma}}$. Assume by contradiction that $(2H_{\widetilde{\Sigma}}+F_{\widetilde{\Sigma}})\cdot\gamma=0$
for an irreducible curve $\gamma\not\subset E_{\widetilde{\Sigma}}\cup F_{\widetilde{\Sigma}}$.
Since $H_{\widetilde{\Sigma}}$ is nef and $\gamma\not\subset F_{\widetilde{\Sigma}}$,
we have $H_{\widetilde{\Sigma}}\cdot\gamma=F_{\widetilde{\Sigma}}\cdot\gamma=0$.
Then, by Lemma \ref{Lem:EandFint} (2) and the conditions that $H_{\widetilde{\Sigma}}\cdot\gamma=0$
and $\gamma\not\subset E_{\widetilde{\Sigma}}$, $\gamma$ is a flipped
curve or the strict transform of a flopped curve. If $\gamma$ is
a flipped curve, then $F_{\widetilde{\Sigma}}\cdot\gamma>0$, a contradiction.
Assume that $\gamma$ is the strict transform of a flopped curve.
If $\gamma$ is disjoint from flipping curves, then $F_{\widetilde{\Sigma}}\cdot\gamma>0$
since a flopped curve is positive for $F_{\Sigma^{+}}$, a contradiction.
Therefore $\gamma$ intersects a flipped curve. Let $\gamma'$ be
the strict transform of $\gamma$ on $\Sigma^{+}$ ($\gamma'$ is
a flopped curve). Since $\gamma$ intersects a flipped curve and the
flipping locus is $E_{\Sigma^{+}}\cap F_{\Sigma^{+}}$, we see that
$\gamma'$ intersects $E_{\Sigma^{+}}$. Since $E_{\Sigma^{+}}\cdot\gamma'=0$,
this implies that $\gamma'\subset E_{\Sigma^{+}}$ and hence $\gamma\subset E_{\widetilde{\Sigma}}$,
a contradiction.

We compute the discrepancy of $E_{\widetilde{\Sigma}}.$ By (\ref{align:KSigma})
and (\ref{eq:EF}), we have $-K_{\widetilde{\Sigma}}=(N+1+\dim A-3)H_{\widetilde{\Sigma}}-(\dim A-3)E_{\widetilde{\Sigma}}$.
Since $H_{\widetilde{\Sigma}}\sim_{\mQ}\mu^{*}M_{\Sigma}-\frac{1}{2}F_{\widetilde{\Sigma}}$,
we have 
\begin{equation}
-K_{\widetilde{\Sigma}}\sim_{\mQ}(N+1+\dim A-3)\mu^{*}M_{\widetilde{\Sigma}}-\frac{N+1+\dim A-3}{2}F_{\widetilde{\Sigma}}-(\dim A-3)E_{\widetilde{\Sigma}}.\label{eq:KSigEF}
\end{equation}
 Therefore the discrepancy of $E_{\widetilde{\Sigma}}$ is equal to
$\dim A-3$. Since this is a positive integer, $\Sigma$ has only
Gorenstein terminal singularities along $\mu(E_{\widetilde{\Sigma}}).$
By (\ref{eq:KSigEF}), we have $-K_{\Sigma}=(\dim A+N-2)M_{\Sigma}$.

We show that $\rho(\Sigma$)=1. Since $\widehat{A}\to A$ is the blow-up
along $l$ disjoint projective spaces, and $\widehat{\Sigma}\to\widehat{A}$
is a projective bundle, we see that $\rho(\widehat{\Sigma})=\rho(A)+l+1$.
Since $\widehat{\Sigma}\dashrightarrow\widetilde{\Sigma}$ is small,
we have $\rho(\widetilde{\Sigma})=\rho(\widehat{\Sigma})=\rho(A)+l+1$.
Finally, since $\widetilde{\Sigma}\to\Sigma$ contracts $l+1$ disjoint
divisors, we have $\rho(\Sigma)\leq\rho(\widetilde{\Sigma})-(l+1)=\rho(A)=1$.
Hence $\rho(\Sigma)=1$.

Finally, we show that $M_{\Sigma}$ is primitive. If $M_{\Sigma}$
were not primitive, then $M_{\Sigma}$ would be written as $M_{\Sigma}=\alpha M'_{\Sigma}$
with a primitive Weil divisor $M'_{\Sigma}$ and positive integer
$\alpha\geq2.$ Since $2M'_{\Sigma}$ are Cartier divisors by (2)
and (3), we have $2H_{\widetilde{\Sigma}}+F_{\widetilde{\Sigma}}=\alpha\mu^{*}(2M'_{\Sigma}).$
Hence there is a Cartier divisor $D$ on $\widehat{\Sigma}$ such
that $2H_{\widehat{\Sigma}}+F_{\widehat{\Sigma}}=\alpha D$. Let $l$
be a flopping curve for $\widehat{\Sigma}\dashrightarrow\Sigma^{+}$.
By (\ref{eq:LBLA}) and the proof of Proposition \ref{Prop:SigmaflopG46},
we have $F_{\widetilde{\Sigma}}\cdot l=-1$. This implies that $\alpha D\cdot l=-1$.
This is impossible for $\alpha\geq2$. Therefore $M_{\Sigma}$ is
primitive. 
\end{proof}

\subsection{Application to the three cases \label{subsec:Compensation}}

In this subsection, we produce the situation as in the subsections
\ref{subsec:Basic-set-up} and \ref{subsec:Prescription} for a $\mQ$-Fano
3-fold of genus 4 or 6.

\subsubsection{\textbf{Genus 4 \label{subsec:Genus-4}}}

In this case, we set
\[
A:=Q^{4}\subset\mP((U^{3})^{*}\oplus U^{3})
\]
 with the same equation as that of $\mathsf{Q}$ in Definition \ref{defSigmabarG4},
and 
\[
\Pi_{1}:=\mP((U^{3})^{*}),\,\Pi_{2}:=\mP(U^{3}),
\]
which are certainly contained in $Q^{4}$. We also set 
\[
B:=B_{6}=\mathbb{P}(\Omega_{\mathbb{P}^{2}}^{1}(1))
\]
 with the equation as in the subsection \ref{subsec:g4Crep}, and
\[
\widehat{A}:=\mathbb{P}_{B_{6}}(\mathcal{O}_{B_{6}}(-1,0)\oplus\mathcal{O}_{B_{6}}(0,-1)).
\]
Finally, we set $b$ as the projection morphism 
\[
\mathbb{P}_{B_{6}}(\mathcal{O}_{B_{6}}(-1,0)\oplus\mathcal{O}_{B_{6}}(0,-1))\to B_{6}.
\]

\begin{lem}
\label{Lem:g4B6} There exists a morphism $a\colon\widehat{A}\to A$
which is the blow-up of $\widehat{A}$ along $\Pi=\Pi_{1}\sqcup\Pi_{2}$
and whose exceptional divisor $F_{a}$ is $\mathbb{P}_{B_{6}}(\mathcal{O}_{B_{6}}(-1,0)\oplus0)\sqcup\mathbb{P}_{B_{6}}(0\oplus\mathcal{O}_{B_{6}}(0,-1))$.
The pull-back of $\sO_{A}(1)$ on $\widehat{A}$ is the tautological
line bundle associated with $\mathcal{O}_{B_{6}}(-1,0)\oplus\mathcal{O}_{B_{6}}(0,-1)$.
The triplet $(Q^{4},\Pi=\Pi_{1}\sqcup\Pi_{2},B_{6})$ satisfies the
condition of $(A,\Pi,B)$ as in the subsection \ref{subsec:Basic-set-up}
by setting $d=l=2$. 
\end{lem}

\begin{proof}
Take a point $\mathsf{p}:=[W^{1}\otimes U^{1}]\in B_{6}$, where $W^{1}\subset(U^{3})^{*}$
and $U^{1}\subset U^{3}$ are $1$-dimensional subspaces such that
$U^{1}\subset(V^{1})^{\perp}$ with respect to the dual pairing. The
fiber of the projection $\widehat{A}\to B_{6}$ over $\mathsf{p}$
is $\mP(W^{1}\oplus U^{1})$, which is a linear subspace of $\mP((U^{3})^{*}\oplus U^{3})$.
Note that, for a point $[\bm{y}+\bm{x}]\in\mP(W^{1}\oplus U^{1})$
with $\bm{y}\in W^{1}$ and $\bm{x}\in U^{1}$, it holds that ${\empty}^{t}\!\bm{y}\bm{x}=0$
since $U^{1}\subset(W^{1})^{\perp}$. Therefore the image of $\widehat{A}\to\mP((U^{3})^{*}\oplus U^{3})$
is contained in $Q^{4}$. We denote by $a$ the induced morphism $\widehat{A}\to A=Q^{4}.$
By Lemma \ref{lem:ABFib} (1), the second assertion follows.

Let $\mathsf{q}:=[\bm{y}+\bm{x}]$ be a point of $Q^{4}$ with $\bm{y}\in(U^{3})^{*}$
and $\bm{x}\in U^{3}$. By Lemma \ref{lem:ABFib} (2), the fiber of
$\widehat{A}\to Q^{4}$ over $\mathsf{q}$ is 
\[
\left\{ \mathsf{q}\right\} \times\left\{ [W^{1}\otimes U^{1}]\mid U^{1}\subset(W^{1})^{\perp},\bm{y}\in W^{1},\bm{x}\in U^{1}\right\} .
\]
If $\bm{y}\not=\bm{o}$ and $\bm{x}\not=\bm{o},$ then $W^{1}$ and
$U^{1}$ are uniquely determined as $W^{1}=\mC\bm{y}$ and $U^{1}=\mC\bm{x}$
(since ${\empty}^{t}\!\bm{y}\bm{x}=0$, it holds that $U^{1}\subset(W^{1})^{\perp}$).
Therefore the morphism $\widehat{A}\to Q^{4}$ is dominant, hence
is surjective, and is an isomorphism outside of $\Pi$. The fiber
of $\widehat{A}\to Q^{4}$ over $\bm{o}+\bm{x}\in\Pi_{1}$ is isomorphic
to $\mP((\mC\bm{x})^{\perp})\simeq\mP^{1}$, and the fiber of $\widehat{A}\to Q^{4}$
over $\bm{y}+\bm{o}\in\Pi_{2}$ is isomorphic to $\mP((\mC\bm{y})^{\perp})\simeq\mP^{1}$.
Note that, since $-K_{\widehat{A}}=2H_{\widehat{A}}+b^{*}L_{B},$
we see that $-K_{\widehat{A}}$ is relatively ample for $\widehat{A}\to Q^{4}$.
Therefore, by \cite[Thm.2.3]{An}, $\widehat{A}\to Q^{4}$ is the
blow-up of $Q^{4}$ along $\Pi$ and the $a$-exceptional divisor
is $\mathbb{P}_{B_{6}}(\mathcal{O}_{B_{6}}(-1,0)\oplus0)\sqcup\mathbb{P}_{B_{6}}(0\oplus\mathcal{O}_{B_{6}}(0,-1))$. 

Assumptions 1--2 are clearly satisfied. We check Assumption 3, equivalently,
the relation (\ref{eq:LBLA}) with $d=2$. Take a hyperplane $L\subset\mathrm{S}^{-1,0,1}U^{3}$.
We consider elements of $\mathrm{S}^{1,0,-1}U^{3}$ and $\mathrm{S}^{-1,0,1}U^{3}$
as $3\times3$ traceless matrices. The dual pairing between $U^{3}\otimes(U^{3})^{*}$
and $(U^{3})^{*}\otimes U^{3}$ induces a natural dual pairing between
$\mathrm{S}^{1,0,-1}U^{3}$ and $\mathrm{S}^{-1,0,1}U^{3}$. Explicitly,
for $Z=(z_{ij})\in\mathrm{S}^{-1,0,1}U^{3}$ and $P=(p_{ij})\in\mathrm{S}^{1,0,-1}U^{3}$,
the dual pairing is defined as $(Z,P)\mapsto\sum z_{ij}p_{ij}$. For
$L$, there exists $M=(m_{ij})\in\mathrm{S}^{1,0,-1}U^{3}$ such that
$L=\left\{ \sum_{1\leq i,j\leq3}m_{ij}z_{ij}=0\right\} .$ The above
construction show that $Q^{4}\setminus\Pi\to B_{6}$ is defined by
$[\bm{y}+\bm{x}]\mapsto[\bm{y}\otimes\bm{x}].$ Therefore we see that
$a_{*}b^{*}(B_{6}\cap\mP(L))=Q^{4}\cap{}$$\left\{ {\empty}^{t}\!\bm{y}M\bm{x}=0\right\} $,
which is a quadric section of $Q^{4}$ containing $\Pi$. We can explicitly
check that a general $a_{*}b^{*}(B_{6}\cap\mP(L))$ is generically
smooth along $\Pi$. Since $b^{*}(B_{6}\cap\mP(L))$ does not contain
the $a$-exceptional divisor $F_{a}$, it is the strict transform
of $a_{*}b^{*}(B_{6}\cap\mP(L))$ for a general $L$. Therefore the
relation (\ref{eq:LBLA}) holds with $d=2$. 
\end{proof}
In the following subsections \ref{subsec:g6QAB} and \ref{subsec:g6Ctwisted},
we use the notation as in the subsection \ref{subsec:g6setup}.

\subsubsection{\textbf{Genus 6, $\mathtt{Q}$-type \label{subsec:g6QAB}}}

In this case, we set 
\[
A:=A_{\mathtt{Q}}=\rG(2,V)\cap\mP(V'\oplus U^{5})
\]
 and $\Pi$ the same as in the subsection \ref{subsec:g6setup}. We
recall that the projection of ${\rm G}(2,V)$ from the $3$-plane
$\Pi$ induces the natural rational map $A_{\mathtt{Q}}\dashrightarrow{\rm G(2,V')}\cap\mP(U^{5})$
and the target ${\rm G(2,V')}\cap\mP(U^{5})$ is the smooth quadric
3-fold $Q^{3}$. We set 
\[
B:=Q^{3}={\rm G(2,V')}\cap\mP(U^{5}),
\]
and 
\[
\widehat{A}:=\widehat{A}_{\mathtt{Q}}:=\mP_{Q^{3}}(\sU|_{Q^{3}}\oplus\sO_{Q^{3}}(-1)),
\]
where $\sU$ is the rank two universal subbundle on $\rG(2,V')$.
Finally we set $b$ as the projection morphism 

\[
\mP_{Q^{3}}(\sU|_{Q^{3}}\oplus\sO_{Q^{3}}(-1))\to Q^{3}.
\]

\begin{lem}
\label{Lem:g6Q} There exists a morphism $\widehat{A}_{\mathtt{Q}}\to A_{\mathtt{Q}}$
which is the blow-up of $\widehat{A}_{\mathtt{Q}}$ along $\Pi$ and
whose exceptional divisor is $\mP_{Q^{3}}(\sU|_{Q^{3}}\oplus0)$.
The pull-back of $\sO_{A_{\mathtt{Q}}}(1)$ on $\widehat{A}_{\mathtt{Q}}$
is the tautological line bundle associated with $\sU|_{Q^{3}}\oplus\sO_{Q^{3}}(-1).$
The triplet $(A_{\mathtt{Q}},\Pi,Q^{3})$ satisfies the condition
of $(A,\Pi,B)$ as in the subsection \ref{subsec:Basic-set-up} by
setting $d=l=1$.
\end{lem}

\begin{proof}
Since we can show the first two assertions in a quite similar way
to Lemma \ref{Prop:g8Sigma'}, we only show that $\widehat{A}_{\mathtt{Q}}\to A_{\mathtt{Q}}$
is the blow-up along $\Pi$. Note that the restriction of the morphism
$\widehat{A}_{\mathtt{Q}}\to A_{\mathtt{Q}}$ over $\Pi$ is $\mathbb{P}_{Q^{3}}(\sU|_{Q^{3}}\oplus0)\to\mP(V'\oplus0)=\Pi\simeq\mP^{3},$
which can be identified with the natural morphism $\mathbb{P}_{Q^{3}}(\sU|_{Q^{3}})\to\mP(V')$
to $\mP(V')\simeq\mP^{3}$ from the total space of lines in $\mP^{3}$
parameterized by $Q^{3}\subset\rG(2,V')$. By \cite[Prop.3.4]{SW},
$\mathbb{P}_{Q^{3}}(\sU|_{Q^{3}})\to\mP^{3}$ is the projectivization
of the null-correlation bundle. Note that, since $-K_{\widehat{A}_{\mathtt{Q}}}=2H_{\widehat{A}_{\mathtt{Q}}}+L_{Q^{3}}$
where $L_{Q^{3}}$ is the pull-back of $\sO_{Q^{3}}(1)$, we see that
$-K_{\widehat{A}}$ is relatively ample for $\widehat{A}_{\mathtt{Q}}\to A_{\mathtt{Q}}$
Therefore, by \cite[Thm.2.3]{An}, $\widehat{A}_{\mathtt{Q}}\to A_{\mathtt{Q}}$
is the blow-up of $A_{\mathtt{Q}}$ along $\Pi$ and the exceptional
divisor is $\mathbb{P}_{Q^{3}}(\sU\oplus0)$. By the construction,
the pull-back of $\sO_{A_{\mathtt{Q}}}(1)$ on $\widehat{A}_{\mathtt{Q}}$
is the tautological line bundle associated with $\mP_{Q^{3}}(\sU|_{Q^{3}}\oplus\sO_{Q^{3}}(-1))$. 

Assumptions 1--2 are clearly satisfied. Since $A_{\mathtt{Q}}\dashrightarrow Q^{3}$
is the restriction of the projection from $\Pi$, the relation (\ref{eq:LBLA})
follows. 
\end{proof}

\subsubsection{\textbf{Genus 6, $\mathtt{C}$-type \label{subsec:g6Ctwisted}}}

In this case, we set 
\[
A:=A_{\mathtt{C}}=\rG(2,V)\cap\mP(U^{8}),
\]
$\Pi$ the same as in the subsection \ref{subsec:g6setup}, and $a\colon\widehat{A}\to A$
the blow-up of $A$ along $\Pi$. We recall that the projection of
${\rm G}(2,V)$ from the 2-plane $\Pi$ induces the natural rational
map $A_{\mathtt{C}}\dashrightarrow\mP(U^{5})$. We set 
\[
B:=\mP(U^{5})\simeq\mP^{4},
\]
 and $b\colon\widehat{A}\to B$ the naturally induced morphism. 
\begin{lem}
\label{Lem:g6B3} The triplet $(A,\Pi,\mP^{4})$ satisfies the condition
of $(A,\Pi,B)$ as in the subsection \ref{subsec:Basic-set-up} by
setting $d=l=1$. 
\end{lem}

\begin{proof}
The assertion is almost clear. 
\end{proof}

\subsubsection{\textbf{Rationality of the key varieties}}
\begin{cor}
In the genus $4$ or $6$ case, $\Sigma$ is rational.
\end{cor}

\begin{proof}
The assertion follows since $\Sigma$ is birational to a projective
bundle over a rational Fano manifold in the genus 4 or 6 case.
\end{proof}

\subsection{Coincidence between $\overline{\Sigma}$'s in the subsection \ref{subsec:Prescription}
and in the section \ref{sec:Extending Mid G46}}
\begin{lem}
\label{lem:coicideG46} In the genus 4 or 6 case, the variety $\overline{\Sigma}$
as in the subsection \ref{subsec:Prescription} is the same as the
variety $\overline{\Sigma}$ defined as in the section \ref{sec:Extending Mid G46}.
The morphism $\varphi_{|H_{\widehat{\Sigma}}|}\colon\widehat{\Sigma}\to\overline{\Sigma}$
is birational. The $\varphi_{|H_{\widehat{\Sigma}}|}$-image of $E_{\widehat{\Sigma}}$
on $\overline{\Sigma}$ is disjoint from $W$.
\end{lem}

\begin{proof}
The variety $\overline{\Sigma}$ as in the subsection \ref{subsec:Prescription}
defined for the triplet $(A,\Pi,B)$ is contained in $\mP(U_{A}\oplus U_{B}^{*})$
by the fact that $H^{0}(a^{*}\sO_{A}(1))=U_{A}$ and $H^{0}(b^{*}(\Omega_{\mathbb{P}(U_{B})}^{1}(1)|_{B}))=U_{B}^{*}$.
Temporarily, we denote by $\overline{\Sigma}'$ the variety $\overline{\Sigma}$
as in the subsection \ref{sec:Extending Mid G46}, which is also contained
in $\mP(U_{A}\oplus U_{B}^{*})$ in each case . 

First we show that $\overline{\Sigma}\subset\overline{\Sigma}'$.
For this, it suffices to check the $\varphi_{|H_{\widehat{\Sigma}}|}$-image
of a general point of $\widehat{\Sigma}$ is mapped in $\overline{\Sigma}'$
since $\widehat{\Sigma}$ is irreducible. Since $\widehat{A}\subset A\times B$
by Lemma \ref{lem:ABFib} (1) and the results in the subsection \ref{subsec:Compensation},
we can express a point of $\widehat{A}$ as $([\bm{x}],[\bm{y}])$
with $\bm{x}\in U_{A}$ and $\bm{y}\in U_{B}$. Note that $\bm{x}$
satisfies the equation of $A$. The fiber of $\pi\colon\widehat{\Sigma}\to\widehat{A}$
over $\mathsf{p}$ is $\mP(\mC\bm{x}\oplus(U_{B}/\mC\bm{y})^{*})$.
We choose a point $\mathsf{p}:=([\bm{x}],[\bm{y}])\in\widehat{A}$
such that $[\bm{x}]\not\in\Pi$. In this case, it holds that $[\bm{y}]=[b(a^{-1}(\bm{x}))]$.
In each of the three cases, we check that the fiber $\pi^{-1}(\mathsf{p})$
is mapped by $\varphi_{|H_{\widehat{\Sigma}}|}$ into $\overline{\Sigma}'$
in the sequel.

\vspace{3pt}

\noindent Genus 4: We can express a point of $\widehat{A}$ as $([\bm{x}_{1}+\bm{x}_{2}],[\bm{y}])$
with $\bm{x}_{1}\in(U^{3})^{*}$, $\bm{x}_{2}\in U^{3},\bm{y}\in\mathrm{S}^{-1,0,1}U^{3}$.
We are choosing a point $\mathsf{p}:=([\bm{x}_{1}+\bm{x}_{2}],[\bm{y}])\in\widehat{A}_{\mathtt{Q}}$
such that $[\bm{x}_{1}+\bm{x}_{2}]\not\in\Pi$, namely, $\bm{x}_{1}\not=\bm{o}$
and $\bm{x}_{2}\not=\bm{o}.$ In this case, it holds that $[\bm{y}]=[\bm{x}_{2}\otimes\bm{x}_{1}]$
by the proof of Lemma \ref{Lem:g4B6} and hence the $\pi$-fiber over
$\mathsf{p}$ is $\mP(\mC(\bm{x}_{1}+\bm{x}_{2})\oplus(U^{5}/\mC$$(\bm{x}_{2}\otimes\bm{x}_{1}))^{*})$.
Therefore, by the definition of $\overline{\Sigma}'$ as in \ref{defSigmabarG4},
we see that $\overline{\Sigma}\subset\overline{\Sigma}'$. 

\vspace{3pt}

\noindent Genus 6, $\mathtt{Q}$-type: We can express a point of
$\widehat{A}_{\mathtt{Q}}$ as $([\bm{x}_{1}+\bm{x}_{2}],[\bm{y}])$
with $\bm{x}_{1}\in V'$, $\bm{x}_{2},\bm{y}\in U^{5}$. We are choosing
a point $\mathsf{p}:=([\bm{x}_{1}+\bm{x}_{2}],[\bm{y}])\in\widehat{A}_{\mathtt{Q}}$
such that $[\bm{x}_{1}+\bm{x}_{2}]\not\in\Pi$, namely, $\bm{x}_{2}\not=\bm{o}.$
In this case, it holds that $[\bm{y}]=[\bm{x}_{2}]$ since $A_{\mathtt{Q}}\dashrightarrow Q^{3}$
is the projection from $\Pi$, and hence the $\pi$-fiber over $\mathsf{p}$
is $\mP(\mC(\bm{x}_{1}+\bm{x}_{2})\oplus(U^{5}/\mC$$\bm{x}_{2})^{*})$.
Therefore, by the definition of $\overline{\Sigma}'$ as in \ref{Def:SigmaBarG6Q},
we see that $\overline{\Sigma}\subset\overline{\Sigma}'$. 

\vspace{3pt}

We can show that $\overline{\Sigma}\subset\overline{\Sigma}'$ in
the case of genus 6 and $\mathtt{C}$-type in a similar way to the
case of genus 6 and $\mathtt{Q}$-type, so we omit a proof.

Now we check that $\overline{\Sigma}=\overline{\Sigma}'$. Since $\dim\widehat{\Sigma}=\dim\overline{\Sigma}'$,
it suffices to show a general $\varphi_{|H_{\widehat{\Sigma}}|}$-fiber
consists of one point. This also implies that $\varphi_{|H_{\widehat{\Sigma}}|}\colon\widehat{\Sigma}\to\overline{\Sigma}$
is birational. We take a point $\mathsf{t}:=[\bm{t}_{1}+\bm{t}_{2}]\in\overline{\Sigma}\setminus\overline{\Pi}$
with $\bm{t}_{1}\in U_{A}\setminus\{\bm{o}\}$ and $\bm{t}_{2}\in U_{B}^{*}$.
Since $\pi^{-1}(F_{a})=\varphi_{|H_{\widehat{\Sigma}}|}^{-1}(\overline{\Pi})$
, we have {\small{}
\[
\varphi_{|H_{\widehat{\Sigma}}|}^{-1}(\mathsf{t})=\left\{ \mathsf{t}\right\} \times\left\{ ([\bm{x}],[b(a^{-1}(\bm{x}))],\mathsf{t})\in\widehat{\Sigma}\mid[\bm{x}]\not\in\Pi,\mathsf{t}_{1}\in\mC\bm{x},\bm{t}_{2}\in(U_{B}/(\mC b(a^{-1}(\bm{x}))))^{*})\right\} 
\]
} by Lemma \ref{lem:ABFib} (2). This is nonempty since we take $\mathsf{t}$
in the $\varphi_{|H_{\widehat{\Sigma}}|}$-image of $\widehat{\Sigma}$.
Moreover, it consist of one point as desired since $[\bm{x}]=[\bm{t}_{1}]$. 

We show the last assertion. Note that the $\varphi_{|H_{\widehat{\Sigma}}|}$-image
of $E_{\widehat{\Sigma}}=\mP(0\oplus b^{*}(\Omega_{\mathbb{P}(U_{B})}^{1}(1)|_{B}))$
on $\overline{\Sigma}$ coincides with $\mP(0\oplus U_{B}^{*})$.
In the genus 4 case, this is equal to $\overline{\Pi}_{1}\cap\overline{\Pi}_{2}=\mP(0\oplus0\oplus\mathrm{S}^{1,0,-1}U^{3})$
and is disjoint from $W$ since $W\cap(\overline{\Pi}_{1}\cap\overline{\Pi}_{2})=\Pi_{1}\cap\Pi_{2}=\emptyset$
by Proposition \ref{Prop:disj}. In the genus 6 case, $\overline{\Sigma}$
has non-hypersurface singularities along $\mP(0\oplus U_{B}^{*})$,
hence must be disjoint from $W$.
\end{proof}

\subsection{Embedding theorem \label{subsec:Embedding-theoremG46}}

Now we show Theorem \ref{thm:main1} for a prime $\mQ$-Fano $3$-fold
$X$ of genus 4 or 6.
\begin{thm}
\label{thm:embg4g6} A $\mQ$-Fano $3$-fold $X$ of genus $4$ or
$6$ is a linear section of $\Sigma$. 
\end{thm}

\begin{proof}
The proof which will be given below is more or less the same as that
of Theorem \ref{thm:emb8} but is slightly involved, so we write it
for readers' convenience. 

Note that $W\cap\Sing\overline{\Sigma}$ is $0$-dimensional since
$W$ has only terminal singularities and $W$ is a linear section
of $\overline{\Sigma}$ with respect to $|\sO_{\overline{\Sigma}}(1)|$.
By Lemmas \ref{lem:coicideG46}, $W$ is disjoint from the image of
$E_{\widetilde{\Sigma}}$. Therefore, since $\widetilde{\Sigma}\to\overline{\Sigma}$
is crepant and small and nontrivial fibers are 1-dimensional over
$W$, the strict transform $W_{\widetilde{\Sigma}}$ of $W$ in $\widetilde{\Sigma}$
is a linear section of $\widetilde{\Sigma}$ with respect to $|H_{\widetilde{\Sigma}}|$
and hence the restriction $W_{\widetilde{\Sigma}}\to W$ of $\widetilde{\Sigma}\to\overline{\Sigma}$
over $W$ is also crepant and small. Since $W$ has only terminal
singularities and $W_{\widetilde{\Sigma}}\to W$ is crepant, we see
that $W_{\widetilde{\Sigma}}$ is normal and has only terminal singularities
by \cite[the proof of Prop.16.4]{CKM}. Note that $F_{\widetilde{\Sigma}}|_{W_{\widetilde{\Sigma}}}$
is the strict transform of $\Pi$ and is relatively ample for $W_{\widetilde{\Sigma}}\to W$.
Since $Y\to W$ in the genus 6 case (resp. $Z\to W$ in the genus
4 case) is the unique small extraction such that the strict transform
of $\Pi$ is relatively ample, we see that $Y=W_{\widetilde{\Sigma}}$
in the genus 6 case (resp. $Z=W_{\widetilde{\Sigma}}$ in the genus
4 case). Since we may write $Y=W_{\widetilde{\Sigma}}=\widetilde{H}_{1}\cap\dots\cap\widetilde{H}_{\dim\widetilde{\Sigma}-3}$
with $\widetilde{H}_{i}\in|H_{\widetilde{\Sigma}}|$ ($1\leq i\leq\dim\widetilde{\Sigma}-3$),
we see that $X=M_{1}\cap\dots\cap M_{\dim\widetilde{\Sigma}-3}$ with
the image $M_{i}\in|M_{\Sigma}|$ of $\widetilde{H}_{i}$ as desired.
\end{proof}

\subsection{Extension of the Sarkisov link}

We have obtained the following diagram in the genus 4 or 6 case:\begin{equation}\label{eq:SarkisovG46} \xymatrix{& \widetilde{\Sigma}\ar@{-->}^{\text{anti-flip}\ }[r]\ar[dl]_{\mu}\ar[dr] & \Sigma^+\ar@{-->}^{\text{flop}}[r]\ar[d]& \widehat{\Sigma}\ar[dr]^{b\circ \pi}\ar[dl]\\
\Sigma & &\overline{\Sigma} & & B.}
\end{equation}
\begin{cor}
\label{cor:ExtSarkisov} In the case of genus 4, the restriction of
(\ref{eq:SarkisovG46}) to 3-folds, we obtain the following diagram:
\begin{equation}\label{eq:SarkisovG4} \xymatrix{& Z\ar@{-->}^{\text{flop}}[rr]\ar[dl]_g\ar[dr] & & Z'\ar[dr]^{g'}\ar[dl]\\
X & & W & & B_6,}
\end{equation}where $Z$ and $W$ defined as in the subsection \ref{subsec:Genus-4mid}
are linear sections of $\widetilde{{\Sigma}}$ and $\overline{\Sigma}$
with respect to $|H_{\widetilde{\Sigma}}|$ and $|\sO_{\overline{\Sigma}}(1)|$
respectively; $Z'$ is defined as the corresponding linear section
of $\widehat{\Sigma}$ with respect to $|H_{\widehat{\Sigma}}|$;
the restriction of the anti-flip to $Z$ is the identity. Moreover,
the following assertions hold:

\noindent $(1)$ The morphism $g'\colon Z'\to B_{6}$ is the blow-up
of $B_{6}$ along a smooth curve $C'$ of genus $8$ isomorphic to
$C$.

\vspace{3pt}

\noindent $(2)$ The curve $C'$ is the complete intersection of
the strict transforms of the $g$-exceptional divisors, which are
divisors of types $(2,1)$ and $(1,2).$

In the case of genus 6, (\ref{eq:SarkisovG46}) is an extension of
(\ref{eq:Sarkisov}), where $Y$,$W$ and $Y'$ are linear sections
of $\widetilde{{\Sigma}}$, $\overline{\Sigma}$ and $\widehat{\Sigma}$
with respect to $|H_{\widetilde{\Sigma}}|$, $|\sO_{\overline{\Sigma}}(1)|$
and $|H_{\widehat{\Sigma}}|$ respectively, and $X'=B=Q^{3}$ in the
case of $\mathtt{Q}$-type (resp. $X'$ is a cubic 3-fold in $\mP(U^{5})$
in the case of $\mathtt{C}$-type). The restriction of the anti-flip
to $Y$ is the identity. 
\end{cor}

\begin{proof}
\noindent The genus 4 case: The restriction of the anti-flip to $Z$
is the identity since the flipped locus in $\widetilde{\Sigma}$ is
contained in $E_{\widetilde{\Sigma}}$ by Proposition \ref{Prop:flipping G46}
(2), and $E_{\widetilde{\Sigma}}$ is disjoint from $Z'$ by Lemma
\ref{lem:coicideG46}. The rest of the assertion except (1) and (2)
easily follow from the proof of Theorem \ref{thm:embg4g6}. 

\vspace{3pt}\noindent $(1)$. Since $\widehat{\Sigma}\to B_{6}$
is a $\mP^{8}$-bundle, and $Z'$ is a linear section of $\widehat{\Sigma}$
with respect to $|H_{\widehat{\Sigma}}|$ of codimension 8, we see
that $Z'\to B_{6}$ is birational and $-K_{Z'}=\pi^{*}a^{*}L_{Q^{4}}|_{Z'}$
by (\ref{align:KSigma}). Since $E_{\widehat{\Sigma}}$ is disjoint
from $Z'$, we have $H_{\widehat{\Sigma}}|_{Z'}\sim(\pi^{*}a^{*}L_{A})|_{Z'}$
by (\ref{eq:EF}). Therefore we have $-K_{Z'}=H_{\widehat{\Sigma}}|_{Z'}$.
Since $H_{\widehat{\Sigma}}$ is relatively ample over $B_{6}$, so
is $-K_{Z'}$ . Since $\rho(Z')=3$ and $\rho(B_{6})=2$, the relative
Picard number of the morphism $Z'\to B_{6}$ is 1. Therefore, by \cite{Mo},
$Z'\to B_{6}$ is the blow-up of $B_{6}$ at a point or along a curve
$C'$. Comparing the Intermediate Jacobians of $Z'$ and the 3-fold
obtained by blowing up of $Y'$ at the $1/2$-singularity, we see
that $Z'\to B_{6}$ is the blow-up of $B_{6}$ along a curve $C'$
such that $C\simeq C'$ as desired. 

\vspace{3pt}\noindent $(2)$. Since the images of $E_{1}$ and $E_{2}$
on $W$ are disjoint by Proposition \ref{Prop:disj}, the strict transforms
$E'_{1}$ and $E_{2}'$ on $Z'$ of $E_{1}$ and $E_{2}$ respectively
are also disjoint. Therefore $C'$ is set-theoretically the intersection
between the strict transforms $E_{1}''$ and $E_{2}''$ on $B_{6}$
of $E_{1}$ and $E_{2}$. Moreover, since $Z'\to B_{6}$ is the blow-up
along $C'$, and $E'_{1}\cap E_{2}'=\emptyset$, it holds that $C'$
is the complete intersection of $E'_{1}$ and $E'_{2}$. Note that
the anticanonical morphism $Z'\to W$ is induced from the restriction
of $\widehat{\Sigma}\to Q^{4}$ since $-K_{Z'}=\pi^{*}a^{*}L_{Q^{4}}|_{Z'}$
as we saw in the proof of (1). Therefore we have $F_{\widehat{\Sigma}}|_{Z'}=E_{1}'\sqcup E'_{2}$.
From this, we obtain $(2\pi^{*}a^{*}L_{Q^{4}})|_{Z'}-E'_{1}-E'_{2}=(\pi^{*}b^{*}L_{B_{6}})|_{Z'}$
by the equation (\ref{eq:LBLA}). Since $-K_{Z'}=H_{\widehat{\Sigma}}|_{Z'}$
as we saw in the proof of (1), we have $-2K_{Z'}-E'_{1}-E'_{2}=(\pi^{*}b^{*}L_{B_{6}})|_{Z'}$.
On the other hand, we have $-K_{Z'}=(2\pi^{*}b^{*}L_{B_{6}})|_{Z'}-E_{C'}$,
where $E_{C'}$ is the exceptional divisor of the blow-up $Z'\to B_{6}$.
Therefore we obtain $E''_{1}+E''_{2}=3L_{B_{6}}$. For the curve $C'$
of genus 4 to be the complete intersection of $E''_{1}$ and $E''_{2}$,
it must holds that $E''_{1}$ and $E''_{2}$ are of types $(2,1)$
and $(1,2).$

\vspace{3pt}

\noindent The genus 6 case: We can show the assertions in a similar
and simpler way as in the genus 4 case.
\end{proof}

\subsection{Singularity of $\Sigma$ along $\mP((U_{B})^{*})$\label{subsec:Singg4g6}}

By Theorem \ref{thm:common}, the birational morphism $\mu\colon\widetilde{\Sigma}\to\Sigma$
contracts $E_{\widetilde{\Sigma}}$ onto $\mP((U_{B})^{*})$. In this
subsection, we describe the morphism $\mu|_{E_{\widetilde{\Sigma}}}\colon E_{\widetilde{\Sigma}}\to\mP((U_{B})^{*}).$
This follows by studying how fibers of the morphism $E_{\widehat{\Sigma}}\to\mP((U_{B})^{*})$
are transformed by the flop $\widehat{\Sigma}\dashrightarrow\Sigma^{+}$
and the flip $\Sigma^{+}\dashrightarrow\widetilde{\Sigma}$. 

We note that the natural morphism $\mP(\Omega_{\mP(U_{B})}(1))\to\mP((U_{B})^{*})$
is the universal family of hyperplanes of $\mP(U_{B})$. Therefore
the naturally induced morphism $E_{\widehat{\Sigma}}\to\mP((U_{B})^{*})$
is the universal family of the members of $|b^{*}L_{B}|$. In particular,
$E_{\widehat{\Sigma}}\to\mP((U_{B})^{*})$ is flat.

Note that the restrictions of the flopping and flipping contractions
on the strict transforms of $E_{\widehat{\Sigma}}$ are defined over
$\mP((U_{B})^{*})$. The strict transform $E_{\Sigma^{+}}$ on $\Sigma^{+}$
of $E_{\widehat{\Sigma}}$ is smooth since $E_{\widehat{\Sigma}}\dashrightarrow E_{\Sigma^{+}}$
is also a flop of Atiyah type by Proposition \ref{Prop:SigmaflopG46}.
By Proposition \ref{Prop:flipping G46} and the construction of the
flip $\Sigma^{+}\dashrightarrow\widetilde{\Sigma}$, we see that $E_{\widetilde{\Sigma}}$
is smooth since $E_{\Sigma^{+}}\dashrightarrow E_{\widetilde{\Sigma}}$
is the blow-up along $l$ $\mu|_{E_{\widetilde{\Sigma}}}$-sections
whose exceptional divisor is $\mP((U_{B})^{*})\times\Pi$. 

By the description of $E_{\widehat{\Sigma}}\to\mP((U_{B})^{*})$ and
$E_{\widehat{\Sigma}}\dashrightarrow E_{\Sigma^{+}}\dashrightarrow E_{\widetilde{\Sigma}}$
as above, the morphism $E_{\widetilde{\Sigma}}\to\mP((U_{B})^{*})$
is also flat.

We denote by $\widehat{\Gamma}$, $\Gamma^{+}$, and $\widetilde{\Gamma}$
a general fiber of $E_{\widehat{\Sigma}}\to\mP((U_{B})^{*})$ and
its strict transforms on $\Sigma^{+}$ and $\widetilde{\Sigma}$ respectively.
By the argument as above, we see that the restriction $\widehat{\Gamma}\dashrightarrow\Gamma^{+}$
to $\widehat{\Gamma}$ of the flop $\widehat{\Sigma}\dashrightarrow\Sigma^{+}$
is also a flop of Atiyah type, and the restriction $\Gamma^{+}\dashrightarrow\widetilde{\Gamma}$
to $\Gamma^{+}$ of the flip $\Sigma^{+}\dashrightarrow\widetilde{\Sigma}$
is the blow-up of $\widetilde{\Gamma}$ at $l$ smooth points whose
exceptional divisor is $\Pi$. Moreover, $\Gamma^{+}$ and $\widetilde{\Gamma}$
are general fibers of $E_{\Sigma^{+}}\to\mP((U_{B})^{*})$ and $E_{\widetilde{\Sigma}}\to\mP((U_{B})^{*})$
respectively. We set $F_{\widehat{\Gamma}}:=F_{\widehat{\Sigma}}|_{\widehat{\Gamma}}$
and $F_{\Gamma^{+}}:=F_{\Sigma^{+}}|_{\Gamma^{+}}$. 

Hereafter we consider separately in each case and determine $\widetilde{\Gamma}$.

\subsubsection{\textbf{Genus 4 }}
\begin{prop}
A general fiber $\widetilde{\Gamma}$ of the morphism $E_{\widetilde{\Sigma}}\to\mP(\mathrm{S}^{1,0,-1}U^{3})$
is $\mP^{1}\times\mP^{1}\times\mP^{1}$. 
\end{prop}

\begin{proof}
Since a general fiber of $\mP(\Omega_{\mP(\mathrm{S}^{-1,0,1}U^{3})}(1)|_{B_{6}})\to B_{6}$
is a smooth sextic del Pezzo surface $S$, $\widehat{\Gamma}$ is
isomorphic to $\mathbb{P}_{S}(\mathcal{O}_{S}(-1,0)\oplus\mathcal{O}_{S}(0,-1))$
by Lemma \ref{Lem:g4B6}. The divisor $F_{\widehat{\Gamma}}$ of $\widehat{\Gamma}$
consists of $\widehat{G}_{1}:=\mathbb{P}_{S}(\mathcal{O}_{S}(-1,0)\oplus0)\simeq S$
and $\widehat{G}_{2}:=\mathbb{P}_{S}(0\oplus\mathcal{O}_{S}(0,-1))\simeq S$.
It is easy to see the assertion as in the following steps:

\begin{itemize}

\item Let $\widehat{\Gamma}\to\widehat{\Gamma}'$ be the flopping
contraction, which is the restriction of $\widehat{\Sigma}\to\widehat{\Sigma}'$.
This induce the morphisms $\widehat{G}_{1}\to\mP(U_{3}^{*}\oplus0)$
and $\widehat{G}_{2}\to\mP(0\oplus U_{3}),$ each of which is a contraction
of three $(-1)$-curves. These can be identified with the restrictions
of $\widehat{G}_{1}\dashrightarrow G_{1}^{+}$ and $\widehat{G}_{2}\dashrightarrow G_{2}^{+}$
respectively, where $G_{1}^{+}$ and $G_{2}^{+}$ are the strict transforms
of $\widehat{G}_{1}$ and $\widehat{G}_{2}$ respectively on $\Gamma^{+}$. 

\item The restriction $\Gamma^{+}\dashrightarrow\widetilde{\Gamma}$
to $\Gamma^{+}$ of the flip $\Sigma^{+}\dashrightarrow\widetilde{\Sigma}$
is the blow-up of $\widetilde{\Gamma}$ at two smooth points whose
exceptional divisor consists of $G_{1}^{+}$ and $G_{2}^{+}$.

\item $\widetilde{\Gamma}$ is isomorphic to $\mP^{1}\times\mP^{1}\times\mP^{1}$.
Indeed, it holds that $\rho(\widetilde{\Gamma})=3$ since $\rho(\widetilde{\Gamma})+2=\rho(\Gamma^{+})=\rho(\widehat{\Gamma})=5.$
Moreover we see that $\widetilde{\Gamma}$ is a sextic del Pezzo 3-folds
as follows: it holds that $(-K_{\widetilde{\Gamma}})^{3}=(-K_{\Gamma^{+}})^{3}+2=(-K_{\widehat{\Gamma}})^{3}+2=6$.
Since $-K_{\widehat{\Gamma}}=2H_{\widehat{\Gamma}}$, we have $-K_{\widetilde{\Gamma}}=2H_{\widetilde{\Gamma}}$
where $H_{\widetilde{\Gamma}}$ is the strict transform of $H_{\widehat{\Gamma}}$.
Since $-K_{\widehat{\Sigma}}$ is nef and big and is numerically trivial
only for flopping curves, $-K_{\Gamma^{+}}$ is nef and big and is
numerically trivial only for flopped curves. Therefore, $-K_{\widetilde{\Sigma}}$
is ample since flopped curves is numerically positive for the exceptional
divisor $G_{1}^{+}\cup G_{2}^{+}$ of the blow-up $\Gamma^{+}\to\widetilde{\Gamma}.$
Since $\widetilde{\Gamma}$ is a sextic del Pezzo 3-folds of $\rho(\widetilde{\Gamma})=3$,
$\widetilde{\Gamma}\simeq\mP^{1}\times\mP^{1}\times\mP^{1}$ by \cite[Thm.5.16]{Fuj1}.

\end{itemize}
\end{proof}

\subsubsection{\textbf{Genus 6, $\mathsf{\mathtt{Q}}$-type}}
\begin{prop}
A general fiber $\widetilde{\Gamma}$ of the morphism $E_{\widetilde{\Sigma}}\to\mP((U^{5})^{*})$
is $\mP^{2}\times\mP^{2}$. 
\end{prop}

\begin{proof}
Since a general fiber of $\mP(\Omega_{\mP(U^{5})}(1)|_{Q^{3}})\to Q^{3}$
is $\mP^{1}\times\mP^{1}$, $\widehat{\Gamma}$ is isomorphic to $\mP_{\mP^{1}\times\mP^{1}}(\sO(-1,0)\oplus\sO(0,-1)\oplus\sO(-1,-1))$
by the definition of $\widehat{A}_{\mathtt{Q}}$ as in the subsection
\ref{subsec:g6QAB}. The divisor $F_{\widehat{\Gamma}}$ of $\widehat{\Gamma}$
is $\mP_{\mP^{1}\times\mP^{1}}(\sO(-1,0)\oplus\sO(0,-1)\oplus0)$.
Referring to \cite[Thm.5.1]{Fuk} for details, we see the assertion
as in the following steps:

\begin{itemize}

\item The flopping locus of the flop $\widehat{\Gamma}\dashrightarrow\Gamma^{+}$
is $\mP_{\mP^{1}\times\mP^{1}}(\sO(-1,0)\oplus0\oplus0)\sqcup\mP_{\mP^{1}\times\mP^{1}}(0\oplus\sO(0,-1)\oplus0)$.
The divisor $F_{\Gamma^{+}}$ of $\Sigma^{+}$ is $\mP^{3}$. 

\item The restriction $\Gamma^{+}\dashrightarrow\widetilde{\Gamma}$
to $\Gamma^{+}$ of the flip $\Sigma^{+}\dashrightarrow\widetilde{\Sigma}$
is the blow-up of $\widetilde{\Gamma}$ at a smooth point whose exceptional
divisor is $F_{\Gamma^{+}}\simeq\mP^{3}$.

\item $\widetilde{\Gamma}$ is isomorphic to $\mP^{2}\times\mP^{2}$.

\end{itemize}
\end{proof}

\subsubsection{\textbf{Genus 6, $\mathtt{C}$-type }}
\begin{prop}
A general fiber $\widetilde{\Gamma}$ of the morphism $E_{\widetilde{\Sigma}}\to\mP((U^{5})^{*})$
is $\mP^{1}\times\mP^{1}\times\mP^{1}$. 
\end{prop}

\begin{proof}
By \cite{Fuj3}, $b\colon\widehat{A}_{\mathtt{C}}\to\mP(U^{5})$ is
the blow-up of $\mP(U^{5})$ along a twisted cubic $\gamma$. Since
a general fiber $H$ of $\mP(\Omega_{\mP(U^{5})}(1))\to\mP((U^{5})^{*})$
is a hyperplane of $\mP(U^{5})$, $\widehat{\Gamma}$ is isomorphic
to the $3$-fold obtained by blowing up $H\simeq\mP^{3}$ along $H\cap\gamma$
which consists of three points $\mathsf{p}_{1},\mathsf{p}_{2},\mathsf{p}_{3}$
in a general position. Referring to \cite[Thm.4.1]{Fuk} for details,
we see the assertion as in the following steps:

\begin{itemize}

\item The flopping locus of the flop $\widehat{\Gamma}\dashrightarrow\Gamma^{+}$
consists of the strict transforms of three lines $l_{ij}$ through
$\mathsf{p}_{i}$ and $\mathsf{p}_{j}$ ($1\leq i<j\leq3).$ The divisor
$F_{\Gamma^{+}}$ of $\Sigma^{+}$ is $\mP^{2}$. 

\item The restriction $\Gamma^{+}\dashrightarrow\widetilde{\Gamma}$
to $\Gamma^{+}$ of the flip $\Sigma^{+}\dashrightarrow\widetilde{\Sigma}$
is the blow-up of $\widetilde{\Gamma}$ at a smooth point whose exceptional
divisor is $F_{\Gamma^{+}}\simeq\mP^{2}$.

\item $\widetilde{\Gamma}$ is isomorphic to $\mP^{1}\times\mP^{1}\times\mP^{1}$.

\end{itemize}
\end{proof}

\subsubsection{\textbf{Comparison of $\Sigma'$ and $\widehat{\Sigma}$} \label{subsec:SigmaSigma'}}

In this subsection, we clarify the relationship between $\widehat{\Sigma}$
as in the subsection \ref{subsec:Prescription} and $\Sigma'$ as
in the subsections \ref{subsec:g4Crep} and \ref{subsec:g6Crep}.
Setting 
\[
\sF:=\begin{cases}
\sO_{B_{6}}(-1,0)\oplus\sO_{B_{6}}(0,-1): & \text{{genus} 4,}\\
\sU|_{Q^{3}}\oplus\sO_{Q^{3}}(-1): & \text{{genus} 6, \ensuremath{\mathtt{Q}}\text{{-type,}}}
\end{cases}
\]
we may write $\widehat{A}=\mP_{B}(\sF)$ and $\Sigma'=\mP_{B}(\sF\oplus(\Omega_{\mP(U_{B})}^{1}(1)|_{B}))$
in each case.
\begin{prop}
\label{Prop:hatTo'} There exists a naturally induced birational morphism
$\tau\colon\widehat{\Sigma}\to\Sigma'$ over $\overline{\Sigma}$
and its exceptional locus coincides with $E_{\widehat{\Sigma}}$.
The morphism $\tau$ is the blow-up of $\Sigma'$ along $\mP_{B}(0\oplus\Omega_{\mP(U_{B})}^{1}(1)|_{B})$. 
\end{prop}

\begin{proof}
By Lemmas \ref{Lem:g4B6} and \ref{Lem:g6Q}, we have a surjection
$b^{*}\sF^{*}\to a^{*}\sO_{A}(1)$, which induces the following natural
morphism $\tau$: 
\begin{align*}
\widehat{\Sigma} & =\mP_{\widehat{A}}(a^{*}\sO_{A}(-1)\oplus b^{*}(\Omega_{\mP(U_{B})}^{1}(1)|_{B}))\\
 & \hookrightarrow\mP_{\widehat{A}}(b^{*}(\sF\oplus\Omega_{\mP(U_{B})}^{1}(1)|_{B})):=\Sigma''\\
 & \to\mP_{B}(\sF\oplus\Omega_{\mP(U_{B})}^{1}(1)|_{B})=\Sigma',
\end{align*}
 where the former is the inclusion morphism of projective bundles,
and the latter is a $\mP^{1}$-bundle since it is the base change
of the $\mP^{1}$-bundle $b\colon\widehat{A}=\mathbb{P}_{B}(\sF)\to B$.
By this construction, the $\tau$-pull-back of $\sO_{\Sigma'}$(1)
coincides with $\sO_{\widehat{\Sigma}}(1)$. Therefore the composite
of the morphism $\tau$ and $\varphi_{|H_{\Sigma'}|}\colon\Sigma'\to\overline{\Sigma}$
coincides with $\varphi_{|H_{\widehat{\Sigma}}|}$. This implies that
$\tau$ is birational since so is $\varphi_{|H_{\widehat{\Sigma}}|}$
by Lemma \ref{lem:coicideG46}.

By Lemma \ref{lem:excepflopE G46}, $E_{\widehat{\Sigma}}$ is contracted
by $\tau$ since $\varphi_{|H_{\Sigma'}|}$ is small. By the description
of $\varphi_{|H_{\widehat{\Sigma}}|}$-fibers as in Lemma \ref{lem:excepflopE G46},
and the description of $\varphi_{|H_{\Sigma'}|}$-fibers as in Propositions
\ref{prop:graph} and \ref{prop:graphg6}, $\tau$ is isomorphic outside
$E_{\widehat{\Sigma}}$. Note that $\tau$ induces $E_{\widehat{\Sigma}}=\mP_{\widehat{A}}(b^{*}(\Omega_{\mP(U_{B})}^{1}(1)|_{B}))\to\mP_{B}(\Omega_{\mP(U_{B})}^{1}(1)|_{B})$
and this is a $\mP^{1}$-bundle. Moreover, by (\ref{align:KSigma}),
$-K_{\widehat{\Sigma}}$ is $\tau$-ample for the morphism. Therefore
$\tau$ is the blow-up of $\Sigma'$ along $\mP_{B}(\Omega_{\mP(U_{B})}^{1}(1)|_{B})$
by \cite[Thm.2.3]{An}.
\end{proof}

\section{\textbf{Embedding theorem in the genus 5 case \label{sec:Embedding-theoremG5}}}

In this section, we treat the genus 5 case. The overall story is similar
to the one of the section \ref{sec:Common-prescription-for} though
details are different. We develop the discussion in this section while
keeping in mind the flow of discussion of the section \ref{sec:Common-prescription-for}.

\subsection{Extending the mid point \label{subsec:Extending-the-midG5}}

By \cite[Thm.6.5 (2) and Prop.7.8]{Mu2}, $W$ is a complete intersection
of three quadrics in $\mP^{6}$. Let $x_{1},\dots,x_{7}$ be coordinates
of $\mP^{6}$. We may assume that the plane $\Pi$ in $W$ is equal
to $\left\{ x_{1}=\cdots=x_{4}=0\right\} $. In this situation, the
equation of $W$ is of the following form:

\[
\left(\begin{array}{cccc}
l_{11} & l_{12} & l_{13} & l_{14}\\
l_{21} & l_{22} & l_{23} & l_{24}\\
l_{31} & l_{32} & l_{33} & l_{34}
\end{array}\right)\left(\begin{array}{c}
x_{1}\\
x_{2}\\
x_{3}\\
x_{4}
\end{array}\right)=\left(\begin{array}{c}
0\\
0\\
0
\end{array}\right),
\]
where $l_{ij}$ are linear forms of $x_{1},\dots,x_{7}$. 

Assume by contradiction that the dimension of the vector space generated
by the linear forms $l_{ij}(0,0,0,0,x_{5},x_{6},x_{7})$ is less than
or equal to $2$. Then $W$ is the cone over a complete intersection
of three quadrics in $\mP^{5}$ with a point $\mathsf{v}$ in $\Pi$
as the vertex. Then the Zariski tangent space of $W$ at $\mathsf{v}$
is dimension 6. This is absurd since $W$ has only Gorenstein terminal
singularities. Therefore, by a coordinate change keeping the equation
of $\Pi$ if necessary, we may assume that some three of $l_{ij}$,
say, $l_{i_{1}j_{1}},l_{i_{2}j_{2}},l_{i_{3}j_{3}}$ are equal to
$x_{5},x_{6},x_{7}$ respectively. 
\begin{defn}[Extension of $W$]
 In the projective space $\mP^{15}$ with coordinates $x_{1},\dots,x_{4}$
and $y_{ij}$ ($1\leq i\leq3,1\leq j\leq4)$, let $\overline{\Sigma}$
be the the following complete intersection of three quadrics:
\[
\overline{\Sigma}:=\left\{ [M_{y},\bm{x}]\in\mP^{15}\mid M_{y}\bm{x}=\bm{o}\right\} ,
\]
where 
\[
M_{y}:=\left(\begin{array}{cccc}
y_{11} & y_{12} & y_{13} & y_{14}\\
y_{21} & y_{22} & y_{23} & y_{24}\\
y_{31} & y_{32} & y_{33} & y_{34}
\end{array}\right),{{\empty}^{t}\!\bm{x}}=\begin{pmatrix}x_{1} & x_{2} & x_{3} & x_{4}\end{pmatrix}.
\]
 We set 
\[
\overline{\Pi}:=\left\{ x_{1}=x_{2}=x_{3}=x_{4}=0\right\} \subset\mP^{15}.
\]
\end{defn}

\begin{prop}
The pair $(W,\Pi)$ is projectively equivalent to the pair of a linear
section $W'$ of $\overline{\Sigma}$ and the $2$-plane $\overline{\Pi}\cap W'$.
\end{prop}

\begin{proof}
By writing $y_{i_{1}j_{1}},y_{i_{2}j_{2}},y_{i_{3}j_{3}}$ as $x_{5},x_{6},x_{7}$
respectively, it holds that the pair $(W,\Pi$) is projectively equivalent
to the pair of 
\[
W'=\overline{\Sigma}\cap\left\{ y_{ij}=l_{ij}\ \text{for}\ (i,j)\not=(i_{1},j_{1}),(i_{2},j_{2}),(i_{3},j_{3})\right\} ,
\]
and $\overline{\Pi}\cap W'.$ 
\end{proof}
We use the following notation:
\[
\overline{\Sigma}_{i}:=\left\{ [M_{y},\bm{o}]\in\overline{\Sigma}\mid\rank M_{y}\leq i\right\} \subset\overline{\Pi}.
\]

By elementary calculations, we determine the singular locus of $\overline{\Sigma}$
as follows:
\begin{prop}
\label{prop:The-singular-locus-G5} The singular locus of $\overline{\Sigma}$
is contained in $\overline{\Pi}$ and is equal to $\overline{\Sigma}_{2}$.
In particular, $\overline{\Sigma}$ is Gorenstein and normal.
\end{prop}

\begin{rem}
The variety $\overline{\Sigma}$ is an example of a variety of a complex
(cf.~\cite[Sect. 5]{Tan}).
\end{rem}

\subsection{Construction of the key variety \label{subsec:ConstructionKey G5}}
\begin{defn}
\label{Def:BundleG5} For $U^{3}\simeq\mC^{3}$ and $U^{4}\simeq\mC^{4}$,
we set 
\begin{eqnarray*}
\Sigma' & := & \mP_{\mP(U^{4})}(U^{3}\otimes\Omega_{\mP(U^{4})}^{1}(1)\oplus\sO_{\mP(U^{4})}(-1)).
\end{eqnarray*}
Note that, by a standard computation, it follows that $-K_{\Sigma'}=10H_{\Sigma}$.
\end{defn}

Under the situation of the subsection \ref{subsec:Extending-the-midG5},
we consider $\bm{x}$ as a coordinate vector of $U^{4}$ and $M_{y}$
as a coordinate matrix of $U^{3}\otimes(U^{4})^{*}.$ Then we can
regard $\overline{\Sigma}$ as a subvariety of $\mP($$U^{3}\otimes(U^{4})^{*}\oplus U^{4})$.
With this identification, we have the following proposition:
\begin{prop}
\label{Prop:G5crep} The following assertions hold:\vspace{3pt}

\noindent $(1)$ The tautological linear system $|H_{\Sigma'}|$
defines a surjective and birational morphism $\Sigma'\to\overline{\Sigma}$,
which we denote by $\varphi_{|H_{\Sigma'}|}$.

\vspace{3pt}

\noindent $(2)$ The morphism $\varphi_{|H_{\Sigma'}|}$ is an isomorphism
outside of $\overline{\Sigma}_{2}=\Sing\overline{\Sigma}$.

\vspace{3pt}

\noindent $(3)$ The $\varphi_{|H_{\Sigma'}|}$-fiber over a point
${\rm {\mathsf{t}}\in\overline{\Sigma}_{2}}$ is
\[
\begin{cases}
\mP^{1}: & {\mathsf{t}}\not\in\overline{\Sigma}_{1},\\
\mP^{2}: & {\rm {\mathsf{t}}\in\overline{\Sigma}_{1}.}
\end{cases}
\]
The morphism $\varphi_{|H_{\Sigma'}|}$ is a crepant small resolution.
\end{prop}

\begin{proof}
Let $\mathsf{{\rm p}}:=[U^{1}]\in\mP(U^{4})$ be a point, where $U^{1}\subset U^{4}$
is a $1$-dimensional vector space. The fiber of the projective bundle
$\Sigma'\to\mP(U^{4})$ over ${\rm p}$ is 
\[
\mP(U^{3}\otimes(U^{4}/U^{1})^{*}\oplus U^{1}),
\]
which is a linear subspace of $\mP(U^{3}\otimes(U^{4})^{*}\oplus U^{4})$.
By Lemma \ref{lem:ABFib} (1), the tautological linear system $|H_{\Sigma'}|$
defines a morphism $\Sigma'\to\mP(U^{3}\otimes(U^{4})^{*}\oplus U^{4})$.
By the description of fibers of $\Sigma'\to\mP(U^{4})$ as above and
the definition of $\overline{\Sigma}$, we see that the image of this
map is contained in $\overline{\Sigma}$. 

Let $\mathsf{t}=[M_{y},\bm{x}]$ be a point of $\overline{\Sigma}$.
By Lemma \ref{lem:ABFib} (2), the fiber of $\Sigma'\to\overline{\Sigma}$
over $\mathsf{t}$ is 
\[
\left\{ \mathsf{t}\right\} \times\left\{ [U^{1}]\mid M_{y}\in U^{3}\otimes(U^{4}/U^{1})^{*},\bm{x}\in U^{1}\right\} .
\]
If $\mathsf{t}\not\in\overline{\Pi}$, namely, $\bm{x}\not=\bm{o}$,
then $U^{1}$ are uniquely determined as $U^{1}=\mC\bm{x}$ (since
$\mathsf{t}\in\overline{\Sigma},$ it holds that $M_{y}\in U^{3}\otimes(U^{4}/\mC\bm{x})^{*})$.
Therefore the morphism $\Sigma'\to\overline{\Sigma}$ is an isomorphism
outside of $\overline{\Pi}$. In particular, the morphism $\Sigma'\to\overline{\Sigma}$
is birational. Assume that $\mathsf{t}\in\overline{\Pi},$ equivalently,
$\bm{x=o}.$ The condition for $U^{1}$ is that $U^{1}\subset\left\{ \bm{x}\in U^{4}\mid M_{y}\bm{x}=\bm{0}\right\} \simeq\mC^{4-\rank M_{y}}$.
Therefore the fiber of $\Sigma'\to\overline{\Sigma}$ over $\mathsf{t}$
is isomorphic to $\mP^{3-\rank M_{y}}$. From this, the description
of the fiber $\Sigma'\to\overline{\Sigma}$ over $\mathsf{t}$ follows.
The morphism $\Sigma'\to\overline{\Sigma}$ is crepant since it holds
that $-K_{\Sigma'}=10H_{\Sigma'}$.
\end{proof}
Let $\Sigma'_{1}\subset\Sigma'$ be the inverse image of $\overline{\Sigma}_{1}$.
Since $\overline{\Sigma}_{1}=\mP(U^{3})\times\mP((U^{4})^{*})$, we
see that 
\begin{equation}
\Sigma_{1}'=\mP(U^{3})\times\mP(\Omega_{\mP(U^{4})}^{1}(1)).\label{eq:Sigma'1}
\end{equation}
 Let $\tau\colon\widehat{\Sigma}\to\Sigma'$ be the blow-up of $\Sigma'$
along $\Sigma'_{1}$ , and $E_{\widehat{\Sigma}}$ the $\tau$-exceptional
divisor. We denote by $\Pi'$ and $\widehat{\Pi}$ the strict transforms
of $\overline{\Pi}$ on $\Sigma'$ and $\widehat{\Sigma}$ respectively.
Note that 
\begin{equation}
\Pi'=\mP(U^{3}\otimes\Omega_{\mP(U^{4})}^{1}(1)\oplus0).\label{eq:Pi'}
\end{equation}

\begin{lem}
\label{Lem:Et} Let $\mathsf{t}=[U^{1}]\times[W^{1}]$ be a point
of $\overline{\Sigma}_{1}=\mP(U^{3})\times\mP((U^{4})^{*})$, where
$U^{1}$ and $W^{1}$ are $1$-dimensional subspaces of $U^{3}$ and
$(U^{4})^{*}$ respectively. The following assertions hold:\vspace{3pt}

\noindent $(1)$ The fiber of $\Sigma'_{1}\to\overline{\Sigma}_{1}$
over $\mathsf{t}$ can be identified with the fiber of $\mP(\Omega_{\mP(U^{4})}^{1}(1))\to\mP((U^{4})^{*})$
over the point $[W^{1}]$ and then with $\mP(W^{1,\perp})$, where
$W^{1,\perp}$ is the subspace of $U^{4}$ orthogonal to $W^{1}$
with respect to the dual pairing. Let $E_{\mathsf{t}}$ be the fiber
of $E_{\widehat{\Sigma}}\to\overline{\Sigma}_{1}$ over $\mathsf{t}$.
It holds that 
\[
E_{\mathsf{t}}=\mP_{\mP(W^{1,\perp})}((U^{3}/U^{1})\otimes\Omega_{\mP(W^{1,\perp})}^{1}(1)\oplus\sO_{\mP(W^{1,\perp})}(-1)),
\]
and 
\[
E_{\mathsf{t}}\cap\widehat{\Pi}=\mP_{\mP(W^{1,\perp})}((U^{3}/U^{1})\otimes\Omega_{\mP(W^{1,\perp})}^{1}(1)\oplus0).
\]
 \vspace{3pt}

\noindent $(2)$ We identify an element of $(U^{3}/U^{1})\otimes(W^{1,\perp})^{*}$
with a $2\times3$ matrix. The linear system $|H_{E_{\mathsf{t}}}|$
defines a morphism $E_{\mathsf{t}}\to\mP((U^{3}/U^{1})\otimes(W^{1,\perp})^{*}\oplus W^{1,\perp})$,
and the image is 
\[
\overline{E}_{\mathsf{t}}:=\left\{ [M,\bm{x}]\mid M\in(U^{3}/U^{1})\otimes(W^{1,\perp})^{*},\bm{x}\in W^{1,\perp},M\bm{x}=\bm{o}\right\} ,
\]
which is a complete intersection of two quadrics.

\vspace{3pt}

\noindent $(3)$ The singular locus of $\overline{E}_{\mathsf{t}}$
is 
\[
\left\{ [M,\bm{o}]\mid M\in(U^{3}/U^{1})\otimes(W^{1,\perp})^{*},\rank M\leq1\right\} ,
\]
 which is $\mP(U^{3}/U^{1})\times\mP((W^{1,\perp})^{*})\simeq\mP^{1}\times\mP^{2}$.
The morphism $E_{\mathsf{t}}\to\overline{E}_{\mathsf{t}}$ is an isomorphism
outside of $\Sing\overline{E}_{\mathsf{t}}$, and the fiber over a
point of ${\rm \Sing}\overline{E}_{\mathsf{t}}$ is $\mP^{1}$. The
morphism $E_{\mathsf{t}}\to\overline{E}_{\mathsf{t}}$ is a crepant
small resolution.

\vspace{3pt}

\noindent $(4)$ The induced morphism $E_{\mathsf{t}}\cap\widehat{\Pi}\to\mP((U^{3}/U^{1})\otimes(W^{1,\perp})^{*}\oplus0)\simeq\mP^{5}$
is the blow-up of $\mP^{5}$ along $\Sing\overline{E}_{\mathsf{t}}\simeq\mP^{1}\times\mP^{2}$.
Let $L_{\mP(W^{1,\perp})}$ be the pull-back to $E_{\mathsf{t}}\cap\widehat{\Pi}$
of a line in $\mP(W^{1,\perp})\simeq\mP^{2}$. The exceptional divisor
of the blow-up of $\mP^{5}$ along $\mP^{1}\times\mP^{2}$ is linearly
equivalent to $2H_{E_{\mathsf{t}}\cap\widehat{\Pi}}-L_{\mP(W^{1,\perp})}$.
\end{lem}

\begin{proof}
We show the assertion (1). The first assertion of (1) easily follows
from (\ref{eq:Sigma'1}). For the second assertion of (2), we have
only to determine the restriction to $\mP(W^{1,\perp})$ of the normal
bundle $\sN_{\Sigma'_{1}/\Sigma'}$. Since $\Sigma_{1}'$ is a sub
$\mP^{2}\times\mP^{2}$-bundle of the $\mP^{8}$-bundle $\Pi'$ over
$\mP(U^{4})$ by (\ref{eq:Sigma'1}) and (\ref{eq:Pi'}), we see that
\begin{equation}
\sN_{\Sigma_{1}'/\Pi'}\simeq T_{\mP(U^{3})}\otimes T_{\mP(\Omega_{\mP(U^{4})}^{1}(1))/\mP(U^{4})}\label{eq:NSigma1}
\end{equation}
 relativising the normal bundle of the Segre embedded $\mP^{2}\times\mP^{2}$
in $\mP^{8}$. Let $p\colon\mP(\Omega_{\mP(U^{4})}^{1}(1))\to\mP(U^{4})$
be the natural morphism. We consider $\mP(W^{1,\perp})$ as the $p$-fiber
over $[W^{1}]\in\mP((U^{4})^{*})$. Restricting to $\mP(W^{1,\perp})$
the relative Euler sequence 
\[
0\to\sO(-H_{\mP(\Omega_{\mP(U^{4})}^{1}(1))})\to p^{*}\Omega_{\mP(U^{4})}^{1}(1)\to T_{\mP(\Omega_{\mP(U^{4})}^{1}(1))/\mP(U^{4})}(-H_{\mP(\Omega_{\mP(U^{4})}^{1}(1))})\to0,
\]
we obtain the exact sequence
\[
0\to\sO_{\mP(W^{1,\perp})}\to\Omega_{\mP(W^{1,\perp})}^{1}(1)\oplus\sO_{\mP(W^{1,\perp})}\to T_{\mP(\Omega_{\mP(U^{4})}^{1}(1))/\mP(U^{4})}|_{\mP(W^{1,\perp})}\to0
\]
 since $H_{\mP(\Omega_{\mP(U^{4})}^{1}(1))}|_{\mP(W^{1,\perp})}=0$.
Therefore we have $T_{\mP(\Omega_{\mP(U^{4})}^{1}(1))/\mP(U^{4})}|_{\mP(W^{1,\perp})}\simeq\Omega_{\mP(W^{1,\perp})}^{1}(1)$.
We also note that $T_{\mP(U^{3})}|_{[U^{1}]}\simeq(U^{3}/U^{1})\otimes(U^{1})^{*}$.
Hence, by (\ref{eq:NSigma1}), we obtain 
\[
\sN_{\Sigma_{1}'/\Pi'}|_{\mP(W^{1,\perp})}\simeq(U^{3}/U^{1})\otimes(U^{1})^{*}\otimes\Omega_{\mP(W^{1,\perp})}^{1}(1)\simeq(U^{3}/U^{1})\otimes\Omega_{\mP(W^{1,\perp})}^{1}(1).
\]
Let $L_{\mP(U^{4})}$ be the pull-back to $\Sigma'$ of a hyperplane
of $\mP(U^{4})$. Since 
\begin{equation}
\Pi'\sim H_{\Sigma'}-L_{\mP(U^{4})},\label{eq:Pi'sim}
\end{equation}
 we have $\sN_{\Pi'/\Sigma'}|_{\mP(W^{1,\perp})}\simeq\sO_{\mP(W^{1,\perp})}(-1)$.
Therefore, by the normal bundle sequence $0\to\sN_{\Sigma_{1}'/\Pi'}\to\sN_{\Sigma_{1}'/\Sigma'}\to\sN_{\Pi'/\Sigma'}|_{\Sigma'_{1}}\to0$,
we obtain 
\[
\sN_{\Sigma_{1}'/\Sigma'}|_{\mP(W^{1,\perp})}\simeq(U^{3}/U^{1})\otimes\Omega_{\mP(W^{1,\perp})}^{1}(1)\oplus\sO_{\mP(W^{1,\perp})}(-1),
\]
and hence the assertion (1) follows.

The assertions (2)--(4) can be proved in a similar way to Proposition
\ref{Prop:G5crep} due to structural similarity between $E_{\mathsf{t}}$
and $\Sigma'$, so we omit a proof. 
\end{proof}
Note that
\begin{equation}
-K_{\widehat{\Sigma}}=\tau^{*}(-K_{\Sigma'})-4E_{\widehat{\Sigma}}=10\tau^{*}H_{\Sigma'}-4E_{\widehat{\Sigma}}=2\tau^{*}H_{\Sigma'}+4(2\tau^{*}H_{\Sigma'}-E_{\widehat{\Sigma}}).\label{eq:G5weakFano}
\end{equation}

By Proposition \ref{Prop:G5crep}, $H_{\Sigma'}$ is nef and big and,
since $\overline{\Sigma}_{1}$ is the intersection of quadrics, $\Bs|2\tau^{*}H_{\Sigma'}-E_{\widehat{\Sigma}}|=\emptyset$.
Therefore $-K_{\widehat{\Sigma}}$ is nef and big. Let $\nu\colon\widehat{\Sigma}\to\widehat{\Sigma}'$
be the anti-canonical model.

\vspace{3pt}

\noindent\textbf{ Flop $\boldmath{\widehat{\Sigma}\dashrightarrow\Sigma^{+}}$.}
\begin{prop}
\label{prop:SigmaflopG5} The anti-canonical model $\nu\colon\widehat{\Sigma}\to\widehat{\Sigma}'$
is defined over $\overline{\Sigma}$ and is a $\widehat{\Pi}$-negative
flopping contraction of Atiyah type. The morphism $\nu|_{E_{\widehat{\Sigma}}}$
is also a flopping contraction of Atiyah type.
\end{prop}

\begin{proof}
Let $l\subset\widehat{\Sigma}$ be an irreducible $\nu$-exceptional
curve. By (\ref{eq:G5weakFano}) and $-K_{\widehat{\Sigma}}\cdot l=0$,
we have $\tau^{*}H_{\Sigma'}\cdot l=(2\tau^{*}H_{\Sigma'}-E_{\widehat{\Sigma}})\cdot l=0,$
and hence $\tau^{*}H_{\Sigma'}\cdot l=E_{\widehat{\Sigma}}\cdot l=0.$
By $\tau^{*}H_{\Sigma'}\cdot l=0$, the morphism $\nu$ is defined
over $\overline{\Sigma}$. Since $\Pi'$ is smooth, we have $\widehat{\Pi}=\tau^{*}\Pi'-E_{\widehat{\Sigma}}$.
Therefore $\widehat{\Pi}$ is negative for any $\nu$-exceptional
curve since $\Pi'$ is negative for any exceptional curve for $\Sigma'\to\overline{\Sigma}$
by (\ref{eq:Pi'sim}). By Proposition \ref{Prop:G5crep} (2), $l$
is contained in the union of $E_{\widehat{\Sigma}}$ and the strict
transform $\widehat{\Sigma}_{2}$ of $\overline{\Sigma}_{2}$. By
Proposition \ref{Prop:G5crep} (3), the $\nu$-fiber over a point
${\rm {\mathsf{s}}\in\widehat{\Sigma}_{2}\setminus E_{\widehat{\Sigma}}}$
is $\mP^{1}$. Note that $E_{\mathsf{t}}$ as in Lemma \ref{Lem:Et}
(1) has two nontrivial contractions; one is the morphism $E_{\mathsf{t}}\to\mP(W^{1,\perp})$,
and nontrivial fibers of another morphism are $\mP^{1}$ by Lemma
\ref{Lem:Et} (3). Since $\nu|_{E_{\mathsf{t}}}$ cannot be the morphism
$E_{\mathsf{t}}\to\mP(W^{1,\perp})$, we see that the nontrivial $\nu$-fiber
over a point ${\rm {\mathsf{s}}\in E_{\widehat{\Sigma}}}$ is also
$\mP^{1}$. Therefore any nontrivial fiber of the morphism $\nu|_{\widehat{\Pi}}\colon\widehat{\Pi}\to\nu(\widehat{\Pi})$
is $\mP^{1}$. In particular, this implies that the relative Picard
number of $\nu|_{\widehat{\Sigma}}$ is one. Note that, since $\widehat{\Pi}$
is $\nu$-negative , $-K_{\widehat{\Pi}}$ is $\nu$-ample. Therefore,
by \cite[Thm.2.3]{An}, $\nu(\widehat{\Pi})$ is smooth and $\nu|_{\widehat{\Pi}}$
is the blow-up of $\nu(\widehat{\Pi})$ along a smooth subvariety
of $\nu(\widehat{\Pi})$ which is the strict transform of $\overline{\Sigma}_{2}$.
This implies that $\sN_{l/\widehat{\Pi}}\simeq\sO_{\mP^{1}}^{\oplus9}\oplus\sO_{\mP^{1}}(-1)$,
and $\widehat{\Pi}\cdot l=K_{\widehat{\Pi}}\cdot l=-1$. Therefore,
by the normal bundle sequence $0\to\sN_{l/\widehat{\Pi}}\to\sN_{l/\widehat{\Sigma}}\to\sN_{\widehat{\Pi}/\widehat{\Sigma}}|_{l}\to0$,
we have $\sN_{l/\widehat{\Sigma}}\simeq\sO_{\mP^{1}}^{\oplus9}\oplus\sO_{\mP^{1}}(-1)^{\oplus2}$,
and hence $\nu$ is a flopping contraction of Atiyah type.

Since $E_{\widehat{\Sigma}}\cdot l=0$, $\nu|_{E_{\widehat{\Sigma}}}$
is also a flopping contraction of Atiyah type.
\end{proof}
Let $\widehat{\Sigma}\dashrightarrow\Sigma^{+}$ be the flop for the
flopping contraction $\nu$. It is well-known that the flop can be
constructed by the blow-up along the $\nu$-exceptional locus and
the blow-down of the exceptional divisor along the other direction.
We denote by $\Pi^{+}$ and $E_{\Sigma^{+}}$ the strict transforms
on $\Sigma^{+}$ of $\widehat{\Pi}$ and $E_{\widehat{\Sigma}}$ respectively.
By the construction of $\widehat{\Sigma}\dashrightarrow\Sigma^{+}$,
we see that the induced map $E_{\widehat{\Sigma}}\dashrightarrow E_{\Sigma^{+}}$
is also the flop and the induced map $\widehat{\Pi}\dashrightarrow\Pi^{+}$
is identified with $\nu|_{\widehat{\Pi}}\colon\widehat{\Pi}\to\nu(\widehat{\Pi})$,
which is the blow-up of $\nu(\widehat{\Pi})$ along the strict transform
of $\overline{\Sigma}_{2}$ on $\nu(\widehat{\Pi})$ by the proof
of Proposition \ref{prop:SigmaflopG5}.

In the following steps, we separate $E_{\Sigma^{+}}$ and $\Pi^{+}$
by a flip, and finally contract their strict transforms.

\vspace{5pt}

\noindent\textbf{ Flip $\boldmath{{\Sigma}^{+}\dashrightarrow\widetilde{\Sigma}}$.}

We set 
\[
G:=E_{\Sigma^{+}}\cap\Pi^{+}.
\]

\begin{lem}
\label{Lem:IntEPiG5} The exceptional locus of $\Pi^{+}\to\overline{\Pi}$
is $G$ and $G$ is a $\mP^{5}$-bundle over $\mP(U^{3})\times\mP((U^{4})^{*})$. 
\end{lem}

\begin{proof}
As we note above, we may identify $\Pi^{+}\to\overline{\Pi}$ with
$\nu(\widehat{\Pi})\to\overline{\Pi}$. Therefore the assertion follows
by Lemma \ref{Lem:Et} (4) and the construction of the flop.
\end{proof}
\begin{prop}
\label{Prop:flippingG5} 
\end{prop}

\begin{enumerate}
\item Let $\Gamma\simeq\mP^{5}$ be a fiber of $G\to\mP(U^{3})\times\mP((U^{4})^{*})$.
It holds that
\[
\sN_{\Gamma/\Sigma^{+}}=\sO_{\mP^{5}}(-1)^{\oplus2}\oplus\sO_{\mP^{5}}^{\oplus5}.
\]
\item There exists a small contraction $\Sigma^{+}\to\overline{\Sigma}^{+}$
contracting $E_{\Sigma^{+}}\cap\Pi^{+}$ onto $\mP(U^{3})\times\mP((U^{4})^{*})$. 
\end{enumerate}
\begin{proof}
(1). Let $l$ be a general line in a fiber of the $\mP^{3}$-bundle
$E_{\mathsf{t}}\cap\widehat{\Pi}\to\mP(W^{1,\perp})$ as in Lemma
\ref{Lem:Et} (1) and $l^{+}$ the strict transform of $l$ on $\Sigma^{+}$.
Since $l$ is contained in a fiber of the blow-up $\tau\colon\widehat{\Sigma}\to\Sigma'$
and $E_{\widehat{\Sigma}}$ is the $\tau$-exceptional divisor, we
have $-K_{\widehat{\Sigma}}\cdot l=4$ and $E_{\widehat{\Sigma}}\cdot l=-1$.
Since both $-K_{\widehat{\Sigma}}$ and $E_{\widehat{\Sigma}}$ are
numerically trivial for flopping curves by the proof of Proposition
\ref{prop:SigmaflopG5}, we have $-K_{\Sigma^{+}}\cdot l^{+}=4$ and
$E_{\Sigma^{+}}\cdot l^{+}=-1$. From the latter equality, we have
$(G\cdot l^{+})_{\Pi^{+}}=E_{\Pi^{+}}\cdot l^{+}=-1$. Therefore,
by the normal bundle sequence $0\to\sN_{\Gamma/G}\to\sN_{\Gamma/\Pi^{+}}\to\sN_{G/\Pi^{+}}|_{\Gamma}\to0$
and $\sN_{\Gamma/G}\simeq\sO_{\Gamma}^{\oplus5}$, we have $\sN_{\Gamma/\Pi^{+}}\simeq\sO_{\mP^{5}}(-1)\oplus\sO_{\mP^{5}}^{\oplus5}$.
Since $-K_{\Sigma^{+}}\cdot l^{+}=4$ and $\Gamma\simeq\mP^{5}$,
we have $\deg\sN_{\Gamma/\Sigma^{+}}=-2$. Therefore, by the normal
bundle sequence $0\to\sN_{\Gamma/\Pi^{+}}\to\sN_{\Gamma/\Sigma^{+}}\to\sN_{\Pi^{+}/\Sigma^{+}}|_{\Gamma}\to0$,
we have $\sN_{\Gamma/\Sigma^{+}}\simeq\sO_{\mP^{5}}(-1)^{\oplus2}\oplus\sO_{\mP^{5}}^{\oplus5}$
as desired.

\noindent (2). Let $H_{\Sigma^{+}}$ be the strict transform on $\Sigma^{+}$
of $\tau^{*}H_{\Sigma^{'}}$. We can show that $(2H_{\Sigma^{+}}-E_{\Sigma^{+}})+\Pi^{+}$
is nef over $\overline{\Sigma}$ and numerically trivial only for
fibers of $G\to\mP(U^{3})\times\mP((U^{4})^{*})$. A proof is quite
similar to the one of Proposition \ref{Prop:flipping G46} (2), so
we omit it. 
\end{proof}
By Proposition \ref{Prop:flippingG5} (1), the contraction $\Sigma^{+}\to\overline{\Sigma}^{+}$
is of flipping type, and the flip can be constructed by the blow-up
along $\Lambda$ and the blow-down of the exceptional divisor along
the other direction . Let $\Sigma^{+}\dashrightarrow\widetilde{\Sigma}$
be the flip. By Proposition \ref{Prop:flippingG5} (1) again, the
flipped locus is a $\mP^{1}$-bundle over $\mP(U^{3})\times\mP((U^{4})^{*})$.
We denote by $E_{\widetilde{\Sigma}}$, $\widetilde{\Pi}$ and $H_{\widetilde{\Sigma}}$
the strict transforms on $\widetilde{\Sigma}$ of $E_{\Sigma^{+}}$,
$\Pi^{+}$ and $H_{\Sigma^{+}}$ respectively. 

The following lemma will describe a part of singularities of the key
variety $\Sigma$ which we are going to construct: 
\begin{lem}
\label{Lem:G25} Let $\mathsf{t}=[U^{1}]\times[W^{1}]$ be a point
of $\overline{\Sigma}_{1}=\mP(U^{3})\times\mP((U^{4})^{*})$ as in
Lemma \ref{Lem:Et}. The fiber $E_{\mathsf{t}}^{+}$ of $E_{\Sigma^{+}}\to\mP(U^{3})\times\mP((U^{4})^{*})$
over $\mathsf{t}$ is the blow-up of the Grassmannian $\rG(2,(U^{3}/U^{1})\oplus((W^{1})^{\perp})^{*})\simeq\rG(2,5)$
at the point $[\wedge^{2}(U^{3}/U^{1})]$. The fiber of $E_{\widetilde{\Sigma}}\to\mP(U^{3})\times\mP((U^{4})^{*})$
over $\mathsf{t}$ is $\rG(2,(U^{3}/U^{1})\oplus((W^{1})^{\perp})^{*})$. 
\end{lem}

\begin{proof}
For simplicity of notation, we set $\overline{U}^{2}:=U^{3}/U^{1},$${\rm Gr}:=\rG(2,\overline{U}^{2}\oplus(W^{1,\perp})^{*})$
and denote by $\widetilde{{\rm Gr}}$ the blow-up of ${\rm Gr}$ at
the point $[\wedge^{2}\overline{U}^{2}].$ By Lemma \ref{Lem:Et},
the fiber $E_{\mathsf{t}}$ of $E_{\widehat{\Sigma}}\to\mP(U^{3})\times\mP((U^{4})^{*})$
over $\mathsf{t}$ has two contractions, one of which is the $\mP^{4}$-bundle
$E_{\mathsf{t}}\to\mP(W^{1,\perp})$ and another of which is the flopping
contraction of Atiyah type $E_{\mathsf{t}}\to\overline{E}_{\mathsf{t}}$.
Note that $\overline{E}_{\mathsf{t}}$ has another unique small resolution
different from $E_{\mathsf{t}}\to\overline{E}_{\mathsf{t}}$, which
we denote by $E_{\mathsf{t}}^{+}\to\overline{E}_{\mathsf{t}}$. Therefore,
to show $E_{\mathsf{t}}^{+}\simeq\widetilde{{\rm Gr}}$, it suffices
to show that $\widetilde{{\rm Gr}}$ has a small contraction onto
$\overline{E}_{\mathsf{t}}$ (note that $\widetilde{{\rm Gr}}$ is
different from $E_{\mathsf{t}}$ since $E_{\mathsf{t}}$ has no contraction
onto $\rG(2,5)$). We note the following decomposition: 
\begin{align*}
\wedge^{2}\left(\overline{U}^{2}\oplus(W^{1,\perp})^{*}\right) & =\wedge^{2}\overline{U}^{2}\oplus\overline{U}^{2}\wedge(W^{1,\perp})^{*}\oplus\wedge^{2}(W^{1,\perp})^{*}\\
 & \simeq\wedge^{2}\overline{U}^{2}\oplus\overline{U}^{2}\otimes(W^{1,\perp})^{*}\oplus W^{1,\perp}.
\end{align*}
Therefore, the linear projection from the point $[\wedge^{2}\overline{U}^{2}]$
maps ${\rm Gr}$ into the projective space $\mP\left(\overline{U}^{2}\wedge(W^{1,\perp})^{*}\oplus W^{1,\perp}\right).$
Let $\bm{e}_{1},\bm{e}_{2}$ and $\bm{e}_{3},\bm{e}_{4},\bm{e}_{5}$
be basis of $\overline{U}^{2}$ and $(W^{1,\perp})^{*}$ respectively,
and $p_{ij}$ the Pl\"ucker coordinates associated to the basis $\bm{e}_{1},\dots,\bm{e}_{5}$
of $\overline{U}^{2}\oplus(W^{1,\perp})^{*}$. The equation of ${\rm Gr}$
is given by $\wedge^{2}(\sum_{1\leq i<j\leq5}p_{ij}\bm{e}_{i}\wedge\bm{e}_{j})=\bm{o}$.
We can check explicitly that the image of the projection of ${\rm Gr}$
is defined by 
\begin{equation}
\left(\sum_{i=1,2,j=3,4,5}p_{ij}\bm{e}_{i}\wedge\bm{e}_{j}\right)\wedge\left(\sum_{3\leq i<j\leq5}p_{ij}\bm{e}_{i}\wedge\bm{e}_{j}\right)=\bm{o}\label{eq:22ci}
\end{equation}
 in $\mP\left(\wedge^{2}\overline{U}^{2}\oplus\overline{U}^{2}\wedge(W^{1,\perp})^{*}\oplus\wedge^{2}(W^{1,\perp})^{*}\right)$,
and the projection is birational onto the image. We identify $\wedge^{2}(W^{1,\perp})^{*}$
with $W^{1,\perp}$ by regarding $\bm{e}_{3}\wedge\bm{e}_{4}$, $\bm{e}_{3}\wedge\bm{e}_{5}$,
and $\bm{e}_{4}\wedge\bm{e}_{5}$ as $\bm{e}_{5}^{*}$, $-\bm{e}_{4}^{*}$,
and $\bm{e}_{3}^{*}\in W^{1,\perp}$ respectively, and also $\overline{U}^{2}\wedge(W^{1,\perp})^{*}$
with $\overline{U}^{2}\otimes(W^{1,\perp})^{*}$ by regarding $\bm{e}_{i}\wedge\bm{e}_{j}$with
$\bm{e}_{i}\otimes\bm{e}_{j}.$ Then the equation (\ref{eq:22ci})
defines $\overline{E}_{\mathsf{t}}$ in $\mP\left(\overline{U}^{2}\wedge(W^{1,\perp})^{*}\oplus W^{1,\perp}\right).$
By a standard property of linear projection, a natural morphism $\rho^{+}\colon\widetilde{{\rm Gr}}\to\overline{E}_{\mathsf{t}}$
is induced. Let $E_{{\rm Gr}}$ be the exceptional divisor of the
blow-up $\widetilde{{\rm Gr}}\to{\rm Gr},$ and $L_{{\rm Gr}}$ the
total transform on $\widetilde{{\rm Gr}}$ of a hyperplane section
of ${\rm Gr}$. We have $-K_{\widetilde{{\rm Gr}}}=5(L_{{\rm Gr}}-E_{{\rm Gr}})$
and $L_{{\rm Gr}}-E_{{\rm Gr}}$ is the total transform of a hyperplane
section of $\overline{E}_{\mathsf{t}}$. Therefore the morphism $\rho^{+}$
is crepant, and hence must be small since $\overline{E}_{\mathsf{t}}$
has only terminal singularities. Now we have shown that $\widetilde{{\rm Gr}}$
has a small contraction onto $\overline{E}_{\mathsf{t}}$ as desired. 

The description of the fiber of $E_{\widetilde{\Sigma}}\to\mP(U^{3})\times\mP((U^{4})^{*})$
over $\mathsf{t}$ follows from the above by the construction of the
flip $\Sigma^{+}\dashrightarrow\widetilde{\Sigma}$.
\end{proof}
\vspace{5pt}

\noindent \textbf{Contracting $E_{\widetilde{\Sigma}}$ and $\widetilde{\Pi}$.} 

\vspace{5pt}

\noindent \noindent By the constructions of the flop $\widehat{\Sigma}\dashrightarrow\Sigma^{+}$
and the flip $\Sigma^{+}\dashrightarrow\widetilde{\Sigma}$, and the
description of $E_{\Sigma^{+}}\cap\Pi^{+}$ as in Lemma \ref{Lem:IntEPiG5},
we see that $E_{\widetilde{\Sigma}}\cap\widetilde{\Pi}=\emptyset$.

By the construction of the flip, we see that $\Sigma^{+}\dashrightarrow\widetilde{\Sigma}$
induces the contraction $\Pi^{+}\to\mP(U^{3}\otimes(U^{4})^{*})$.
Thus $\widetilde{\Pi}\simeq\mP^{11}$.
\begin{lem}
\label{Lem:Pi|PiG5} The normal bundle $\sN_{\widetilde{\Pi}/\widetilde{\Sigma}}$
is $\sO_{\mP^{11}}(-2)$ .
\end{lem}

\begin{proof}
A proof is similar to those of Lemmas \ref{lem:P/PG8} and \ref{lem:F|F G46},
so we omit it. 
\end{proof}
\begin{lem}
\label{Lem:2H+PiG5} $2H_{\widetilde{\Sigma}}+\widetilde{\Pi}$ is
semiample. 
\end{lem}

\begin{proof}
A proof is quite similar to the one of Lemma \ref{Lem:2H+F0}, so
we omit it. 
\end{proof}
\begin{thm}
\label{Thm:keyg5} Let $\mu\colon\widetilde{\Sigma}\to\Sigma$ be
the contraction defined by a sufficient multiple of $2H_{\widetilde{\Sigma}}+\widetilde{\Pi}$.
The following assertions hold:

\vspace{3pt}

\noindent $(1)$ The $\mu$-exceptional locus is the union of $E_{\widetilde{\Sigma}}$
and $\widetilde{\Pi}$. 

\vspace{3pt}

\noindent $(2)$ $\mu(\widetilde{\Pi})$ is a $1/2$-singularity.

\vspace{3pt}

\noindent $(3)$ The discrepancy of $E_{\widetilde{\Sigma}}$ is
$4$ and $\mu(E_{\widetilde{\Sigma}})\simeq\mP(U^{3})\times\mP((U^{4})^{*})$.
Any fiber of $E_{\widetilde{\Sigma}}\to\mu(E_{\widetilde{\Sigma}})$
is $\rG(2,5).$ In particular $\Sigma$ has Gorenstein terminal singularities
along $\mP(U^{3})\times\mP((U^{4})^{*})$.

\vspace{3pt}

\noindent $(4)$ The variety $\Sigma$ is a $12$-dimensional rational
$\mQ$-Fano variety with $\rho(\Sigma)=1$. 

\vspace{3pt}

\noindent $(5)$ The image $M_{\Sigma}$ of $H_{\Sigma}$ is a primitive
integral ample Weil divisor $M_{\Sigma}$ and it holds that $-K_{\Sigma}=10M_{\Sigma}$. 
\end{thm}

\begin{proof}
The assertion (3) follows from Lemma \ref{Lem:G25}. The rest assertion
can be proved similarly to Theorem \ref{thm:common}, so we omit a
proof.
\end{proof}

\subsection{Embedding theorem \label{subsec:Embedding-theorem G5}}

Now we arrive at Theorem \ref{thm:main1} for a prime $\mQ$-Fano
$3$-fold $X$ of genus 5. 
\begin{thm}
\label{thm:embg5} A $\mQ$-Fano $3$-fold $X$ of genus $5$ is a
linear section of $\Sigma$. 
\end{thm}

\begin{proof}
We only remark that $W$ is disjoint from $\overline{\Sigma}_{1}$
since $\overline{\Sigma}$ has non-hypersurface singularities along
$\overline{\Sigma}_{1}$. The rest of the proof is similar to the
one of Theorem \ref{thm:embg4g6}, so we omit it.
\end{proof}

\subsection{Extension of the Sarkisov link}

We have obtained the following diagram: 

\begin{equation}\label{eq:SarkisovG5} \xymatrix{& \widetilde{\Sigma}\ar@{-->}^{\text{anti-flip}\ }[r]\ar[dl]_{\mu}\ar[dr] & \Sigma^+\ar@{-->}^{\text{flop}}[r]\ar[d]& \widehat{\Sigma}\ar[dr]\ar[dl]\\
\Sigma & &\overline{\Sigma} & & \mP^3.}
\end{equation}

By the proof of Theorem \ref{thm:embg5}, we obtain the following:

\begin{cor}
\label{cor:ExtSarkisovG5}The diagram (\ref{eq:SarkisovG5}) is an
extension of (\ref{eq:Sarkisov}) in the case of genus 5, where $Y$,$W$
and $Y'$ are linear sections of $\widetilde{{\Sigma}}$, $\overline{\Sigma}$
and $\widehat{\Sigma}$ with respect to $|H_{\widetilde{\Sigma}}|$,
$|\sO_{\overline{\Sigma}}(1)|$ and $|H_{\widehat{\Sigma}}|$ respectively.
The restriction of the anti-flip to $Y$ is the identity. 
\end{cor}

\begin{proof}
The restriction of the anti-flip to $Y$ is the identity since the
image in $\overline{\Sigma}$ of the flipped locus in $\widetilde{\Sigma}$
is contained in the image of $E_{\widetilde{\Sigma}}$ on $\overline{\Sigma}$
by Proposition \ref{Prop:flippingG5} (2), and the latter is disjoint
from $W$ as we remarked in the proof of Theorem \ref{thm:embg5}.
The rest follows easily.
\end{proof}


\begin{thebibliography}{Fuj1}
\bibitem[An]{An} T.~Ando, \textit{On extremal rays of the higher-dimensional
varieties}, Invent. Math. 81 (1985), no. 2, 347--357. 

\bibitem[CKM]{CKM} H.~Clemens, J.~Koll\'ar, S.~Mori, \textit{Higher-dimensional
complex geometry}, Ast\'erisque No. 166 (1988), 144 pp. (1989).

\bibitem[Fuj1]{Fuj1} T.~Fujita, \textit{On the structure of polarized
manifolds with total deficiency one. I}, J. Math. Soc. Japan 32 (1980),
no. 4, 709--725. 

\bibitem[Fuj2]{Fuj2} T.~Fujita, \textit{On the structure of polarized
manifolds with total deficiency one. II}, J. Math. Soc. Japan 33 (1981),
no. 3, 415--434.

\bibitem[Fuj3]{Fuj3}T.~Fujita, \textit{Projective varieties of $\Delta$-genus
one}, Algebraic and topological theories (Kinosaki, 1984), 149--175,
Kinokuniya, Tokyo, 1986. 

\bibitem[Fuk]{Fuk} T.~Fukuoka, \textit{Relative linear extensions
of sextic del Pezzo fibrations over curves}, Ann. Sc. Norm. Super.
Pisa Cl. Sci. (5) Vol. XXI (2020), 1371--1409.

\bibitem[Gu]{Gu} N.~P.~Gushel\textasciiacute , \textit{Fano varieties
of genus 8}, Uspekhi Mat. Nauk 38 (1983), no. 1(229), 163--164

\bibitem[Ha]{Ha}U.~Hayat, \textit{The Cramer varieties ${\rm Cr}(r,r+s,s)$},
J. Geom. Phys. 79 (2014), 53--58. 

\bibitem[Il]{Il} A.~Iliev,\textit{ Lines on the Gushel\textasciiacute{}
threefold,} Indag. Math. (N.S.) 5 (1994), no. 3, 307--320. 

\bibitem[Is]{Is}V.~A.~Iskovskih, \textit{Fano threefolds. II},
Izv. Akad. Nauk SSSR Ser. Mat. 42 (1978), no. 3, 506--549. 

\bibitem[IsP]{IsP}V.~A.~Iskovskikh, Y.~G.~Prokhorov, \textit{Fano
varieties. Algebraic geometry, V}, 1--247, Encyclopaedia Math. Sci.,
47, Springer, Berlin, 1999.

\bibitem[KMM]{KMM} Y.~Kawamata, K. Matsuda, K.~Matsuki, \textit{Introduction
to the minimal model problem}, Algebraic geometry, Sendai, 1985, 283--360,
Adv. Stud. Pure Math., 10, North-Holland, Amsterdam, 1987. 

\bibitem[Ka1]{Ka1} Y.~Kawamata, \textit{Small contractions of four-dimensional
algebraic manifolds}, Math. Ann. 284 (1989), no. 4, 595--600. 

\bibitem[Ko]{Ko} J.~Koll\'ar, \textit{Flops}, Nagoya Math. J. 113
(1989), 15--36.

\bibitem[L]{L} A.~Langer, \textit{Fano 4-folds with scroll structure,}
Nagoya Math. J. 150 (1998), 135--176. 

\bibitem[Mo]{Mo} S.~Mori, \textit{Threefolds whose canonical bundles
are not numerically effective}, Ann. of Math. (2) 116 (1982), no.
1, 133--176. 

\bibitem[Mu1]{Mu1} S.~Mukai, \textit{Biregular classification of
Fano 3-folds and Fano manifolds of coindex 3}, Proc. Nat. Acad. Sci.
U.S.A. 86 (1989), no. 9, 3000--3002.

\bibitem[Mu2]{Mu2} S.~Mukai, \textit{New developments in the theory
of Fano threefolds: vector bundle method and moduli problems,} Sugaku
Expositions 15 (2002), no. 2, 125--150.

\bibitem[R]{R}M.~Reid, \textit{Graded rings and birational geometry},
in Proc. of algebraic geometry symposium (Kinosaki, Oct.~2000), 1--72. 

\bibitem[SW]{SW} M. Szurek and J.Wi\'sniewski, \textit{Fano bundles
over $\mP^{3}$ and $Q^{3}$}, Pacific J. Math. 141 (1990), no. 1,
197--208. 

\bibitem[Tak1]{Tak1}H.~Takagi, \textit{On classification of $\mathbb{Q}$-Fano
$3$-folds of Gorenstein index $2$. I, II, }Nagoya Math. J. 167 (2002),
117--155, 157 -- 216. 

\bibitem[Tak2]{Tak2}H.~Takagi, \textit{Classification of primary
$\mathbb{Q}$-Fano threefolds with anti-canonical Du Val $K3$ surfaces.}
I, J.Alg.Geom.\textbf{15} (2006), no. 1, 31--85. 

\bibitem[Tak3]{Tak3} H.~Takagi, \textit{Duality Related with Key
Varieties of $\mQ$-Fano $3$-folds. I, }in preparation.

\bibitem[Tak4]{Tak4} H.~Takagi, \textit{Key varieties for prime
$\mQ$-Fano threefolds of codimension five related with $\mP^{2}\times\mP^{2}$-
or $\mP^{1}\times\mP^{1}\times\mP^{1}$-fibration,} in preparation. 

\bibitem[Tan]{Tan} R.~Tange, \textit{On embeddings of certain spherical
homogeneous spaces in prime characteristic}, Transform. Groups 17
(2012), no. 3, 861--888.
\end{thebibliography}
\end{document}